\documentclass[10pt,a4paper]{amsart} 
\usepackage{verbatim}
\usepackage{amssymb}
\usepackage{amsbsy}
\usepackage{amscd}
\usepackage{amsmath}
\usepackage{amsthm}
\usepackage[mathscr]{eucal}
\usepackage{epic,eepic}
\usepackage{graphicx}

\newtheorem{theorem}{Theorem}[section]
\newtheorem{assumption}[theorem]{Assumption}
\newtheorem{lemma}[theorem]{Lemma}
\newtheorem{proposition}[theorem]{Proposition}
\newtheorem{conjecture}[theorem]{Conjecture}
\newtheorem{corollary}[theorem]{Corollary}

\theoremstyle{definition} 
\newtheorem{Picture}[theorem]{Picture} 
\newtheorem{remark}[theorem]{Remark} 
\newtheorem{example}[theorem]{Example} 
\newtheorem{definition}[theorem]{Definition}
\newcommand{\C}{\mathbb C}
\newcommand{\R}{\mathbb R}
\newcommand{\Z}{\mathbb Z}
\newcommand{\Q}{\mathbb Q}

\newcommand{\F}{\mathbb F}
\newcommand{\Flim}{\F_{\rm lim}}

\newcommand{\bR}{\boldsymbol{R}}
 
\newcommand{\Proj}{\mathbb P}

\newcommand{\U}{\mathbb U}

\newcommand{\age}{\operatorname{age}}

\newcommand{\End}{\operatorname{End}}

\newcommand{\Ker}{\operatorname{Ker}}

\newcommand{\Image}{\operatorname{Im}}

\newcommand{\id}{\operatorname{id}}

\newcommand{\seminf}{$\frac{\infty}{2}$} 
\newcommand{\Res}{\operatorname{Res}}
\newcommand{\rank}{\operatorname{rank}}

\newcommand{\unit}{\operatorname{\boldsymbol{1}}}

\newcommand{\Tr}{\operatorname{Tr}}
\newcommand{\Eff}{\operatorname{Eff}}

\newcommand{\ch}{\operatorname{ch}} 
\newcommand{\tch}{\widetilde{\operatorname{ch}}} 
 
\newcommand{\tTd}{\widetilde{\operatorname{Td}}}
\newcommand{\codim}{\operatorname{codim}} 
 
\newcommand{\pr}{\operatorname{pr}}

\newcommand{\Aut}{\operatorname{Aut}} 
\newcommand{\Gr}{\operatorname{Gr}} 
\newcommand{\inv}{\operatorname{inv}}
\newcommand{\Log}{\operatorname{Log}}

\newcommand{\Spec}{\operatorname{Spec}}

\newcommand{\bzero}{\boldsymbol{0}} 

\newcommand{\sfT}{\mathsf{T}}

\newcommand{\cA}{\mathcal{A}}
\newcommand{\cU}{\mathcal{U}}
\newcommand{\cO}{\mathcal{O}}

\newcommand{\cL}{\mathcal{L}}
\newcommand{\cX}{\mathcal{X}}
\newcommand{\cH}{\mathcal{H}}

\newcommand{\cM}{\mathcal{M}}
\newcommand{\tcM}{\widetilde{\cM}}

\newcommand{\cZ}{\mathcal{Z}}

\newcommand{\cV}{\mathcal{V}}

\newcommand{\Sol}{\mathcal{S}}

\newcommand{\hF}{\widehat{F}} 
\newcommand{\hGamma}{\widehat{\Gamma}}

\newcommand{\ov}{\overline}

\newcommand{\iu}{\mathtt{i}}

\newcommand{\fract}[1]{\langle #1 \rangle}

\def\pair#1#2{\langle #1,#2\rangle}

\def\parfrac#1#2{\frac{\partial{#1}}{\partial #2}}
\def\corr#1{\left\langle #1 \right\rangle}

\title[Ruan's conjecture and integral structures]
{Ruan's conjecture and integral structures in quantum cohomology}

\author[Hiroshi Iritani]{Hiroshi Iritani}

\address{}

\email{iritani@math.kyushu-u.ac.jp}
\email{h.iritani@imperial.ac.uk} 

\begin{document} 

\maketitle
\begin{abstract}
This is an expository article on 
the recent studies \cite{CIT,Coates-R,I:real,Coates} 
of Ruan's crepant 
resolution/flop conjecture \cite{Ruan:crepant1, Ruan:crepant2} 
and its possible relations to 
the $K$-theory integral structure 
\cite{I:real, KKP} in quantum cohomology. 
\end{abstract}

\section{Introduction}
The small quantum cohomology is a family $(H^*(X),\circ_\tau)$ 
of commutative ring structures on $H^*(X)$   
parametrized by $\tau\in H^{1,1}(X)$.  
The quantum product $\circ_\tau$ 
goes to the cup product in the \emph{large radius limit}: 
$-\Re \left(\int_C \tau\right) 
\to \infty$ for every effective curve $C\subset X$. 

Roughly speaking, Yongbin Ruan's conjecture 
says that, for a pair $(X_1,X_2)$ of 
birational varieties in some ``crepant" relationships 
(like flops or crepant resolutions), 
the small quantum cohomologies 
$(H^*(X_1),\circ_{\tau_1})$ and $(H^*(X_2),\circ_{\tau_2})$ 
are isomorphic under analytic continuation 
of the parameter $\tau$. 
\begin{figure}[htbp]
\begin{center}
\begin{picture}(300,100) 
\put(20,5){\line(1,0){230}}
\put(20,5){\line(1,3){30}} 
\put(250,5){\line(1,3){30}} 
\put(50,95){\line(1,0){230}}
\put(270,90){\makebox(0,0){$\cM$}}
\put(75,45){\rotatebox{20}{\shade\ellipse{60}{50}}}
\put(225,55){\rotatebox{10}{\shade\ellipse{60}{50}}}
\put(63,55){\makebox(0,0){$V_1$}}
\put(237,70){\makebox(0,0){$V_2$}} 
%\multiput(80,50)(2,0.2){80}{\line(1,0){1}}
%\multiput(80,38)(2,0.13){70}{\line(1,0){1}}
{\linethickness{0.2pt}
\qbezier(80,38)(110,30)(150,40)
\qbezier(150,40)(190,50)(220,48)
}
\put(150,40){\vector(4,1){0}}
%\put(150,42.5){\vector(1,0){0}}
\put(155,35){\makebox(0,0){$\gamma(t)$}} 
\put(150,53){\arc{10}{-3}{2.3}} 
\path(144,56)(145,53)(147,55)
\put(159,57){\makebox(0,0){$c$}}
%\put(183,80){\makebox(0,0){$D$}}
\put(75,43){\makebox(0,0){$\bullet$}}
\put(223.5,55){\makebox(0,0){$\bullet$}}
\put(82,52){\makebox(0,0){$\bzero_1$}}
\put(230,49){\makebox(0,0){$\bzero_2$}}

{
\thicklines 
%\qbezier(37,37)(150,58)(265,55) % central line
\qbezier(37,37)(93,48)(153,52)
\qbezier(156.5,52)(207,57)(265,55)
\qbezier(70,80)(77,50)(75,15)   % V1 
\qbezier(230,88)(220,55)(225,20) % V2
\qbezier(105,80)(140,90)(150,65)  %discriminant 
\qbezier(127,72)(150,90)(175,80)
}
\end{picture}
\end{center} 
\caption{K\"{a}hler moduli space $\cM$ containing 
cusp neighborhoods $V_i \subset H^{1,1}(X_i,\C)$, $i=1,2$.
The \emph{global quantum $D$-module} over $\cM$ 
develops singularities along thick lines.}
\label{fig:GKM}
\end{figure}
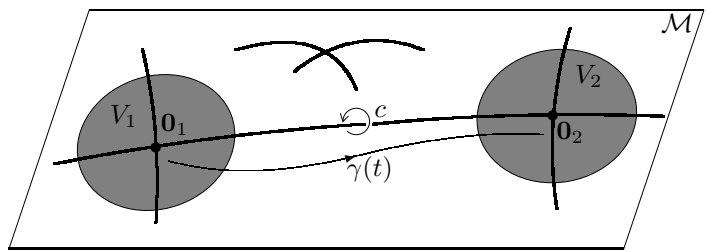 
The conjectural space 
where the quantum product $\circ_\tau$ is analytically 
continued is known as 
\emph{K\"{a}hler moduli space} $\cM$ 
(Figure \ref{fig:GKM}) in physics. 
% Roughly speaking, Yongbin Ruan's conjecture states that 
% (small) quantum cohomology of birational varieties 
% $X_1,X_2$ in some crepant relationships  
% (like flops or crepant resolutions)    
% have isomorphic quantum cohomology 
% under analytic continuation in quantum parameters. 
% This conjecture can be explained by the 
% picture of \emph{K\"{a}hler moduli space} $\cM$ 
% in physics. 
In our situation, this space $\cM$ 
has two limit points (cusps) $\bzero_1$, $\bzero_2$ 
corresponding to the large radius limit points of 
$X_1$ and $X_2$ respectively. 
A neighborhood $V_i$ of $\bzero_i$ 
is identified with an open subset of $H^{1,1}(X_i)$. 
A weak form of Ruan's conjecture asserts that 
\emph{there exists a family $(F,\circ_\tau)$ of 
commutative rings over $\cM$ 
such that its restriction to $V_i$ 
is isomorphic to the small quantum cohomology of $X_i$}.    
In particular, the cohomology rings $H^*(X_1)$, 
$H^*(X_2)$ are connected 
through quantum deformations. 

In a more precise picture, 
the family of rings should come 
from a $D$-module $(F,\nabla)$ 
(a meromorphic flat connection) over $\cM$ 
--- a \emph{global quantum $D$-module}. 
This $D$-module restricted to $V_i$ is 
identified with the quantum $D$-module 
given by the 
\emph{Dubrovin connection} (\ref{eq:Dub_conn}): 
\[
\nabla_\alpha = \parfrac{}{t^\alpha} + \frac{1}{z}\phi_\alpha \circ_\tau, 
\quad \text{where $z\in \C^*$ is a parameter}.  
\] 
The Dubrovin connection $z \nabla_\alpha$ recovers the 
quantum product $\phi_\alpha\circ_\tau$ in the limit $z\to 0$, 
but the $D$-module structure contains 
much more information than a family of rings. 
In fact, the global quantum $D$-module 
$(F,\nabla)$ together with additional data --- 
opposite subspace and dilaton shift ---  
yields a \emph{flat (or Frobenius) structure} on the 
(extended) K\"{a}hler moduli space\footnote{In the 
later formulation, $\cM$ will be extended so that 
it supports ``big" quantum cohomology 
whose deformation parameter $\tau$ lies in the 
total cohomology ring $H^*(X_i)$.}. 
Moreover, the local monodromy around each cusp 
determines a canonical choice 
of the opposite subspace 
and recovers the flat structure on $V_i$ 
coming from the vector space $H^{1,1}(X_i)$.  
Here, as the example in \cite{CIT} suggests,  
the flat structures from the different 
cusps $\bzero_1$ and $\bzero_2$ 
do not necessarily coincide.

% Also, the $D$-module gives rise to 
% many interesting monodromies 
% around its singularities. 
% Over the cusp neighborhood $V_i$, 
% $\nabla$ has logarithmic poles   
% along normal crossing divisors 
% intersecting at $\bzero_i$ and 
% the local monodromies around 
% these divisors are quasi-unipotent. 
% Some of the local monodromies  
% are common to $X_1$ and $X_2$ (\emph{e.g.} around $c$)
% and play an important role in our discussion. 

% In our picture, different 
% birational varieties arise from different 
% singularities of the global quantum $D$-module. 
% Also there may be some locus $D$ 
% which does not correspond to quantum cohomology of 
% the ordinary commutative geometry. 
% This picture will be formulated 
% more precisely in Section \ref{sec:GQDM}. 
%We will see that the $D$-module produces 
%a rich geometry (as Frobenius manifold 
%structures) on the extended 
%K\"{a}hler moduli space

In this article, we moreover postulate that 
the global quantum $D$-module is underlain 
by an integral local system. 
We also conjecture that, 
over $V_i$, the integral local system in question 
comes from the $K$-theory of $X_i$. 
This has the following physical explanation. 
Quantum cohomology is 
part of the A-model topological string theory.  
A chiral field in the A-model (\emph{i.e.} 
a section of the quantum $D$-module) 
should have a pairing with a B-type D-brane 
(\emph{i.e.} an object 
of the derived category $D^b_{\rm coh}(X_i)$) 
(see \emph{e.g.} \cite{Hor-Vaf}). 
This suggests that \emph{a vector bundle  
on $X_i$ should give a flat section of the quantum $D$-module}. 
In mirror symmetry, this is mirror to 
the fact that a holomorphic $n$-form 
has a pairing with a (real) Lagrangian 
$n$-cycle by integration. 
Based on mirror symmetry for toric orbifolds, 
the author \cite{I:real} proposed  
a formula (\ref{eq:Kgroupframing}) 
which assigns a flat section of the quantum $D$-module 
to an element of the $K$-group. 
Katzarkov-Kontsevich-Pantev \cite{KKP} 
also found a similar formula 
for a rational structure independently. 
The flat sections arising from the $K$-group
define an integral local system over $V_i$. 
Via the analytic continuation of 
$K$-theory flat sections along a path $\gamma(t)$ 
connecting $V_1$ and $V_2$ (see Figure \ref{fig:GKM}),  
we obtain an isomorphism of $K$-groups: 
\[
\U_{K,\gamma} \colon K(X_1) 
\overset{\cong}{\longrightarrow} K(X_2).  
\]
The isomorphism $\U_{K,\gamma}$ contains complete 
information about relationships 
between genus zero Gromov-Witten theories 
(quantum cohomology) of $X_1$ and $X_2$.  
We expect that $\U_{K,\gamma}$ is 
given by a certain Fourier-Mukai transformation. 
% Our proposal is a conjecture between conjectures: 
% A conjectural isomorphism in $K$-groups 
% gives a description of Ruan's conjecture. 

The paper is structured as follows. 
In Section \ref{sec:Kint}, we review 
orbifold quantum cohomology/$D$-module 
and introduce the $K$-theory integral structure on it. 
% We also describe the Galois actions 
% (or monodromy actions) on quantum $D$-modules.  
In Section \ref{subsec:GQDM}, 
we formulate a precise picture (Picture \ref{Pic:GQDM}) 
of the global quantum $D$-module sketched above. 
% In the rest of the section, 
% \emph{assuming} this picture, 
% we will see what follows from this:   
%\begin{itemize}
% \item \S \ref{subsec:Fmanifold}. An isomorphism of 
% quantum cohomology as \emph{$F$-manifolds}.  
%
% \item \S \ref{subsec:seminf}, 
% \S \ref{subsec:opp_Frob}. 
% A construction of Frobenius 
% manifold structures 
% from the \seminf VHS associated to the 
% global quantum $D$-module. 
%
% \item \S \ref{subsec:opp_cusp}.  
% A characterization of the Frobenius 
% manifold structure on the cusp neighborhood $V_i$ 
% (using the monodromy invariance 
% and the Deligne extension). 
%
% \item \S \ref{subsec:symp}. 
% Symplectic transformation between 
% the Givental spaces (connection 
% to the other formulations \cite{CIT,Coates-R}). 
%
% \item \S \ref{subsec:HL}. Hard Lefschetz condition which 
% guarantees that the Frobenius structures 
% arising from different cusps match. 
%
% \item \S \ref{subsec:intperiod}, \S \ref{subsec:locex}.  
% Integral periods on the K\"{a}hler 
% moduli space and the crepant resolution 
% conjecture for $\C^n/G$, $n=2,3$ 
% where $G\subset SL(n, \C)$ is a finite 
% subgroup. 
% \end{itemize} 
In Sections \ref{subsec:Fmanifold}--\ref{subsec:HL},  
we discuss what follows from the picture 
\emph{without} using integral structures. 
The main observation here is the fact 
that \emph{each cusp determines a (possibly different) 
Frobenius/flat structure on $\cM$}. 
%The Frobenius structure 
%arising from the cusp $\bzero_i$  
%is the one canonically associated to 
%the quantum cohomology of $X_i$. 
The Hard Lefschetz condition in Section 
\ref{subsec:HL} is a sufficient condition 
for the Frobenius structures from different cusps 
to match. 
These facts were found in \cite{CIT}, 
but the present article contains a complete proof 
of the characterization of Frobenius structures at cusps 
(Theorem \ref{thm:uniqueness_cusp}, announced in \cite{CIT}) 
and a \emph{generalized} Hard Lefschetz condition 
(Theorem \ref{thm:HL}). 
In Sections \ref{subsec:intperiod}, \ref{subsec:locex}, 
we use integral structures to 
study the crepant resolution conjecture 
for Calabi-Yau orbifolds and 
give an explicit prediction  
(Conjecture \ref{conj:localex}) 
for the change of co-ordinates 
in local examples. 
Readers who want to know a role of integral 
structures in Ruan's conjecture can safely skip 
Sections \ref{subsec:Fmanifold}--\ref{subsec:HL} 
and go directly to Sections \ref{subsec:intperiod}
or \ref{subsec:locex}. 

\vspace{5pt} 
\noindent
{\bf Acknowledgments} 
Thanks are due to Masa-Hiko Saito 
for the invitation to the conference 
``New developments in Algebraic Geometry,
Integrable Systems and Mirror symmetry" 
held at RIMS, January 2008. 
The author is grateful to  
Tom Coates, 
Alessio Corti,  
Hsian-Hua Tseng 
for helpful conversations 
during joint works \cite{CIT,CCIT:comp,CCIT:toric} 
and in many other occasions.  
He is also grateful to 
Samuel Boissiere 
for explaining the McKay correspondence 
for a finite subgroup of $SO(3)$ 
and to 
Jim Bryan, 
Martin Guest and 
Yongbin Ruan  
for helpful comments on 
the present article. 
He also thanks  
Claus Hertling,  
Shinobu Hosono,  
Tony Pantev,   
Shunsuke Takagi 
for helpful discussions. 
This research is supported by 
EPSRC(EP/E022162/1).  

\tableofcontents

\section{$K$-theory integral structure in quantum cohomology} 
\label{sec:Kint}
In this section, we review the 
orbifold quantum cohomology for 
smooth Deligne-Mumford stacks 
and introduce the $K$-theory integral structure on it. 
Assuming the convergence of structure constants, 
quantum cohomology defines a flat connection, 
called \emph{Dubrovin connection}, on some cohomology 
bundle over a neighborhood of 
the ``large radius limit point". 
This is called \emph{quantum $D$-module}. 
We will see that the $K$-group defines 
an integral lattice in the space of (multi-valued) 
flat sections of the quantum $D$-module. 
The key definition will be given 
in Definition \ref{def:Ktheoryintstr}. 
The true origin of this integral structure 
is yet to be known, but it has a number of 
good properties: 
\begin{itemize}
\item This is invariant under every local monodromy 
around the large radius limit point. 

\item The pairing on quantum cohomology 
is translated into the Mukai pairing on the $K$-group. 

\item This gives a real structure 
which is \emph{pure and polarized} in a neighborhood of 
the large radius limit point \cite{I:real}. 
In particular, we have \emph{$tt^*$-geometry} \cite{CV,Her} 
on quantum cohomology.

\item This looks compatible with many 
computations done in the context of mirror symmetry 
\cite{Hosono, Bor-Hor}.  
Especially this matches with the integral structure on 
the Landau-Ginzburg mirror in the case of 
toric orbifolds \cite{I:real}. 

\item Thus in toric case, this integral structure 
is compatible also with the \emph{Stokes structure}. 
\end{itemize} 
In this article, we will not explain 
the last three items. 
See \cite{I:real, HS, KKP} for the properties 
``pure and polarized" or 
``compatibility with Stokes structure".

\subsection{Orbifold quantum cohomology} 
We start from the notation on orbifolds. 
Let $\cX$ be a smooth Deligne-Mumford stack 
with projective coarse moduli space $X$.  
Let $I\cX$ be the inertia stack of $\cX$. 
A point on $I\cX$ is given by a pair $(x,g)$ 
of a point $x\in\cX$ and an element $g$ 
of the automorphism group (local group) 
$\Aut_\cX(x)$ at $x$. 
The element $g\in \Aut_\cX(x)$ 
is also called a stabilizer. 
Let 
\[
I\cX = \bigsqcup_{v\in \sfT} \cX_v = 
\cX_0 \sqcup \bigsqcup_{v\in \sfT'} \cX_v
\] 
be the decomposition of $I\cX$ into connected components. 
Here $\sfT$ is a finite set parametrizing 
connected components of $I\cX$. 
$\sfT$ contains a distinguished element 
$0\in \sfT$ which corresponds to 
the trivial stabilizer $g=1$ and 
we set $\sfT = \{0\} \cup \sfT'$. 
Then $\cX_0$ is isomorphic to $\cX$. 
At each point $(x,g)$ in $I\cX$, we can define a 
rational number $\iota_{(x,g)}$ called \emph{age}. 
The element $g$ of the automorphism group 
acts on the tangent space $T_x\cX$ and 
decomposes it into eigenspaces: 
\[
T_x\cX = \bigoplus_{0\le f<1} (T_x\cX)_f 
\]
where $g$ acts on $(T_x\cX)_f$ by $\exp(2\pi\iu f)$. 
The age $\iota_{(x,g)}$ is defined to be 
\[
\iota_{(x,g)} = \sum_{0\le f<1} f \dim_\C (T_x\cX)_f. 
\]
The age $\iota_{(x,g)}$ is constant along the connected 
component $\cX_v$ of $I\cX$, so we denote 
by $\iota_v$ the age $\iota_{(x,g)}$ 
at any point $(x,g)$ in $\cX_v$. 
The \emph{orbifold or Chen-Ruan cohomology group}  
$H_{\rm CR}^*(\cX)$ is a $\Q$-graded vector space defined by 
\[
H_{\rm CR}^p(\cX) = \bigoplus_{v\in\sfT} 
H^{p-2\iota_v}(\cX_v,\C), \quad p\in \Q. 
\]
This is the same as $H^*(I\cX,\C)$ as a vector space, 
but the grading is shifted by the age. 
In this paper, 
we only consider the \emph{even parity part} 
of $H_{\rm CR}^*(\cX)$, \emph{i.e.} 
the summands satisfying $p-2\iota_v \equiv 0 \mod 2$ 
in the above decomposition. 
Unless otherwise stated, we denote 
by $H_{\rm CR}^*(\cX)$ the even parity part. 
The inertia stack has an involution 
$\inv\colon I\cX\to I\cX$ which sends 
$(x,g)$ to $(x,g^{-1})$. 
This induces an involution $\inv\colon \sfT \to \sfT$ 
on the index set $\sfT$ and 
$\inv^*\colon H_{\rm CR}^*(\cX) 
\to H_{\rm CR}^*(\cX)$ on the cohomology. 
The \emph{orbifold Poincar\'{e} pairing} on 
$H_{\rm CR}^*(\cX)$ is defined by 
\[
(\alpha,\beta)_{\rm orb} = \int_{I\cX} \alpha \cup \inv^*(\beta) 
= \sum_{v\in \sfT} \int_{\cX_v} \alpha_v \cup \beta_{\inv(v)},   
\]
where $\alpha_v,\beta_v$ are 
the $v$-components of $\alpha,\beta$. 
This pairing is symmetric, non-degenerate 
and of degree $-2\dim_\C\cX$. 

Gromov-Witten theory for manifolds 
has been extended to the class of 
symplectic orbifolds or smooth Deligne-Mumford stacks. 
This was done by Chen-Ruan \cite{CR} in the symplectic 
category and by Abramovich-Graber-Vistoli \cite{AGV} 
in the algebraic category.  
The formal properties of 
the genus zero Gromov-Witten theory 
hold in orbifold theory as well:  
the genus zero orbifold Gromov-Witten theory 
defines a cohomological field theory 
(see \emph{e.g.} \cite{Manin})  
on the metric vector space 
$(H_{\rm CR}^*(\cX),(\cdot,\cdot)_{\rm orb})$. 
In particular, we have the following 
correlation functions (Gromov-Witten invariants): 
\begin{equation}
\label{eq:primaryGW} 
\corr{\cdot,\dots,\cdot}_{0,m,d} \colon 
(H_{\rm CR}^*(\cX))^{\otimes m} \to \C 
\end{equation} 
defined for $m\ge 0$ and $d\in H_2(X,\Z)$. 
This is zero when 
$d$ is not in the semigroup $\Eff_\cX\subset H_2(X,\Z)$
generated by classes of effective curves 
or $m\le 2$ and $d=0$. 
Also these correlation functions satisfy 
the so-called \emph{WDVV equation} or the 
\emph{splitting axiom} 
(see \emph{e.g.} \cite[Theorem 6.4.3]{AGV}).   
The genus zero Gromov-Witten invariant 
is homogeneous with respect 
to the grading of $H_{\rm CR}^*(\cX)$. 
More precisely, $
\corr{\alpha_1,\dots,\alpha_m}_{0,m,d}= 0$  
unless $p_1+\dots+ p_m = 2 (\dim_\C \cX  
+ \pair{c_1(\cX)}{d}+m-3)$, 
where $\alpha_i\in H^{p_i}_{\rm CR}(\cX)$.  

The genus zero Gromov-Witten invariants  
define a quantum product $\bullet_\tau$ 
on $H_{\rm CR}^*(\cX)$ 
parametrized by $\tau\in H_{\rm CR}^*(\cX)$: 
\begin{equation}
\label{eq:quantumproduct}
(\alpha\bullet_\tau\beta,\gamma)_{\rm orb} 
= \sum_{d\in \Eff_\cX,l\ge 0}
\frac{1}{m!} \corr{\alpha,\beta,\gamma,
\overbrace{\tau,\dots,\tau}^{\text{$m$ times}}}_{0,m+3,d} Q^d.  
\end{equation} 
Here $Q^d$ denotes the element of the group ring 
$\C[\Eff_\cX]$ corresponding to 
$d\in \Eff_\cX\subset H_2(X,\Z)$. 
The right-hand side belongs to 
$\C[\![\tau]\!][\![\Eff_\cX]\!]$ 
(a certain completion\footnote{For example, 
one completion is given by the additive 
valuation $v$ on $\C[\Eff_\cX]$ defined by 
$v(Q^d) = \int_d \omega$, 
where $\omega$ is a K\"{a}hler form on $\cX$. }
of $\C[\![\tau]\!]\otimes \C[\Eff_\cX]$)  
and defines the element $\alpha\bullet_\tau\beta$ 
in $H_{\rm CR}^*(\cX)\otimes \C[\![\tau]\!][\![\Eff_\cX]\!]$ 
because the orbifold Poincar\'{e} pairing is non-degenerate.  
By extending $\bullet_\tau$ 
as a $\C[\![\tau]\!][\![\Eff_\cX]\!]$-bilinear map, 
we have an associative commutative ring 
$(H_{\rm CR}^*(\cX)\otimes \C[\![\tau]\!][\![\Eff_\cX]\!],\bullet_\tau)$. 
Here the associativity of the product $\bullet_\tau$ 
follows from the WDVV equation. 
This is the \emph{orbifold quantum cohomology} of $\cX$. 

Using the \emph{Divisor equation} (see \emph{e.g.} 
\cite[Theorem 8.3.1]{AGV}), we can write 
\[
(\alpha\bullet_\tau\beta,\gamma)_{\rm orb} 
= \sum_{d\in \Eff_\cX, m\ge 0}
\frac{1}{m!} \corr{\alpha,\beta,\gamma,
\overbrace{\tau',\dots,\tau'}^{\text{$m$ times}}}_{0,m+3,d} 
e^{\pair{\tau_{0,2}}{d}} Q^d, 
\] 
where we put 
\begin{equation}
\label{eq:decomp_tau}
\tau=\tau_{0,2} + \tau', \quad 
\tau_{0,2}\in H^2(\cX_0), \quad 
\tau'\in \bigoplus_{p\neq 2} H^p(\cX_0) \oplus 
\bigoplus_{v\in\sfT'}H^*(\cX_v). 
\end{equation} 
This shows that the parameters $\tau$ and $Q$ 
in the product $\bullet_\tau$ are redundant. 
In fact $\bullet_\tau$ depends only on 
$\tau'$ and $e^{\tau_{0,2}}Q$.   
We put 
\[
\circ_\tau:=\bullet_\tau|_{Q=1}. 
\]
The product $\circ_\tau$ is 
a formal power series in $\tau'$ and 
a formal Fourier series in $\tau_{0,2}$. 
It is clear from the formula that 
$\circ_\tau$ recovers $\bullet_\tau$. 
In what follows, we will study 
$\circ_\tau$ instead of $\bullet_\tau$ 
and assume that 
\emph{the product $\circ_\tau$ is convergent 
in some open set $U$ of $H_{\rm CR}^*(\cX)$}.  
\begin{assumption}
\label{as:convergence}
The orbifold quantum product $\circ_\tau$ is convergent 
on a simply-connected open set $U$ 
containing the following set 
\[
\{\tau\in H^*_{\rm CR}(\cX)\;;\; 
\Re(\pair{d}{\tau_{0,2}})<-M, \forall d\in \Eff_\cX\setminus\{0\}, 
\ \|\tau'\|\le e^{-M}. \} 
\]
where $\tau=\tau_{0,2}+\tau'$ is the decomposition 
in (\ref{eq:decomp_tau}), $M>0$ is sufficiently big  
and $\|\cdot\|$ is a suitable norm 
on $H_{\rm CR}^*(\cX)$. 
\end{assumption} 
\begin{remark}
Working over a certain formal 
power series ring, we could discuss 
the $K$-theory integral structure 
without this assumption. 
However, when considering Ruan's conjecture later, 
we cannot avoid the convergence problem 
of quantum cohomology. 
\end{remark} 

The open set $U$ above is considered to be a 
neighborhood of the ``\emph{large radius limit point}" 
which is the limit point of the sequence 
\begin{equation}
\label{eq:largeradius}
\tau=\tau_{0,2}+\tau' \colon 
\quad \Re(\pair{d}{\tau_{0,2}})\to -\infty, \quad 
\tau' \to 0.  
\end{equation} 
(This notion will be made more precise later.) 
In this limit, the orbifold quantum product $\circ_\tau$ 
goes to the \emph{Chen-Ruan orbifold cup product} 
$\cup_{\rm CR}$. 
This product $\cup_{\rm CR}$ is the same as the cup 
product when $\cX$ is a manifold, but in the orbifold case, 
this is different from the cup product on $I\cX$. 

\subsection{Quantum $D$-modules with Galois actions} 
\label{subsec:QDM} 
Let $\{\phi_i\}$ be a homogeneous 
$\C$-basis of $H_{\rm CR}^*(\cX)$ 
and $\{t^i\}$ be the linear 
co-ordinate system on $H_{\rm CR}^*(\cX)$ dual to 
the basis $\{\phi_i\}$. 
Denote by $\tau = \sum_{i=1}^N t^i \phi_i$  
a general point on $H_{\rm CR}^*(\cX)$. 
The quantum $D$-module is a meromorphic flat 
connection on the trivial 
$H_{\rm CR}^*(\cX)$-bundle over $U\times \C$. 
Denote by $(\tau,z)$ a general point on the base space 
$U\times \C$. 
Let $(-)\colon U\times \C\to U\times \C$ be the 
map sending $(\tau,z)$ to $(\tau,-z)$.

\begin{definition}
\label{def:QDM} 
\emph{Quantum $D$-module} 
$QDM(\cX)=(F,\nabla,(\cdot,\cdot)_F)$ 
is the trivial holomorphic vector bundle 
$F:=H_{\rm CR}^*(\cX)\times (U\times \C) \to (U\times \C)$ 
endowed with the meromorphic flat connection 
$\nabla$: 
\begin{align}
\label{eq:Dub_conn}
\begin{split}
&\nabla_i = \nabla_{\parfrac{}{t^i}} = 
\parfrac{}{t^i} + \frac{1}{z} \phi_i \circ_\tau, \\  
&\nabla_{z\partial_z} = 
z\parfrac{}{z} - \frac{1}{z} E\circ_\tau + \mu, 
\end{split}
\end{align} 
and the $\nabla$-flat pairing 
\[
(\cdot,\cdot)_F \colon (-)^*\cO(F) 
\otimes \cO(F) \to \cO_{U\times \C}
\]
induced from the orbifold Poincar\'{e} pairing 
$F_{(\tau,-z)}\times F_{(\tau,z)} 
= H_{\rm CR}^*(\cX)\times H_{\rm CR}^*(\cX) \to \C$. 
Here $E$ is the \emph{Euler vector field} on $U$ 
given by 
\begin{equation}
\label{eq:Euler}
E := c_1(T\cX) + 
\sum_{i=1}^N (1 - \frac{1}{2}\deg \phi_i) t^i \phi_i  
\end{equation} 
and $\mu\in \End(H_{\rm CR}^*(\cX))$ 
is the \emph{Hodge grading operator} defined by 
\begin{equation}
\label{eq:Hodgegrading}
\mu(\phi_i) := (\frac{1}{2} \deg\phi_i -\frac{n}{2})\phi_i, 
\quad n = \dim_\C \cX.    
\end{equation} 
The flat connection $\nabla$ is called \emph{Dubrovin connection} 
or \emph{the first structure connection}. 
Note that $\nabla_i$ has a pole of order 1 along $z=0$ 
and $\nabla_{\partial_z}$ has a pole of order 2 along $z=0$.  
The flatness of $\nabla$ follows from 
the WDVV equations and the homogeneity 
of Gromov-Witten invariants. 
\end{definition} 

\begin{remark} 
By \emph{$D$-module} one means a module over 
the ring of differential operators. 
In our case, 
the ring $\cO_{\cM\times \C^*}\langle 
\partial_{t^i}, z\partial_z\rangle$ 
of differential operators on $\cM\times \C^*$ 
acts on the space of sections of $F$ 
via the flat connection: 
$\partial_{t^i} \mapsto \nabla_i$, 
$z\partial_z \mapsto \nabla_{z\partial_z}$. 
This explains the name ``quantum $D$-module". 
\end{remark} 

The quantum $D$-module admits certain discrete 
symmetries (Galois actions). 
Firstly, since $\circ_\tau$ depends only on 
$e^{\tau_{0,2}}$ and $\tau'$, 
it is clear that $\circ_\tau$ is invariant 
under the following translation:  
\[
\tau_{0,2}\mapsto \tau_{0,2} - 2\pi\iu \xi, \quad 
\xi \in H^2(X,\Z).  
\]
This is a consequence of the Divisor equations 
and is familiar in ordinary Gromov-Witten theory. 
Interestingly, we have a finer symmetry for orbifold theory. 
Let $H^2(\cX,\Z)$ be the sheaf cohomology of the constant 
sheaf $\Z$ on the stack $\cX$ (not on the coarse 
moduli space $X$). This group is identified with 
the set of isomorphism classes of 
topological orbifold line bundles on $\cX$. 
Then $H^2(X,\Z)$ is identified with 
the subset of $H^2(\cX,\Z)$ 
consisting of line bundles which 
are pulled back from the coarse moduli space $X$. 
For $\xi\in H^2(\cX,\Z)$, let $L_\xi$ be the 
corresponding topological orbifold line bundle on 
$\cX$ and $\xi_0 :=c_1(L_\xi)\in H^2(X,\Q)$ be 
the first Chern class. 
For $v\in \sfT$, define $0\le f_v(\xi)<1$ 
to be the rational number such that  
the stabilizer $g$ at $(x,g)\in \cX_v$ 
acts on the fiber $L_{\xi,x}$ 
by $\exp(2\pi\iu f_v(\xi))$.  
\begin{lemma}[{\cite[Proposition 3.1]{I:real}}]
\label{lem:Galois}  
The flat connection $\nabla$ and the 
pairing $(\cdot,\cdot)_F$ of the quantum $D$-module 
$QDM(\cX)=(F,\nabla,(\cdot,\cdot)_F)$ 
is invariant under the following map 
given for $\xi\in H^2(\cX,\Z)$: 
\begin{align*} 
& H^*_{\rm CR}(\cX)\times (U\times\C) \to 
H^*_{\rm CR}(\cX) \times (U\times \C)  \\ 
&(\phi, \tau, z) \longmapsto 
(dG(\xi)(\phi), G(\xi)(\tau), z).  
\end{align*} 
Here $G(\xi), dG(\xi) \colon H_{\rm CR}^*(\cX) \to H_{\rm CR}^*(\cX)$ 
is defined by 
\begin{align*}
G(\xi)(\tau_0\oplus \bigoplus_{v\in \sfT'} \tau_v) 
&= (\tau_0 - 2\pi\iu \xi_0)\oplus \bigoplus_{v\in \sfT'}  
e^{2\pi\iu f_v(\xi)} \tau_v \\ 
dG(\xi)(\tau_0\oplus \bigoplus_{v\in \sfT'} \tau_v)
&=\tau_0 \oplus \bigoplus_{v\in \sfT'}  
e^{2\pi\iu f_v(\xi)} \tau_v  
\end{align*} 
where we used the decomposition 
$H^*_{\rm CR}(\cX) = H^*(\cX_0)\oplus 
\bigoplus_{v\in\sfT'} H^*(\cX_v)$ and 
$\tau_v\in H^*(\cX_v)$. 
(Here we implicitly assume that 
$U$ is invariant under the map $G(\xi)$, but 
we can assume this without loss of generality). 
\end{lemma}

When $\xi\in H^2(X,\Z)$, the above symmetry 
is the same as the aforementioned one. 
Note that the new symmetry can act non-trivially 
on the fiber of the quantum $D$-module. 
The quantum $D$-module descends 
to a flat connection on 
$F/H^2(\cX,\Z) \to (U/H^2(\cX,\Z))\times \C$. 
We call this flat connection on the quotient space 
also the \emph{quantum $D$-module}. 
In view of this, we can refer to the symmetries in 
Lemma \ref{lem:Galois} as \emph{Galois actions} 
or \emph{local monodromies at the large radius limit}. 

We can construct a partial 
compactification $\ov{V}$ 
of the quotient $V=U/H^2(\cX,\Z)$ such that $\ov{V}$ 
contains the large radius limit point 
and that the quantum $D$-module on $V$ 
extends to a $D$-module on $\ov{V}$ 
with a logarithmic singularity 
along the (\'{e}tale locally) 
normal crossing divisor $\ov{V}\setminus V$. 
Choose a $\Z$-basis $p_1,\dots,p_r$ 
of $H^2(X,\Z)/{\rm torsion}$ such that $p_a$ intersects every 
effective curve class $d\in \Eff_\cX$ 
non-negatively (\emph{i.e.} $p_a$ is nef). 
Then we have the embedding 
\[
U/H^2(X,\Z) \hookrightarrow \C^r \times W, \quad 
\left[\sum_{a=1}^r t^a p_a + \tau'\right] 
\mapsto (e^{t_1},\dots,e^{t_r}, \tau')
\] 
where $W = \bigoplus_{p\neq 2} H^p(\cX_0) \oplus 
\bigoplus_{v\in \sfT'} H^*(\cX_v)$ and $\tau'\in W$. 
By Assumption \ref{as:convergence} and the 
choice of $p_a$, 
the image of this embedding contains 
the open set 
$((\C^*)^r \times W) \cap \Delta_M$ 
for a sufficiently big $M>0$, where 
\[
\Delta_M = 
\{(q^1,\dots,q^r,\tau')\in \C^r\times W\;;\; 
|q^a|<e^{-M}, \|\tau'\|<e^{-M} \}. 
\]
We set $\ov{U/H^2(X,\Z)} 
:= (U/H^2(X,\Z))\cup \Delta_M \subset 
\C^r \times W$. 
For $\tau_{0,2} = \sum_{a=1}^r t^a p_a$, 
we have $e^{\pair{\tau_{0,2}}{d}} = 
(e^{t^1})^{\pair{p_1}{d}}
\cdots (e^{t^r})^{\pair{p_r}{d}}$. 
Therefore by the formula (\ref{eq:quantumproduct}), 
since $p_a$ is nef,  
the quantum product $\circ_\tau$ 
on $U/H^2(X,\Z)$ extends to $\ov{U/H^2(X,\Z)}$. 
The Dubrovin connection on $\Delta_M$ 
in the direction of $q^a=e^{t^a}$ can be written as 
\[
\nabla_{\parfrac{}{t^a}} = 
q^a\parfrac{}{q^a} + \frac{1}{z}p_a \circ_\tau. 
\] 
Hence it has a logarithmic pole along $q^1\cdots q^r =0$. 
We can now define $\ov{V}$ as the quotient space 
(or stack): 
\[
\ov{V} := \ov{U/H^2(X,\Z)}/(H^2(\cX,\Z)/H^2(X,\Z)). 
\]
This contains both 
$U/H^2(\cX,\Z)$ and 
the large radius limit point $q=\tau'=0$. 

\begin{remark}
The partial compactification $\ov{V}$ 
depends on the choice of a nef basis $p_a$. 
The most canonical choice of a partially compactified 
base space might be the possibly singular stack 
$\ov{V}_{\rm can} =  
(\Spec\C[\Eff_\cX]\times W)/(H^2(\cX,\Z)/H^2(X,\Z))$.   
Then we always have a map $\ov{V}\to \ov{V}_{\rm can}$. 
\end{remark} 

\begin{remark} 
Due to the new discrete symmetries, 
the large radius limit point in $\ov{V}$ 
can have an orbifold singularity 
when $\cX$ is an orbifold. 
Also, the quantum $D$-module $F/H^2(\cX,\Z)$ 
on the quotient space may not be trivialized 
in the standard way. 
In other words, an element of $H_{\rm CR}^*(\cX)$ 
gives a possibly multi-valued section of 
$F/H^2(\cX,\Z)$. 
\end{remark} 

\subsection{Fundamental solution $L(\tau,z)$ 
and the space $\Sol(\cX)$ of flat sections} 

We introduce a \emph{fundamental solution} 
for $\nabla$-flat sections of the quantum 
$D$-module $(F,\nabla)$. 
Orbifold Gromov-Witten theory also has 
(gravitational) descendant invariants (as opposed to the primary 
invariants (\ref{eq:primaryGW})) of 
the form 
\[
\corr{\alpha_1\psi_1^{k_1},\dots,\alpha_m \psi_m^{k_m}}_{0,m,d} 
\]
where $\alpha_i\in H_{\rm CR}^*(\cX)$, $d\in \Eff_\cX$ 
and $k_i$ is a non-negative integer. 
The symbol $\psi_i$ represents the first Chern class 
of the line bundle on the moduli space of stable maps 
formed by the cotangent lines  
at the $i$-th marked point of the coarse domain  curve. 
As is well-known in manifold Gromov-Witten theory 
(see \emph{e.g.}\ \cite[Proposition 2]{Pand}), 
we can write the fundamental solution to 
the equation $\nabla s =0$ by 
using descendant invariants. 
Let $\pr \colon I\cX \to \cX$ be the natural projection. 
For $\tau_0\in H^*(\cX_0)$, we define 
the action of $\tau_0$ on $H_{\rm CR}^*(\cX)$ 
as 
\[
\tau_0 \cdot \alpha = \pr^*(\tau_0) \cup \alpha  
\]
where the right-hand side is the cup product on 
$H^*(I\cX)$. (This is known to be the same as 
the orbifold cup product $\tau_0\cup_{\rm CR}\alpha$). 
Let $\{\phi_k\}_{k=1}^N$ and $\{\phi^k\}_{k=1}^N$ 
be bases of $H_{\rm CR}^*(\cX)$ dual with respect to 
the orbifold Poincar\'{e} pairing, \emph{i.e.}  
$(\phi_i,\phi^j)_{\rm orb}=\delta_i^j$.

\begin{proposition}[See \emph{e.g.}  
{\cite[Proposition 3.3]{I:real}}] 
\label{prop:fundamentalsol} 
Let $L(\tau,z)$ be the following $\End(H_{\rm CR}^*(\cX))$-valued 
function on $U\times \C^*$: 
\begin{equation}
\label{eq:fundamentalsol}
L(\tau,z) \phi = 
e^{-\tau_{0,2}/z} \phi  - 
\sum_{\substack{(d,m)\neq (0,0), \\ d\in \Eff_\cX, 1\le k\le N} } 
\frac{e^{\pair{\tau_{0,2}}{d}}}{m!} \phi_k
\corr{\phi^k, \tau',\dots,\tau',
\frac{e^{-\tau_{0,2}/z}\phi}{z+\psi_{m+2}}}_{0,m+2,d}, 
\end{equation} 
where $\tau =\tau_{0,2}+\tau'$ is the 
decomposition in (\ref{eq:decomp_tau}) 
and $1/(z+\psi_{m+2})$ in the correlator 
should be expanded in the $z^{-1}$-series 
$\sum_{k\ge 0} (-1)^k z^{-k-1}\psi_{m+2}^k$. 
Set $\rho := c_1(\cX)\in H^2(\cX_0)$ and 
\[
z^{-\mu}z^\rho := \exp(-\mu\log z) \exp(\rho\log z), \quad 
\text{$\mu$ is given in (\ref{eq:Hodgegrading}).} 
\]
Then we have 
\begin{align*} 
&\nabla_i (L(\tau,z) z^{-\mu}z^\rho \phi) = 0, 
\quad
\nabla_{z\partial_z} (L(\tau,z) z^{-\mu}z^\rho \phi) =0, \\ 
& (L(\tau,-z)\phi_i, L(\tau,z)\phi_j)_{\rm orb} 
=(\phi_i,\phi_j)_{\rm orb}. 
\end{align*} 
In particular, $s_i(\tau,z) =L(\tau,z) z^{-\mu}z^\rho \phi_i$, 
$1\le i\le N$,   
form a basis of multi-valued $\nabla$-flat sections of $F$ 
satisfying the asymptotic initial condition 
at the large radius limit (\ref{eq:largeradius}):  
\[
s_i(\tau,z) \sim z^{-\mu} z^{\rho} e^{-\tau_{0,2}}\phi_i. 
\]
\end{proposition} 

\begin{remark}
The convergence of the fundamental solution $L(\tau,z)$ 
is not a priori clear. 
From the Assumption \ref{as:convergence}, 
we know that $L(\tau,z)$ also converges 
on $U\times \C^*$ because this is 
a solution to the 
linear partial differential equations $\nabla s=0$. 
\end{remark} 

\begin{definition} 
\label{def:spaceofflatsections} 
Define $\Sol(\cX)$ to be the space of multi-valued 
$\nabla$-flat sections of the quantum $D$-module 
$QDM(\cX)=(F,\nabla,(\cdot,\cdot)_F)$: 
\[
\Sol(\cX) := \{s(\tau,z) \in 
\Gamma(U\times \widetilde{\C^*}, \cO(F))
\; ;\; \nabla s =0\}.  
\]
This is a $\C$-vector space with $\dim_\C\Sol(\cX) 
= \dim_\C H_{\rm CR}^*(\cX)$. 
$\Sol(\cX)$ is endowed with the 
pairing $(\cdot,\cdot)_\Sol$: 
\[
(s_1,s_2)_{\Sol} := (s_1(\tau,e^{\pi\iu} z), s_2(\tau,z))_{\rm orb} 
\in \C,   
\]
where $s_1(\tau,e^{\pi\iu} z)$ denotes the parallel translate 
of $s_1(\tau,z)$ along the counterclockwise 
path $[0,1]\ni \theta\mapsto e^{\iu\pi\theta} z$. 
Because $s_1,s_2$ are flat sections, 
the right-hand side is a complex number 
independent of $(\tau,z)$. 
$\Sol(\cX)$ is also equipped with 
the automorphism $G^\Sol(\xi)$ for $\xi\in H^2(\cX,\Z)$ 
induced from the Galois action in 
Lemma \ref{lem:Galois}: 
\[
G^\Sol(\xi) \colon \Sol(\cX) \to \Sol(\cX), \quad 
s(\tau,z) \mapsto dG(\xi) (s(G(\xi)^{-1}\tau,z)).  
\] 
\end{definition} 

In general, $(\cdot,\cdot)_\Sol$ is neither symmetric 
nor anti-symmetric. When $\cX$ is Calabi-Yau, \emph{i.e.} 
$\rho=c_1(\cX)=0$, $(\cdot,\cdot)_\Sol$ is 
symmetric when $n=\dim_\C\cX$ is even 
and is anti-symmetric when $n$ is odd. 

The fundamental solution in Proposition 
\ref{prop:fundamentalsol} gives the 
\emph{cohomology framing} 
$\cZ_{\rm coh}$ of $\Sol(\cX)$: 
\begin{equation}
\label{eq:cohframing}
\cZ_{\rm coh} \colon H_{\rm CR}^*(\cX) 
\overset{\cong}{\longrightarrow} \Sol(\cX), 
\quad \phi \mapsto L(\tau,z) z^{-\mu}z^\rho \phi. 
\end{equation} 
In terms of this cohomology framing $\cZ_{\rm coh}$, 
it is easy to check that the pairing and 
Galois actions on $\Sol(\cX)$ can be 
written as follows: 
\begin{align}
\label{eq:cohframing_Galois} 
\begin{split}
&(\cZ_{\rm coh}(\alpha),\cZ_{\rm coh}(\beta))_{\Sol} = 
(e^{\pi\iu \rho} \alpha, e^{\pi\iu\mu} \beta)_{\rm orb} \\ 
&G^{\Sol}(\xi) (\cZ_{\rm coh}(\alpha)) 
= \cZ_{\rm coh}(
(\bigoplus_{v\in \sfT} e^{-2\pi\iu \xi_0} 
e^{2\pi\iu f_v(\xi)})\alpha)  
\end{split}
\end{align} 
where we used the decomposition $H_{\rm CR}^*(\cX)
= \bigoplus_{v\in\sfT} H^*(\cX_v)$ in the second line. 
(See the paragraph before Lemma \ref{lem:Galois} 
for $\xi_0\in H^2(\cX_0)$ and $f_v(\xi)\in [0,1)$.)

\subsection{$K$-theory integral lattice of flat sections}

We will introduce an integral lattice 
in the space $\Sol(\cX)$ of flat sections 
using the $K$-group and the characteristic class 
called $\hGamma$-class. 
Let $K(\cX)$ be the Grothendieck group of topological 
orbifold vector bundles over $\cX$ 
(see \emph{e.g.} \cite{AR} for orbifold vector bundles 
and orbifold $K$-theory). 
For simplicity, we assume that 
$\cX$ is isomorphic to a quotient orbifold 
$[M/G]$ as a topological orbifold, 
where $M$ is a compact manifold and 
$G$ is a compact Lie group acting on $M$ 
with at most finite stabilizers. 
Under this assumption, $K(\cX)$ is isomorphic to the 
$G$-equivariant $K$-theory $K_G^0(M)$ and 
is a finitely generated abelian group \cite{AR}. 
For an orbifold vector bundle $V$ on $I\cX$ 
and a component $\cX_v$ of $I\cX$, 
we denote the eigenbundle decomposition 
of $V|_{\cX_v}$ with respect to 
the stabilizer action as follows: 
\[
V|_{\cX_v} = \bigoplus_{0\le f<1} V_{v,f},  
\]
where the stabilizer of $\cX_v$ acts on 
$V_{v,f}$ by $\exp(2\pi\iu f)$. 
The Chern character $\tch\colon K(\cX) \to H^*(I\cX)$ 
for orbifold vector bundles is defined as follows: 
\[
\tch(V) := \bigoplus_{v\in \sfT} 
\sum_{0\le f<1} e^{2\pi\iu f} \ch((\pr^*V)_{v,f}),   
\]
where $\pr\colon I\cX \to \cX$ is the natural projection. 
For an orbifold vector bundle $V$ on $\cX$, 
let $\delta_{v,f,i}$, $i=1,\dots,l_{v,f}$ be 
the Chern roots of the vector bundle 
$(\pr^*V)_{v,f}$ on $\cX_v$
(where $l_{v,f}=\rank (\pr^*V)_{v,f}$). 
The Todd class $\tTd(V)$ is defined by 
\[
\tTd(V) := \bigoplus_{v\in\sfT} 
\prod_{0<f<1, 1\le i\le l_{v,f}} 
\frac{1}{1-e^{-2\pi\iu f} e^{-\delta_{v,f,i}}} 
\prod_{f=0, 1\le i\le l_{v,0}} 
\frac{\delta_{v,0,i}}{1-e^{-\delta_{v,0,i}}}. 
\]
When the orbifold vector bundle $V$ admits the 
structure of a holomorphic orbifold vector bundle, 
the holomorphic Euler characteristic 
$\chi(V) := \sum_{i=1}^n (-1)^i \dim_\C H^i(\cX,V)$ 
is given by the Kawasaki-Riemann-Roch 
formula \cite{Kaw}: 
\begin{equation}
\label{eq:KRR}
\chi(V) = \int_{I\cX} \tch(V) \cup \tTd(T\cX).  
\end{equation} 
Note that $\chi(V)$ is an integer by definition. 
For a (not necessarily holomorphic) topological 
orbifold vector bundle $V$ on $\cX$, we \emph{define}  
$\chi(V)$ to be the right-hand side of the 
above formula (\ref{eq:KRR}). 
It follows from 
Kawasaki's index theorem \cite{Kaw2} for 
elliptic operators on orbifolds that 
$\chi(V)$ is an integer for any $V$. 
In fact, the right-hand side of 
(\ref{eq:KRR}) equals the index of an elliptic operator 
$\ov{\partial}+\ov{\partial}^* \colon 
V\otimes \cA_\cX^{0,\rm even} 
\to V\otimes \cA_\cX^{0,\rm odd}$, 
where $\ov{\partial}$ is a not necessarily integrable 
$(0,1)$-connection on $V$ and $\ov{\partial}^*$ 
is its adjoint with respect to a hermitian metric on $V$. 

Define a multiplicative characteristic class 
$\hGamma\colon K(\cX) \to H^*(I\cX)$ as follows:  
\[
\hGamma(V) := \bigoplus_{v\in \sfT} 
\prod_{0\le f<1} \prod_{i=1}^{l_{v,f}} 
\Gamma(1-f+\delta_{v,f,i}). 
\]
Here $\delta_{v,f,i}$ is the same as above. 
The Gamma function in the right-hand side 
should be expanded in Taylor series at $1-f>0$. 
%We set $\hGamma_\cX := \hGamma(T\cX)$. 
The $\hGamma$-class can be viewed as 
a funny ``square root" of the Todd class 
(more precisely, $\widehat{A}$-class). 
In fact, using the Gamma function equality 
$\Gamma(z) \Gamma(1-z) = \pi/\sin(\pi z)$, 
we find 
\begin{align*}
[(e^{\pi\iu \deg/2} \hGamma(V)) \cup &
\inv^* \hGamma(V)]_v \cup  
e^{\pi\iu (c_1(\pr^*V)|_{\cX_v} + \age_v(V))}  \\
&=(2\pi\iu)^{\sum_{f\neq 0} l_{v,f}} 
[(2\pi\iu)^{\deg/2} \tTd(V)]_v, 
\end{align*} 
where $\deg\colon H^*(I\cX)\to H^*(I\cX)$ is 
the ordinary grading operator defined by 
$\deg = p$ on $H^{p}(I\cX)$, 
$\age_v(V)=\sum_{0<f<1} f l_{v,f}$ is the 
age of $V$ along $\cX_v$,  
and $[\cdots]_v$ is the $H^*(\cX_v)$-component. 
In this sense, the $K$-group framing $\cZ_K\colon 
K(\cX) \to \Sol(\cX)$ below 
can be considered as a ``Mukai vector" 
in quantum cohomology. 

\begin{definition}
\label{def:Ktheoryintstr} 
We define the \emph{$K$-group framing} 
$\cZ_K\colon K(\cX) \to \Sol(\cX)$ 
of the space $\Sol(\cX)$ of flat sections by 
the formula: 
\begin{equation}
\label{eq:Kgroupframing} 
\begin{split} 
&\cZ_K(V) :=  
\cZ_{\rm coh}(\Psi(V)) = 
L(\tau,z)z^{-\mu}z^\rho \Psi(V), \\ 
&\text{where }
\Psi(V) 
:=  
(2\pi)^{-\frac{n}{2}} \hGamma(T\cX)\cup (2\pi\iu)^{\deg/2} 
\inv^* \tch(V).
\end{split} 
% L(\tau,z) z^{-\mu}z^\rho 
% \left( 
% (2\pi)^{-n/2} \hGamma(T\cX)\cup (2\pi\iu)^{\deg/2} 
% \inv^* \tch(V) 
% \right). 
\end{equation} 
Here $\cZ_{\rm coh}$ is the cohomology framing 
(\ref{eq:cohframing}), 
$L(\tau,z)z^{-\mu}z^\rho$ 
is the fundamental solution in Proposition 
\ref{prop:fundamentalsol}, 
$(2\pi\iu)^{\deg/2}\in \End(H^*(I\cX))$ 
is defined by $(2\pi\iu)^{\deg/2}|_{H^{2p}(I\cX)} 
=(2\pi\iu)^p$ and 
$\hGamma(T\cX)\cup$ is the cup product in 
$H^*(I\cX)$. 
The image $\Sol(\cX)_\Z:= \cZ_K(K(\cX)) \subset \Sol(\cX)$ 
of the $K$-group framing is called  
\emph{the $K$-theory integral structure} 
on the quantum cohomology.  
\end{definition} 

The notation $\cZ_K$ for the $K$-group framing 
is motivated by the \emph{central charge} in physics. 
Conjecturally, the integral 
\begin{equation}
\label{eq:centralcharge}
Z(V) := c(z) \int_\cX \cZ_K(V)(\tau,z) 
= c(z) (\unit, \cZ_K(V)(\tau,z))_{\rm orb} 
\end{equation}
with $c(z) = (2\pi z)^{\frac{n}{2}}/(2\pi\iu)^n$, $n=\dim \cX$  
gives the central charge  
of a B-type D-brane in the class $V$ 
at the point $\tau$ of the (extended) 
K\"{a}hler moduli space. 
This plays a central role in 
stability conditions on the derived category 
$D^b_{\rm coh}(\cX)$ \cite{Douglas, Bridgeland}. 
It would be very interesting 
to find an intrinsic explanation for 
the formula (\ref{eq:Kgroupframing}) 
from this point of view. 
In the language of quantum $D$-modules, 
$Z(V)$ is a coefficient of the unit section 
$\unit$ expressed in a $\nabla$-flat frame.

\begin{proposition}[
{\cite[Definition-Proposition 3.16]{I:real}}] 
{\rm (i)}  
The image $\Sol(\cX)_\Z$ of the $K$-group framing $\cZ_K$ 
is a lattice in $\Sol(\cX)$: 
\[
\Sol(\cX)_\Z \otimes_\Z \C = \Sol(\cX). 
\]

{\rm (ii)} 
The pairing $(\cdot,\cdot)_\Sol$ on $\Sol(\cX)$ 
corresponds to the Mukai pairing on $K(\cX)$ 
through the $K$-group framing $\cZ_K$: 
\[
(\cZ_K(V_1),\cZ_K(V_2))_{\Sol} = \chi(V_1 \otimes V_2^\vee) 
\]
Therefore, we have a $\Z$-valued pairing 
$\Sol(\cX)_\Z \times \Sol(\cX)_\Z \to \Z$. 

{\rm (iii)} 
For $\xi\in H^2(\cX,\Z)$, 
the Galois action $G^\Sol(\xi)$ on $\Sol(\cX)$ corresponds to 
the tensor by the orbifold line bundle $L_\xi^\vee$ 
(corresponding to $-\xi$) on $K(\cX)$:  
\[
\cZ_K(L_\xi^\vee\otimes V) = G^\Sol(\xi) (\cZ_K(V)).  
\] 
In particular, the lattice $\Sol(\cX)_\Z$ 
is invariant under the Galois action.  
\end{proposition} 
The statement (i) follows from the Adem-Ruan 
decomposition theorem \cite[Theorem 5.1]{AR}, 
which implies that $\tch\colon K(\cX)\to H^*(I\cX)$ 
is an isomorphism when tensored with $\C$. 
The statements (ii) and (iii) 
follow from straightforward calculations. 
It is somewhat surprising that many complicated terms 
finally give the Mukai pairing in (ii) via the 
Kawasaki-Riemann-Roch (\ref{eq:KRR}). 

\begin{remark}
The formula (\ref{eq:Kgroupframing}) 
arose in \cite{I:real} 
from the study of mirror symmetry for toric orbifolds. 
The mirror Landau-Ginzburg model has 
the natural integral structure and we can 
shift it to the quantum cohomology.  
Katzarkov-Kontsevich-Pantev \cite{KKP} 
also proposed essentially the same definition 
(for a rational structure)  
when $\cX$ is a manifold. 
Closely related results have been observed 
in the context of mirror symmetry. 
Calculations and conjecture of Hosono 
\cite{Hosono:IIA}, \cite[Conjecture 6.3]{Hosono} 
are compatible with the integral structure above; 
The works of Horja \cite{Horja1,Horja2} 
and Borisov-Horja \cite{Bor-Hor} strongly 
suggest a relation between $K$-group 
and quantum $D$-module. 
\end{remark} 

\begin{example}
\label{ex:integral} 
(i) $\cX=\Proj^1$. 
Let $\omega\in H^2(\Proj^1,\Z)$ be the integral 
K\"{a}hler class. We take $1,\omega$ as a basis 
of $H^*(\Proj^1)$. 
In terms of the cohomology framing 
$\cZ_{\rm coh}\colon H^*(\Proj^1) \cong \Sol(\Proj^1)$ 
in (\ref{eq:cohframing}), 
the Galois action and the pairing on 
$\Sol(\Proj^1)$ is represented by the matrices: 
\[
G^{\Sol}(\omega) = 
\begin{bmatrix}
1 & 0 \\
-2\pi\iu & 1
\end{bmatrix}, \quad 
(\cdot,\cdot)_{\Sol} = 
\begin{bmatrix}
2\pi & \iu \\ 
-\iu & 0 
\end{bmatrix}.  
\]
If an integral lattice $L$ in 
$H^*(\Proj^1)\cong \Sol(\Proj^1)$ 
is invariant under $G^{\Sol}(\omega)$ 
and if the restriction of $(\cdot,\cdot)_\Sol$ to 
$L$ gives a perfect pairing $L\times L\to \Z$, 
then $L$ must take the following form:  
\[
L = \Z \sqrt{\frac{n}{2\pi}}(1 + c \omega)  
\oplus \Z \iu \sqrt{\frac{2\pi}{n}} \omega
\]
for some $n\in \Z\setminus \{0\}$ and $c\in \C$. 
The $K$-theory integral structure corresponds 
to the choice $n=1$ and $c=-2\gamma$, where 
$\gamma=0.57721...$ is Euler's constant 
(coming from the $\hGamma$-class 
$\hGamma(T\Proj^1)=1-2\gamma \omega$). 

(ii) When $\cX=X$ is a Calabi-Yau threefold, 
the $\hGamma$ class is given by 
\[
\hGamma(TX) = 
1-\frac{\pi^2}{6} c_2(X) - \zeta(3) c_3(X)   
\]
where $\zeta(3)$ is the special value of 
Riemann's zeta function. From this, it follows that 
the central charges (\ref{eq:centralcharge}) 
of $\cO_{\rm pt}$, $\cO_{C}$, $\cO_{S}$ and $\cO$ 
(where $C,S$ are smooth curve and surface) 
restricted to $H^2(X)$ are 
\begin{align*}
Z(\cO_{\rm pt}) &= 1,  \\ 
Z(\cO_{C})  &=  
((1-g(C)) - \frac{\tau}{2\pi\iu} \cap [C], \\ 
Z(\cO_S) &=  
\frac{[S]^3}{8} +\frac{\chi(S)}{24} 
+ \frac{\tau}{2\pi\iu}\cap \frac{[S]^2}{2} 
+ \frac{d_{[S]}F_0(\tau)}{(2\pi\iu)^2}, \\
Z(\cO) &= 
-\frac{\zeta(3)}{(2\pi\iu)^3} \chi(X) 
- \frac{\tau}{2\pi\iu}\cdot \frac{c_2(X)}{24} + 
\frac{H(\tau)}{(2\pi\iu)^3}, 
\end{align*} 
where $\tau=\tau_{0,2}\in H^2(X)$, 
$g(C)$ is the genus of $C$, and 
$\chi(X)$ and $\chi(S)$ are the Euler numbers of $X$ and $S$. 
$F_0(\tau)$ is the genus zero potential of $X$ 
\[
F_0(\tau) := \frac{1}{6} \int_X \tau^3 + 
\sum_{d\in \Eff_X \setminus\{0\}} 
\corr{\phantom{,}}_{0,0,d} e^{\pair{\tau}{d}},  
\]
$d_{[S]}F_0$ is its derivative in the direction 
of the Poincar\'{e} dual of $[S]$ and 
$H(\tau) :=  2F_0(\tau)-\sum_i t^i \partial_i F_0(\tau)$. 
The zeta value $\zeta(3)$ also appeared 
in the quintic mirror calculation of 
Candelas-de la Ossa-Green-Parkes \cite{CDGP}.

(iii) When $\cX$ is a weak Fano compact toric orbifold, 
it is shown in \cite[Theorem 4.17]{I:real} that 
the central charge of the structure sheaf 
can be written as an oscillating integral 
of the mirror Landau-Ginzburg model 
$W_\tau \colon (\C^*)^n\to \C$:  
\[
Z(\cO_\cX)(\tau,z) 
= \frac{1}{(2\pi \iu)^{n}} 
\int_{\Gamma_\R\subset (\C^*)^n} 
e^{-W_{\tau}(y)/z} \frac{dy}{y}, \quad 
n= \dim_\C\cX.  
\]
Here $dy/y$ is an invariant holomorphic 
$n$-form on $(\C^*)^n$ and 
$\Gamma_\R$ is a non-compact cycle 
(Lefschetz thimble) in $(\C^*)^n$. 
(Strictly speaking, we need a ``mirror map" 
between $\tau\in H^2_{\rm CR}(\cX)$ in the left-hand side 
and the parameter $\tau$ in the Landau-Ginzburg 
potential $W_\tau$.) 
This shows that the integral structure 
in Definition \ref{def:Ktheoryintstr} is 
compatible with (and actually the same as) 
that of the mirror given by the lattice of Lefschetz thimbles. 
The Lefschetz thimble $\Gamma_\R$ corresponds to 
the structure sheaf $\cO_\cX$ 
and the oscillating form $e^{-W_\tau/z}(dy/y)$ 
corresponds to the unit section $\unit$ 
of the quantum $D$-module. 
See \cite{I:real} for more details. 

(iv) The $\hGamma$-class 
contains odd zeta values 
$\zeta(3), \zeta(5),\dots$ 
and products of Gamma values. 
When $\cX$ is holomorphic symplectic, however, 
the $\hGamma$-class is defined over $\Q(\zeta)[\pi]$ 
for some root of unity $\zeta$.  
This might be related to the fact that 
there is no quantum correction. 
\end{example} 

\begin{remark} 
We can consider the Grothendieck group of 
algebraic vector bundles or coherent sheaves on $\cX$ 
instead of topological $K$-groups.  
In this case, the $K$-theory integral structure 
is defined on the algebraic part of the orbifold 
cohomology $H_{\rm CR}^*(\cX)$, 
\emph{i.e.} cohomology classes on $I\cX$ 
which can be written as 
linear combinations of Poincar\'{e} duals 
of algebraic cycles with complex coefficients. 
The algebraic part of orbifold quantum cohomology 
makes sense due to the algebraic construction 
of orbifold Gromov-Witten theory \cite{AGV}. 
A theoretical difficulty is that we do not 
know if the orbifold Poincar\'{e} pairing is non-degenerate 
when restricted to the algebraic part of $H_{\rm CR}^*(\cX)$: 
This would be a consequence of 
the famous Hodge conjecture/Grothendieck standard conjecture. 
Apart from this point, many discussions in this paper 
can be equally applied to \emph{algebraic} $K$-theory 
integral structures. 
\end{remark} 

\subsection{Remark on non-compact case} 

Even when the space $\cX$ is non-compact, 
we can sometimes define the (orbifold) 
quantum cohomology. 
Non-compact local cases are important 
in the study of Ruan's conjecture. 
One standard way is to use the \emph{torus-equivariant} 
Gromov-Witten theory. 
If $\cX$ admits a torus action and 
the fixed point set is compact, 
we can define torus-equivariant 
orbifold Gromov-Witten invariants 
using the Atiyah-Bott style localization on 
the moduli space of stable maps \cite{GP}. 
In good cases, we can take the non-equivariant 
limit and have the non-equivariant quantum cohomology. 
In general, we can define Gromov-Witten invariants 
if the moduli space of stable maps to $\cX$ is compact\footnote{
However, the degree zero moduli space 
always has a non-compact component, 
so we indeed need that the evaluation map is proper as stated. 
This is particularly relevant to the orbifold case 
where degree zero moduli spaces give a lot of 
non-trivial Gromov-Witten invariants.}. 
More generally, even when the moduli space may not be compact,  
if the evaluation map from the moduli space 
to the inertia stack $I\cX$ is proper, 
we can define the quantum product 
by the push-forward by the evaluation map 
at the ``last" marked point. 
As suggested in \cite{BG}, 
this happens for example when $X$ is semi-projective, 
\emph{i.e.} projective over an affine scheme. 
In this section, assuming the existence of 
a well-defined orbifold quantum cohomology for 
a non-compact space, we describe a possible  
framework for $K$-theory integral structures 
in this case. 

Assume that the (non-equivariant) quantum cohomology 
of $\cX$ is well-defined. 
Quantum cohomology defines the Dubrovin connection 
and the quantum $D$-module 
in the same fashion as in Definition \ref{def:QDM}. 
The discrete Galois symmetry in Lemma 
\ref{lem:Galois} is also well-defined.  
A problem in non-compact case is 
that the orbifold Poincar\'{e} pairing on 
$H_{\rm CR}^*(\cX)$ is degenerate. 
However, we have a non-degenerate 
pairing between $H_{\rm CR}^*(\cX)$ and 
the \emph{compactly supported} orbifold cohomology 
$H_{\rm CR,c}^*(\cX)$, which is defined to be 
the direct sum of compactly supported cohomology 
groups of the inertia components $\cX_v$ (with 
the same grading shift as before):   
\[
(\cdot,\cdot)_{\rm orb} \colon 
H_{\rm CR,c}^*(\cX) \times H_{\rm CR}^*(\cX) \to \C. 
\] 
This pairing defines the \emph{dual Dubrovin connection} 
on the $H_{\rm CR,c}^*(\cX)$-bundle 
$F_{\rm c}   := H_{\rm CR,c}^*(\cX)\times (U\times \C) \to U\times \C$: 
\begin{align*}
&\nabla_i = \parfrac{}{t^i} + 
\frac{1}{z} (\phi_i \circ_\tau)^\dagger, \\ 
&\nabla_{z\partial_z} 
= z \parfrac{}{z} - \frac{1}{z} (E\circ_\tau)^\dagger + \mu 
\end{align*} 
where $(\phi_i \circ_\tau)^\dagger, 
(E\circ_\tau)^\dagger \in \End(H_{\rm CR,c}^*(\cX))$ 
are the adjoint operators with respect to 
$(\cdot,\cdot)_{\rm orb}$.  
We call $(F_{\rm c} , \nabla)$ the 
\emph{compactly supported quantum $D$-module}. 
Note that the dual product $(\phi_i\circ_\tau)^\dagger$
is defined by essentially the same formula as 
the original product: 
$(\alpha\circ_\tau \beta, \gamma)_{\rm orb} = 
(\alpha, (\beta\circ_\tau)^\dagger \gamma)_{\rm orb}$ may be 
defined by the right-hand side of (\ref{eq:quantumproduct}) 
with $\alpha,\beta\in H_{\rm CR}^*(\cX)$ 
and $\gamma\in H_{\rm CR,c}^*(\cX)$ 
(under the assumption that the evaluation map is proper). 
Tautologically, one has a $\nabla$-flat pairing: 
\[
(-)^*\cO(F_{\rm c}  ) \otimes \cO(F) \to \cO_{U\times \C} 
\]
induced from the orbifold Poincar\'{e} pairing, 
where recall that $(-)\colon U\times \C\to U\times \C$ 
is the map sending $(\tau,z)$ to $(\tau,-z)$. 
One has a natural map 
\[
(F_{\rm c} , \nabla) \to (F,\nabla) 
\]
induced from $H_{\rm CR,c}^*(\cX) \to H_{\rm CR}^*(\cX)$. 
The fundamental solution in Proposition 
\ref{prop:fundamentalsol} also makes sense. 
We have two fundamental solutions 
$\tilde{L}(\tau,z)$, $L(\tau,z)$ taking values 
in $\End(H_{\rm CR,c}^*(\cX))$ and $\End(H_{\rm CR}^*(\cX))$
respectively such that 
\begin{gather*} 
\nabla (\tilde{L}(\tau,z)z^{-\mu} z^\rho \varphi) =0, \quad 
\nabla (L(\tau,z) z^{-\mu} z^\rho \phi) =0, \\ 
(\tilde{L}(\tau,-z)\varphi, L(\tau,z)\phi)_{\rm orb} = 
(\varphi,\phi)_{\rm orb},  
\end{gather*}  
where $\varphi \in H^*_{\rm CR,c}(\cX)$ 
and $\phi\in H^*_{\rm CR}(\cX)$. 
Here again, $\tilde{L}(\tau,z)$ and $L(\tau,z)$ 
can be defined by the same formula (\ref{eq:fundamentalsol}), 
with different domains of definitions\footnote{ 
Here one of the dual pairs $\{\phi_k\}, \{\phi^k\}$ 
in (\ref{eq:fundamentalsol}) 
should be taken from $H_{\rm CR,c}^*(\cX)$ 
and the other from $H_{\rm CR}^*(\cX)$.   
We take $\phi^k\in H_{\rm CR,c}^*(\cX)$ 
when defining $L(\tau,z)$ and 
$\phi_k\in H_{\rm CR,c}^*(\cX)$ when defining 
$\tilde{L}(\tau,z)$.}.  
The spaces $\Sol(\cX)$, $\Sol_{\rm c}(\cX)$ 
of multi-valued flat sections of $F$, $F_{\rm c}$  
are defined in the same way as 
in Definition \ref{def:spaceofflatsections}. 
The symmetries in Lemma \ref{lem:Galois}   
act on these spaces as automorphisms 
preserving the pairing: 
\[
(\cdot,\cdot)_{\Sol}\colon 
\Sol_{\rm c} (\cX) \times \Sol(\cX) \to \C, \quad 
(s_1,s_2)\mapsto (s_1(\tau,e^{\pi\iu}z), s_2(\tau,z))_{\rm orb} 
\]
Likewise, the formula (\ref{eq:Kgroupframing}) 
defines $K$-group framings  
\[ 
\cZ_K \colon K(\cX) \to \Sol(\cX), \quad 
\cZ_{K,\rm c} \colon K_{\rm c}(\cX) \to \Sol_{\rm c}(\cX)
\] 
where $K_{\rm c}(\cX)$ 
is the \emph{compactly supported} $K$-group.  
(We need to use $\tilde{L}(\tau,z)$ instead of 
$L(\tau,z)$ in (\ref{eq:Kgroupframing}) 
for the compact support version.) 
For example, when $\cX$ is of the form $M/G$, 
one can define $K_{\rm c}(\cX)$ as the $G$-equivariant 
reduced $K$-group $\widetilde{K}_G^0(M^+)$ of the 
one-point compactification $M^+$ of $M$ 
(as in \cite{Seg}).   
One can also use the 
Grothendieck group $K_Z(\cX)$ of coherent sheaves 
on $\cX$ supported on a compact set $Z$. 
In non-compact case, 
the definition of $K(\cX)$ may be subject to change 
\emph{e.g.} we may need to include perfect complexes 
or infinite dimensional bundles 
\emph{c.f.} \cite{Totaro}. 
We will not pursue a more precise formulation here. 
Note that we have a well-defined central charge 
$Z(V):=c(z) \int_\cX \cZ_{K,\rm c}(V)$ 
for $V\in K_{\rm c}(\cX)$.

\begin{example}[\emph{c.f.} {\cite[Example 6.5]{I:real}}]
(i) $\cX=[\C^2/G]$ where $G$ is a finite subgroup of $SL(2,\C)$. 
The inertia stack $I\cX$ is given by 
\[
I\cX = \cX \sqcup \bigsqcup_{(g)\neq 1} \cX_{(g)}, \quad 
\cX_{(g)} = [\{ 0 \}/C(g)] \quad (g\neq 1),   
\] 
where $(g)$ is a conjugacy class of $G$, $g\in G$,  
and $C(g)$ is the centralizer of $g$ in $G$. 
Let $\unit$ be the unit class supported on $\cX$ 
and $\unit_{(g)}\in H^*_{\rm CR}(\cX)$ be 
the unit class supported on $\cX_{(g)}$. 
The grading is given by 
\[
\deg \unit =0, \quad \deg \unit_{(g)}=2 \quad (g\neq 1).  
\] 
Since $\cX$ is holomorphic symplectic, 
there is no quantum deformation and 
$\circ_\tau$ is trivial: 
$\unit\circ_\tau\unit_{(g)}=\unit_{(g)}$ 
and all other products are zero. 
(We can get non-trivial quantum cohomology 
by considering the equivariant version.)  
The $\hGamma$-class is given by 
\[
\hGamma(T\cX) = 
\unit \oplus \bigoplus_{(g)\neq (1)}
\frac{\pi}{\sin(\pi f_g)} \unit_{(g)} 
\in H^0(I\cX) 
\]
where $0\le f_g\le 1/2$ is the rational number 
such that the eigenvalues of $g\in SL(2,\C)$ 
are $\exp(\pm 2\pi\iu f_g)$.   
Let $\beta, \unit_{(g)}$ ($g\neq 1$) be 
compactly supported cohomology classes 
on $\cX$, $\cX_{(g)}$   
such that 
\[
(\beta, \unit)_{\rm orb} = \frac{1}{|G|}, \quad 
(\unit_{(g)}, \unit_{(g^{-1})})_{\rm orb} 
=\frac{1}{|C(g)|}  \quad (g\neq 1).  
\] 
Here $\deg \beta =4$. 
We consider the Grothendieck group $K_0^G(\C^2)$ 
of $G$-equivariant coherent sheaves on $\C^2$ 
supported at the origin. 
A finite dimensional representation $\varrho$ of $G$ 
defines a $G$-equivariant sheaf 
$\cO_0\otimes \varrho$ on $\C^2$.  
These sheaves generate $K^G_0(\C^2)$ 
and the Galois action corresponds to the 
tensor product by a one-dimensional representation.   
By the equivariant Koszul resolution: 
\[
0 \to 
\cO_{\C^2}\otimes \varrho \to  
\cO_{\C^2}\otimes \varrho\otimes Q^\vee \to  
\cO_{\C^2}\otimes \varrho \to 
\cO_0\otimes \varrho \to 0, 
\]
where $Q=\C^2$ is the standard $G$-representation 
defined by the inclusion $G\subset SL(2,\C)$, 
we compute the Chern character as 
\[
\tch(\cO_0\otimes \varrho) = (\dim\varrho) \beta \oplus 
\bigoplus_{(g)\neq (1)} 
\Tr(g|\varrho\otimes (\C^2-Q)) \unit_{(g)} 
\in H_{\rm c}^*(I\cX). 
\] 
Here $\Tr(g|\varrho\otimes (\C^2-Q))$ is the 
trace of $g$ on the virtual representation 
$\varrho\otimes (\C^2-Q)$. 
Therefore, using 
$\tilde{L}(\tau,z)=\exp(-(\tau\circ_\tau)^\dagger/z)$, 
we find 
\begin{equation}
\label{eq:charge_C2}
Z(\cO_0\otimes \varrho) = e^{-t^0/z} 
\left(\frac{\dim \varrho}{|G|} + 
\sum_{(g)\neq 1} \frac{\Tr(g|\varrho) 
\sin(\pi f_g)}{|C(g)| \pi}t^{(g)}\right),    
\end{equation} 
where we put $\tau =t^0 \unit + 
\sum_{(g)\neq 1} t^{(g)} \unit_{(g)}$.   
The simplest central charge is given by 
the regular representation $\varrho_{\rm reg}$: 
\[
Z(\cO_0\otimes \varrho_{\rm reg}) = e^{-t^0/z}. 
\] 
The vector $[\cO_0\otimes \varrho_{\rm reg}]\in K^G_0(\C^2)$ 
is invariant under every Galois action. 
% Moreover, when $G=\Z/n\Z$, this vector 
% is characterized as a primitive vector 
% invariant under every Galois action. 

(ii) $\cX = \C^3/G$ where $G$ is 
a finite subgroup of $SL(3,\C)$. 
This case can have a non-trivial 
(non-equivariant) quantum cohomology.  
The inertia stack $I\cX$ is given by 
\[
I\cX = \cX \sqcup \bigsqcup_{(g)\neq (1)} \cX_{(g)}, \quad 
\cX_{(g)} = [(\C^3)^g/C(g)],   
\]
where $(\C^3)^g\subset \C^3$ 
is the subspace fixed by $g$. 
The ordinary and compactly supported 
orbifold cohomology are 
\begin{align*}
H_{\rm CR}^*(I\cX) &= \C \unit \oplus 
\bigoplus_{(g)\neq 1} \C\unit_{(g)}, \\ 
H_{\rm CR,c}^*(I\cX) &= \C \alpha \oplus 
\bigoplus_{(g):n_g=1} \C \beta_{(g)} 
\oplus \bigoplus_{(g):n_g=0} \C \unit_{(g)},  
\end{align*} 
where $n_{g}=\dim \cX_{(g)}$. 
Here $\unit_{(g)}$ is the unit class 
supported on $\cX_{(g)}$ and 
$\alpha$, $\beta_{(g)}$ are top classes 
on $\cX$, $\cX_{(g)}$ respectively (with $n_g=1$) 
such that 
\[
(\alpha, \unit)_{\rm orb} = \frac{1}{|G|}, \quad 
(\beta_{(g)},\unit_{(g^{-1})})_{\rm orb}
=(\unit_{(g)},\unit_{(g^{-1})})_{\rm orb} = 
\frac{1}{|C(g)|}. 
\]
Note that $\deg \unit_{(g)} = 2\iota_{(g)}$,  
$\iota_{(g)}=1$ if $n_{g}=1$, 
$\deg \alpha=6$ and $\deg \beta_{(g)}= 4$. 
When $n_g=1$, let $0< f_g\le 1/2$ be a rational number 
such that $1,e^{\pm 2\pi\iu f_g}$ are the eigenvalues 
of $g\in SL(3,\C)$. 
When $n_g=0$, 
let $0< f_{g,1}\le f_{g,2}\le f_{g,3}< 1$ be rational 
numbers such that $e^{2\pi\iu f_{g,j}}$, $j=1,2,3$,  
are the eigenvalues of $g$. 
Consider again the Grothendieck group $K^G_0(\C^3)$ 
of $G$-equivariant coherent sheaves 
supported at the origin. 
A finite dimensional representation $\varrho$ of $G$ 
gives a class $[\cO_0\otimes \varrho]\in K^G_0(\C^3)$.  
This yields a dual flat section 
$\cZ_{K,\rm c}(\cO_0\otimes \varrho) = 
\tilde{L}(\tau,z)z^{-\mu} \Psi(\cO_0\otimes \varrho)$ 
with $\Psi(\cO_0\otimes \varrho)$ given by 
\[
(\dim\varrho) \alpha \oplus 
\bigoplus_{(g):n_g=1} 
(-1)A^\varrho_{g^{-1}} \beta_{(g)}\oplus 
\bigoplus_{(g):n_g=0} 
(-1)^{1+\iota_{(g)}} B^\varrho_{g^{-1}} \unit_{(g)}.    
\]
Here 
\[
A^\varrho_g = \Tr(g|\varrho) \frac{\sin(\pi f_g)}{\pi}, \quad 
B^\varrho_g = 
\frac{\Tr(g|\varrho)}{\prod_{j=1}^3 \Gamma(1-f_{g,j})}. 
\]
The corresponding central charge restricted to 
$H^2_{\rm CR}(\cX)$ is 
\begin{equation}
\label{eq:charge_C3}
Z(\cO_0\otimes \varrho) 
=\frac{\dim\varrho}{|G|} + 
\sum_{(g):n_g=1} \frac{A^\varrho_g}{|C(g)|} t^{(g)}
+ \sum_{(g):n_g=0} 
B^\varrho_g F_{0,(g^{-1})}(\tau), 
\end{equation} 
where $\tau = \sum_{\iota_{(g)}=1} t^{(g)} \unit_{(g)} 
\in H^2_{\rm CR}(\cX)$ and   
\begin{equation}
\label{eq:C3_pot}
F_{0,(g^{-1})}(\tau) = 
\begin{cases} 
t^{(g)}/|C(g)|, & \iota_{(g)}=1, \\   
\sum_{m\ge 2} \frac{1}{m!} 
\corr{\unit_{(g^{-1})}, \tau,\dots,\tau}_{0,m+1,0}, 
& \iota_{(g)}=2. 
\end{cases} 
\end{equation} 
This follows from 
$Z(\cO_0\otimes \varrho) = 
(\tilde{L}(\tau,z)^\dagger \unit, 
z^{-\mu}\Psi(\cO_0\otimes \varrho))_{\rm orb}$ 
and the formula for the \emph{$J$-function} 
$J(\tau,-z) = \tilde{L}(\tau,z)^\dagger \unit$: 
\[
J(\tau,-z) = \unit - \frac{\tau}{z} + 
\sum_{\iota_{(g)}=2} 
F_{0,(g^{-1})}(\tau) |C(g)|\frac{\unit_{(g)}}{z^2}. 
\] 
Again the regular representation $\rho_{\rm reg}$ 
gives the simplest charge $1$. 
The $\Gamma$-product 
$\prod_{j=1}^3 \Gamma(1-f_{g,j})$ in the central charge 
may have something to do with the 
Chowla-Selberg formula \cite{CS}. 
\end{example} 
 
% \subsection{Aside: Hodge theoretic properties}
% The $K$-theory integral structure also satisfy 
% nice properties from the viewpoint of Hodge theory 
% \cite{Sab, Her, HS, KKP}. 
% This section can be read independently of other sections. 
% The readers who are not interested in Hodge 
% theoretic properties can safely skip this section. 
 
\section{Ruan's conjecture} 
\label{sec:Ruan}

We incorporate our $K$-theory picture 
into the Ruan's conjecture \cite{Ruan:crepant1, Ruan:crepant2} 
and discuss what follows from this. 
We propose the picture that a conjectural 
isomorphism between $K$-theory induces 
an isomorphism of quantum $D$-modules 
via the $K$-group framing (\ref{eq:Kgroupframing}). 

Ruan's conjecture can be discussed in many situations.  
It basically asserts that two birational spaces 
$\cX_1$, $\cX_2$ in a ``crepant" relationship 
have isomorphic (orbifold) quantum cohomology 
under a suitable identification of quantum parameters. 
One of such relationships is a \emph{crepant resolution}. 
Let $\cX$ be a Gorenstein orbifold 
without generic stabilizers, 
\emph{i.e.} the automorphism group at every point $x$ 
is contained in $SL(T_x\cX)$.  
Then the canonical line bundle $K_\cX$ of $\cX$ 
becomes the pull-back of $K_X$ of the coarse moduli space $X$.  
A resolution of singularity $\pi\colon Y\to X$ 
is called \emph{crepant} if $\pi^*K_X \cong K_Y$.  
We can regard $Y$ and $\cX$ as two  
different crepant resolutions of the same space $X$: 
\[
\begin{CD}
\cX @>>> X @<<< Y.  
\end{CD}
\] 
In this case, Ruan's conjecture for a pair $(\cX,Y)$ 
is called the \emph{crepant resolution conjecture} 
and has been studied in many literatures  
\cite{BGP,Perroni,BG,CIT,B-Gh:ADE,BMP,B-Gh:poly,Coates}. 
Ruan's conjecture have been discussed also for \emph{flops}. 
Li-Ruan \cite{LR} showed that the 
quantum cohomology is invariant under flops 
between Calabi-Yau 3-folds. 
Recently, this was generalized 
to the case of simple $\Proj^r$-flops 
and Mukai flops \cite{LLW} in any dimension.   
The case of certain singular flops between 
orbifolds are also studied in \cite{CLZ,CLZZ}. 

More generally, Ruan's conjecture 
may hold for \emph{$K$-equivalences}. 
We say that two smooth Deligne-Mumford stacks 
$\cX_1$, $\cX_2$ are $K$-equivalent if there exist  
a smooth Deligne-Mumford stack $\cX$ 
and a diagram of projective birational morphisms 
\begin{equation}
\label{eq:Kequiv}
\begin{CD}
\cX_1 @<{p_1}<< \cX @>{p_2}>> \cX_2  
\end{CD} 
\end{equation} 
such that $p_1^*K_{\cX_1} \cong p_2^*K_{\cX_2}$. 
% In dimension three, two 
% $K$-equivalent smooth varieties are 
% connected by a sequence of flops. 
The most general form of Ruan's conjecture 
would be the invariance of quantum cohomology 
under \emph{$D$-equivalences}, 
\emph{i.e.} the equivalence of derived categories of coherent sheaves. 
It is conjectured in \cite{Kaw:KD} that 
$K$-equivalence is equivalent to $D$-equivalence 
for smooth birational varieties,  
but $D$-equivalence does not imply 
birational equivalence in general. 
An interesting example is reported  
\cite{Ro, Ho-Ko} where the Gromov-Witten theories 
of non-birational but $D$-equivalent Calabi-Yau 3-folds 
have the same mirror family and, 
in particular, should be equivalent. 

One striking feature in Ruan's conjecture 
is that we need the \emph{analyticity} 
of the quantum cohomology. In the crepant 
resolution conjecture, the orbifold quantum cohomology 
is identified with the \emph{expansion} of the manifold 
quantum cohomology around a point where 
the quantum parameter $q=e^{\tau_{0,2}}$ is 
a root of unity. 
In the flop conjecture, two quantum cohomology 
are identified under the transformation 
$q\mapsto q^{-1}$, where $q$ is the parameter 
of the exceptional curve. 

\subsection{A picture of the global quantum $D$-module} 
\label{subsec:GQDM} 

Let $\cX_1,\cX_2$ be a pair of smooth Deligne-Mumford 
stacks for which Ruan's conjecture is expected to hold. 
For a complex analytic space $\cM$, 
let $\pi\colon \cM\times \C\to \cM$ be 
the projection to the first factor, $z$ 
be the co-ordinate on the $\C$ factor 
and $(-)\colon \cM\times \C\to \cM\times \C$ 
be the map sending $(\tau,z)$ to $(\tau,-z)$ 
as before.

\begin{Picture}[Global quantum $D$-modules: 
See Figure \ref{fig:GKM}] 
\label{Pic:GQDM} 
There exists a \emph{global quantum $D$-module} 
$(F,\nabla,(\cdot,\cdot)_F,F_\Z)$ 
over a \emph{global K\"{a}hler moduli space} $\cM$ 
given by the following data: 

---A connected complex analytic space $\cM$;

---A holomorphic vector bundle $F$ of rank $N$ 
over $\cM \times\C$;   

---A meromorphic flat connection $\nabla$ 
on $F$ (with poles along $z=0$): 
\[
\nabla \colon \cO(F) \to \cO(F)(\cM\times \{0\})
\otimes_{\cO_{\cM\times \C}} 
(\pi^*\Omega_{\cM}^1 \oplus 
\cO_{\cM\times \C} \frac{dz}{z});  
\]

---A non-degenerate, $\nabla$-flat pairing $(\cdot,\cdot)_F$: 
\[
(\cdot,\cdot)_F \colon 
(-)^* \cO(F) \otimes \cO(F) \to \cO_{\cM\times \C}; 
\]

---An integral local system ($\Z^N$-subbundle) 
$F_\Z\to \cM\times \C^*$  
underlying the flat vector bundle $F|_{\cM\times \C^*}$ 
such that 
\[
F_\Z \subset \Ker(\nabla), \quad 
F|_{\cM\times \C^*} = F_\Z \otimes \C, \quad 
((-)^*F_\Z, F_\Z)_F \subset \Z.   
\]
% We have a connected complex analytic space $\cM$ 
% and a holomorphic vector bundle $F$ of rank $N$ 
% over $\cM\times \C$ endowed 
% with a meromorphic flat connection $\nabla$ 
% \[
% \nabla \colon \cO(F) \to \cO(F)(\cM\times \{0\})
% \otimes_{\cO_{\cM\times \C}} 
% (\pi^*\Omega_{\cM}^1 \oplus 
% \cO_{\cM\times \C} \frac{dz}{z}) 
% \]
% and a non-degenerate, $\nabla$-flat 
% pairing $(\cdot,\cdot)_F$ 
% \[
% (\cdot,\cdot)_F \colon 
% (-)^* \cO(F) \otimes \cO(F) \to \cO_{\cM\times \C}.  
% \]
% Here $\cM$ is the conjectural  
% \emph{global K\"{a}hler moduli space} 
% and $(F,\nabla,(\cdot,\cdot)_F)$ is 
% the \emph{global quantum $D$-module}. 
% The flat bundle $F$ restricted to $\cM\times \C^*$ 
% contains an integral local system ($\Z^N$-subbundle) 
% $F_\Z \to \cM\times \C^*$ such that 
% \[
% F_\Z \subset \Ker(\nabla), \quad 
% F|_{\cM\times \C^*} = F_\Z \otimes \C, \quad 
% ((-)^*F_\Z, F_\Z)_F \subset \Z.   
% \]
We postulate that the tuple 
$(F,\nabla,(\cdot,\cdot)_F,F_\Z)$ satisfies 
the following. 

{\rm (i)} 
There exist open subsets $V_i\subset \cM$, $i=1,2$,  
such that $V_i$ is identified with the base space 
of the quantum $D$-module $QDM(\cX_i)$: 
\[
V_i \cong U_i/H^2(\cX_i,\Z),  
\]
and that the restriction of $(F,\nabla, (\cdot,\cdot)_F)$ 
to $V_i\times \C$ is isomorphic to $QDM(\cX_i)$:  
\[
(F,\nabla,(\cdot,\cdot)_F)|_{V_i\times \C} \cong QDM(\cX_i), 
\quad i=1,2.   
\]
Here $U_i\subset H_{\rm CR}^*(\cX_i)$  
is the convergence domain of the quantum product 
in Assumption \ref{as:convergence} 
and $U_i/H^2(\cX,\Z)$ is the quotient by the Galois action. 
Moreover, this isomorphism matches 
the integral local system $F_\Z$ with 
the $K$-theory integral structure of $QDM(\cX_i)$ 
in Definition \ref{def:Ktheoryintstr}. 

% {\rm (ii)} 
% The flat vector bundle $(F,\nabla)$ 
% restricted to $\cM\times \C^*$ 
% underlies an integral local system. 
% This integral structure on $F|_{V_i\times \C^*}$ 
% and the $K$-theory integral structure on $QDM(\cX_i)$ 
% (Definition \ref{def:Ktheoryintstr}) 
% coincide under the isomorphism in {\rm (i)}. 
% {\rm (ii)} 
% Let $\Sol_i$ be the space of multi-valued $\nabla$-flat 
% sections of $F$ over $V_i\times \C^*$. 
% There exists a path $\gamma\colon [0,1]\to \cM$ 
% from $x_1=\gamma(0)\in V_1$ 
% to $x_2=\gamma(1)\in V_2$ such that 
% the analytic continuation $P_\gamma$ 
% of $\nabla$-flat sections 
% along the path $\hat{\gamma}=(\gamma,1)\colon 
% [0,1]\to \cM\times \C^*$ is induced from 
% the isomorphism $\U_K$ in (\ref{eq:Kiso}). 
% \begin{align*} 
% \begin{CD}
% K(\cX_1) @>{\U_K}>> K(\cX_2) \\ 
% @V{\cZ_K}VV           @V{\cZ_K}VV  \\ 
% \Sol_1 @>{P_\gamma}>> \Sol_2.   
% \end{CD} 
% \end{align*} 
% By a slight abuse of notation,  
% the vertical map $\cZ_K$ denotes the composite  
% of the $K$-theory framing (\ref{eq:Kgroupframing}) 
% and the isomorphism $\Sol(\cX_i) \cong \Sol_i$. 

(ii) 
Assume that $\cX_1$ and $\cX_2$ are $K$-equivalent 
(\ref{eq:Kequiv}) and also related by a birational 
correspondence 
\begin{equation} 
\label{eq:contract} 
\begin{CD} 
\cX_1 @>{\pi_1}>> Z @<{\pi_2}<< \cX_2 
\end{CD} 
\end{equation} 
such that $\pi_1\circ p_1 = \pi_2 \circ p_2$. 
Take base points $x_i\in V_i$. 
For a line bundle $L$ on $Z$, 
denote by $l_i(L)\in \pi_1(V_i,x_i)$ 
the homotopy class of a loop 
given by the class 
$[\pi_i^*(L)]\in H^2(\cX_i,\Z)$. 
(Recall that $V_i \cong U_i/H^2(\cX_i,\Z)$.)
There exists a path $\gamma\colon [0,1]\to \cM$ 
from $\gamma(0)=x_1$ to $\gamma(1)=x_2$ such that 
$\gamma_*(l_1(L)) = l_2(L)$ 
for any line bundle $L$ on $Z$. 
Here $\gamma$ is independent of $L$. 
\end{Picture} 

As far as the author knows, all the concrete 
examples of global quantum $D$-modules 
arise from mirror symmetry. 
For example, in the case of 
toric flops or toric crepant resolutions 
(and complete intersections in them), 
we can construct a global quantum $D$-module 
using the mirror Landau-Ginzburg model 
and $\cM$ is identified with the complex moduli 
space of the mirror \cite{CIT, CCIT:comp, Coates}. 
The space of stability conditions on the 
derived category $D^b_{\rm coh}(\cX_i)$ 
due to Douglas and Bridgeland \cite{Douglas,Bridgeland} 
gives a candidate for the universal 
cover of $\cM$. 
Another conjectural candidate 
(though being infinitesimal) 
is the space of $A_\infty$-deformations 
of the derived Fukaya category of $\cX_i$. 

We assume the existence of a global quantum $D$-module 
$F$ connecting $QDM(\cX_1)$ and $QDM(\cX_2)$. 
Choosing a path $\gamma\colon [0,1]\to \cM$  
from a point $x_1\in V_1$ to a point $x_2\in V_2$, 
we have an analytic continuation map $P_\gamma$ 
of flat sections 
\begin{equation}
\label{eq:analyticcont}
P_\gamma\colon \Sol(\cX_1) \to \Sol(\cX_2) 
\end{equation} 
along the path $\hat{\gamma}=(\gamma,1)
\colon [0,1]\to \cM\times \C^*$. 
Here by (i), we identified the space of flat sections of 
$F$ over $V_i\times \C^*$ with $\Sol(\cX_i)$. 
This preserves the $K$-theory integral structures 
$P_\gamma(\Sol(\cX_1)_\Z) = \Sol(\cX_2)_\Z$  
and the pairing $(\cdot,\cdot)_\Sol$. 
Then it would be natural to conjecture the following. 

\begin{conjecture}  
\label{conj:Ktheoryisom} 
For each path $\gamma$, 
there exists an isomorphism of $K$-groups 
\begin{equation}
\label{eq:Kiso}
\U_{K,\gamma}\colon K(\cX_1) \to K(\cX_2) 
\end{equation} 
which induces the analytic continuation map 
$P_\gamma$ in (\ref{eq:analyticcont}) 
through the $K$-group framing (\ref{eq:Kgroupframing}). 
$\U_{K,\gamma}$ preserves the Mukai pairing 
$\chi(\U_{K,\gamma}(V_1)\otimes \U_{K,\gamma}(V_2)^\vee) =
\chi(V_1\otimes V_2^\vee)$. 
Note that $\U_{K,\gamma}$ 
gives the full relationships between 
$QDM(\cX_1)$ and $QDM(\cX_2)$ 
modulo the problem of analytic continuation. 
\end{conjecture} 

We expect that the $K$-group isomorphisms $\U_{K,\gamma}$ 
are given by geometric correspondences 
such as Fourier-Mukai transformations  
\cite{BKR, Kaw:logcrepant}. 
This conjecture is compatible with 
Borisov-Horja's result \cite{Bor-Hor}, 
where they identified the $K$-group of 
toric Calabi-Yau orbifold with 
the space of solutions to the GKZ system 
and also identified 
the analytic continuation 
of GKZ solutions with 
the Fourier-Mukai transformations 
between $K$-groups. 
If the path $\gamma$ is the same as what  
appeared in (ii) of Picture \ref{Pic:GQDM}, 
we also expect that $\U_{K,\gamma}$ commutes with 
the actions of line bundles pulled back from $Z$, 
\emph{i.e.} $\U_{K,\gamma}(\pi_1^*(L)\otimes V) = 
\pi_2^*(L)\otimes \U_{K,\gamma}(V)$ for a line bundle $L$ on $Z$. 
This is compatible with (ii) in Picture \ref{Pic:GQDM} 
and the fact that 
the tensor by $\pi_i^*L$ on $K(\cX_i)$ corresponds 
to the monodromy (Galois) action on $\Sol(\cX_i)$ 
along the loop $l_i(L)$. 

\begin{remark} 
(i) Unlike the original quantum $D$-module, 
the \emph{global} quantum $D$-module $F$ is 
not a priori trivialized in the standard way. 
This is an important point in this formulation. 
In fact, for the crepant resolution of $\C^3/\Z_3$ 
(or its compactification $\Proj(1,1,1,3)$), 
$F$ has different trivializations over $V_1$ and $V_2$ 
\cite{ABK, CIT}.   
Here different trivializations  
correspond to different Frobenius/flat structures  
on the base $\cM$. 

(ii) The flat connection can have poles along $z=0$. 
For a local section $s$ of $F$ around $z=0$, 
$\nabla_X s$ has a pole of order $\le 1$ along $z=0$ 
for $X\in T\cM$ and $\nabla_{\partial_z} s $ has 
a pole of order $\le 2$ along $z=0$. 

% (iii) 
% The condition (iii) of Picture \ref{Pic:GQDM} 
% means that the equivariance 
% of $\U_K$ with respect to pulled-back line bundles 
% is induced from the monodromy equivariance of 
% the analytic continuation $P_\gamma$. 
% Here recall that the tensor product by $\pi_i^*L$  
% corresponds to the monodromy along $l_i(L)$ 
% under the $K$-group framing (\ref{eq:Kgroupframing}). 

(iii) The $K$-theory isomorphism (\ref{eq:Kiso}) 
depends on the choice of a path $\gamma$. 
It would be very interesting to study 
the global monodromy of $(F,\nabla,(\cdot,\cdot)_F,F_\Z)$.  
\end{remark} 

\begin{remark}
In the context of Ruan's conjecture,  
the picture of the global quantum $D$-module  
has been proposed in \cite{CIT}, \cite{Coates-R} 
in terms of the Givental formalism.  
An integral structure was incorporated 
in this picture in \cite{I:real}. 
The structure analogous to 
the global quantum $D$-module 
$(F,\nabla,(\cdot,\cdot)_F, F_\Z)$ 
first emerged in singularity theory \cite{SaitoK} 
and have been studied under various names:  
\emph{Frobenius manifolds} \cite{Dub:2D};  
\emph{semi-infinite Hodge structures} \cite{Bar}; 
\emph{TE(R)P structures} \cite{Her, HS}; 
\emph{twistor structures} \cite{Simpson, Sabbah:twistor}; 
\emph{non-commutative Hodge structures} \cite{KKP} etc. 
\end{remark} 

\subsection{Family of algebras: 
isomorphism of $F$-manifolds} 
\label{subsec:Fmanifold}

We explain that 
Picture \ref{Pic:GQDM} implies the 
deformation equivalence of quantum cohomology. 
In a local frame of $F$, the connection operator 
$\nabla_X$ with $X\in T\cM$ can be written as  
\[
\nabla_X = X + \frac{1}{z}\cA_X(\tau,z).  
\]
The residual part $\cA_X(\tau,0) = [z \nabla_X]|_{z=0}$ 
defines a well-defined endomorphism 
of $F|_{\cM\times \{0\}}$. 
The flatness of the connection $\nabla$ 
implies the commutativity of these operators 
$[\cA_X(\tau,0), \cA_Y(\tau,0)]=0$. 
Note that on $V_i\subset \cM$, $\cA_X(\tau,0)$ is 
identified with the quantum product $X\circ_\tau$.    
(Here we identify the tangent vector $X$ 
with an element of $H_{\rm CR}^*(\cX_i)$.) 
We call that $(F,\nabla)$ is \emph{miniversal} 
at a point $\tau\in \cM$ if there exists a 
vector $v\in F_{(\tau,0)}$ such that the map 
\begin{equation}
\label{eq:KS}
T_\tau \cM \to F_{(\tau,0)}, \quad X \mapsto \cA_X(\tau,0) v 
\end{equation} 
is an isomorphism. 
This property clearly holds at $\tau\in V_i$ 
since we can choose $v$ to be the unit 
$\unit \in H_{\rm CR}^*(\cX)$. 
The miniversality may fail along a complex analytic 
subvariety of $\cM$. 
In the sequel, by deleting such locus if necessary,  
we assume that $(F,\nabla)$ is miniversal 
everywhere on $\cM$. 
Then we can define the product $\circ_\tau$ 
on the tangent space $T_\tau\cM$ by the formula: 
\[
\cA_{X\circ_\tau Y} (\tau,0) v = \cA_X(\tau,0) (\cA_Y(\tau,0) v), 
\]
where $v\in F_{(\tau,0)}$ is a vector which makes 
the map (\ref{eq:KS}) an isomorphism. 
The unit vector $e\in T_\tau\cM$ is defined by 
\[
\cA_e(\tau,0) v = v. 
\]
Then $(T_\tau M, \circ_\tau, e)$ becomes an 
associative commutative ring by the commutativity 
of $\cA_X(\tau,0)$. 
This definition does not depend on the choice of $v$. 
In fact, the inclusion 
\[
T_\tau \cM \hookrightarrow \End(F_{(\tau,0)}), \quad 
X \mapsto \cA_X(\tau,0) 
\]
becomes a homomorphism of rings. 
This product $\circ_\tau$ endows the base space $\cM$ 
with the structure of an $F$-manifold \cite{Her-Man}. 

The $F$-manifold $\cM$ here 
admits the \emph{Euler vector field}.  
In a local frame of $F$, 
we can write the connection 
in the $z$-direction as 
\begin{equation}
\label{eq:nabla_z_localframe}
\nabla_{z \partial_z} = z \partial_z  
- \frac{1}{z} \cU(\tau) + \cV(\tau,z), \quad 
\text{$\cV(\tau,z)$ is regular at $z=0$}. 
\end{equation} 
The residual part 
$\cU(\tau) = [z^2 \nabla_{\partial_z}]|_{z=0}$ 
again defines a well-defined endomorphism of 
the bundle $F|_{\cM\times \{0\}}$. 
The flatness of $\nabla$ implies that the endomorphism 
$\cU(\tau)$ commutes with $\cA_X(\tau,0)$ for every $X\in T\cM$. 
From this (and miniversality) it follows that 
there exists a unique vector field $E\in \Gamma(\cM,T\cM)$ 
such that 
\[
\cU(\tau) = \cA_{E}(\tau,0).  
\]
This satisfies the axiom of the Euler vector field: 
\begin{equation}
\label{eq:Euleraxiom} 
[E, X\circ_\tau Y] = [E,X]\circ_\tau Y + 
X \circ_\tau [E,Y] + X \circ_\tau Y. 
\end{equation} 

\begin{proposition}
\label{prop:Fmanifold} 
Under the Picture \ref{Pic:GQDM}, 
the quantum cohomology rings of $\cX_1$ and $\cX_2$ 
are deformation equivalent. 
They underlie the same $F$-manifold $\cM$ 
with the Euler vector field $E$. 
\end{proposition} 

\subsection{Semi-infinite variation of Hodge structures} 
\label{subsec:seminf}
The deformation equivalence explained in 
the previous section is a rather weak relationship. 
The global quantum $D$-module $F$ has 
much more information than just a family of algebras. 
We consider the \emph{semi-infinite variation of Hodge 
structures} or \seminf VHS associated to $F$.  
This notion was introduced by Barannikov \cite{Bar}. 
The information of \seminf VHS is in fact 
equivalent to that of the meromorphic 
flat connection $(F,\nabla,(\cdot,\cdot)_F)$, 
but the analogy with the ordinary Hodge theory 
may be clearer in this language.

We will work over the universal cover $\tcM$ 
of $\cM$. Let $\cH$ be the space 
of flat sections of $F$ over $\tcM\times \C^*$: 
\[
\cH := \{ s\in \Gamma(\tcM\times \C^*, \cO(F)) \;;\; 
\nabla_X s = 0, \ \forall X\in T\cM\}.
\]
Note that $s\in \cH$ is flat only
in the direction of $\cM$ and 
can be arbitrary in the $z$ direction. 
This is infinite dimensional over $\C$. 
For $\tau\in \tcM$, every section 
$s(\tau,\cdot)\in \Gamma(\{\tau\}\times \C^*, F)$ 
can be uniquely extended to a flat section 
over $\tcM\times \C^*$. 
Therefore $\cH$ is isomorphic to 
$\Gamma(\{\tau\}\times \C^*,F)$ and 
is a free $\cO(\C^*)$-module of rank $N$, 
where $\cO(\C^*)$ is the space of holomorphic functions 
on $\C^*$ and $N$ is the rank of $F$. 
The pairing on $\cH$ is defined by 
\[
(s_1, s_2)_{\cH} := (s_1(\tau,-z), s_2(\tau,z))_F \in \cO(\C^*).   
\]
Note that the right-hand side does not depend 
on $\tau$ since $s_1,s_2$ are flat in the $\cM$-direction. 
This pairing satisfies $(s_2,s_1)_{\cH} = (-)^*(s_1,s_2)_{\cH}$. 
For $\tau\in \tcM$, the space of sections of 
$F$ over $\{\tau\}\times \C$ 
is naturally embedded into $\cH$ 
(via the $\nabla$-flat extension of sections): 
\[
\Gamma(\{\tau\}\times \C, F) 
\hookrightarrow \cH.  
\] 
We denote by $\F_\tau$ the image of this embedding. 
Recall that the image of $\Gamma(\{\tau\}\times \C^*,F)$ 
gives the whole space $\cH$. $\F_\tau$ consists of 
flat sections $s\in \cH$ such that 
$s(\tau,\cdot)$ is regular at $z=0$.  
We call $\F_\tau$ the \emph{semi-infinite Hodge structure}. 
$\F_\tau$ is a free $\cO(\C)$-submodule of $\cH$ and 
can be regarded as a point on the 
Segal-Wilson Grassmannian \cite{Pre-Seg} of $\cH$ 
as follows: 
Fix an $\cO(\C^*)$-basis $e_1,\dots,e_N$ of $\cH$.  
An $\cO(\C)$-basis $s_1,\dots,s_N$ of $\F_\tau$ 
can be written as 
$s_j = \sum_{i=1}^{N} e_i c_{ij}(\tau,z)$. 
By restricting $z$ to lie on $S^1$, 
the $N\times N$ matrix $(c_{ij}(\tau,z))$ defines 
an element of the loop group $LGL(N,\C)$.  
A change of the basis $s_j$ changes 
the matrix $(c_{ij})$ by the left multiplication by 
an element of the positive loop group $LGL^+(N,\C)$ 
(whose entries are holomorphic functions on $\C$). 
Thus the subspace $\F_\tau$ is identified with 
an element $[(c_{ij}(\tau,z))]$ of 
$LGL(N,\C)/LGL^+(N,\C)=:\Gr_{\frac{\infty}{2}}(\cH)$. 
We call the map 
\[
\tcM \ni \tau \longmapsto \F_\tau \in 
\Gr_{\frac{\infty}{2}}(\cH)  
\]
the \emph{semi-infinite period map}. 

\begin{proposition}[{\cite[Proposition 2.9]{CIT}}]  
\label{prop:seminf}
The semi-infinite period map $\tau\mapsto \F_\tau$ 
satisfies the following: 
\begin{align*} 
&X \F_\tau \subset z^{-1}\F_\tau, \quad X\in T_\tau\cM, \\
&(\F_\tau,\F_\tau)_{\cH} \subset \cO(\C), \\ 
&(\nabla_{z\partial_z} +E) \F_\tau  \subset \F_\tau,   
\end{align*} 
where we used the fact that 
$\nabla_{z\partial_z}$ acts on $\cH$ 
as a $\C$-endomorphism.  
The first property is an analogue of  
Griffiths transversality and 
the second is the Hodge-Riemann bilinear relation. 
\end{proposition}

\subsection{Opposite subspace and Frobenius manifolds} 
\label{subsec:opp_Frob} 
As we remarked, the global quantum $D$-module 
is not a priori trivialized. 
A good trivialization is given by the choice of 
an \emph{opposite subspace} to the \seminf VHS. 
The choice of an opposite subspace and a dilaton shift 
defines a Frobenius structure on 
the universal cover of $\cM$. 
The Frobenius/flat structure was discovered by 
K. Saito \cite{SaitoK} as a structure on a miniversal 
deformation of isolated hypersurface singularities   
and the use of opposite subspaces 
goes back to M. Saito's work \cite{SaitoM} 
in that context. 
Let $\cO(\Proj^1\setminus\{0\})$ be the space 
of holomorphic functions on $\Proj^1\setminus\{0\}$. 
This is contained in $\cO(\C^*)$. 

\begin{definition}
An \emph{opposite subspace} $\cH_-$ at $\tau\in \tcM$ 
is a free $\cO(\Proj^1\setminus\{0\})$-submodule of 
$\cH$ such that the natural map 
\begin{equation}
\label{eq:opposite}
\cH_- \oplus \F_\tau \to \cH
\end{equation} 
is an isomorphism. 
$\cH_-$ is said to be \emph{homogeneous} if 
\[
\nabla_{z\partial_z} \cH_- \subset \cH_- 
\]
and \emph{isotropic} if 
\[
(\cH_-,\cH_-)_\cH \subset z^{-2} \cO(\Proj^1\setminus\{0\}).  
\]
\end{definition} 

In terms of the loop Grassmannian $LGL(N,\C)/LGL^+(N,\C)$, 
$\cH_-$ is opposite 
at $\tau$ if $\F_\tau$ lies on the ``big cell": 
an open orbit of $LGL^-(N,\C)$. 
Therefore, the opposite property 
((\ref{eq:opposite}) is an isomorphism) 
is an open condition:  
If $\cH_-$ is opposite at $\tau$, then 
it is opposite in a neighborhood of $\tau$. 
Given an opposite subspace $\cH_-$ at some point, 
the opposite property may fail 
along a complex analytic subvariety of $\tcM$.

We explain that a homogeneous 
opposite subspace corresponds to 
an extension of $(F,\nabla)$ across $z=\infty$ 
such that the connection $\nabla$ 
has a logarithmic singularity along $z=\infty$. 

\begin{lemma}
\label{lem:opposite} 
For a point $\tau\in \tcM$, 
the following are equivalent: 

{\rm (i)} $\cH_-$ is a homogeneous opposite subspace 
at $\tau$. 

{\rm (ii)} $\cH_-$ is homogeneous and 
one of the natural maps  
\[
\begin{CD}
z \cH_-/\cH_- @<<< z \cH_- \cap \F_\tau @>>> \F_\tau/z\F_\tau  
\end{CD} 
\]
is an isomorphism of finite dimensional $\C$-vector spaces. 

{\rm (iii)} 
Define an extension $\hF_\tau \to \{\tau\}\times \Proj^1$ 
of the vector bundle $F|_{\{\tau\}\times \C}$ to 
$\{\tau\}\times \Proj^1$ as follows: 
We define a section $s\in \Gamma(\{\tau\}\times \C^*, F)$ 
to be regular at $z=\infty$ 
if the image of $s$ in $\cH$ lies in $z\cH_-$. 
Then the extension $(\hF_\tau,\nabla)$ is a trivial 
vector bundle over $\Proj^1$ and 
$\nabla$ has a logarithmic singularity 
at $z=\infty$. 
\end{lemma} 
\begin{proof}
(i) $\Rightarrow$ (ii). The injectivity of the maps in (ii) 
is obvious. For $[v] \in z\cH_-/\cH_-$ with $v\in z\cH_-$, 
write $v=v_0 + v_-$ where $v_0\in \F_\tau$ and 
$v_0 \in \cH_-$. Then $v_0 = v - v_- \in z \cH_-\cap \F_\tau$ 
and $[v] = [v_0]$. This shows the surjectivity 
of $z\cH_-\cap \F_\tau \to z\cH_-/\cH_-$. 
For $[v]\in \F_\tau/z\F_\tau$ with $v\in \F_\tau$, 
write $z^{-1}v = v_- + v_0$, where $v_-\in \cH_-$ 
and $v_0\in \F_\tau$. Then 
$z v_- = v - zv_0 \in \F_\tau \cap z \cH_-$ 
and $[v]=[zv_-]$. This shows the surjectivity 
of $z\cH_-\cap \F_\tau \to \F_\tau/z\F_\tau$. 

(ii) $\Rightarrow$ (iii). 
Consider the extension $\hF_\tau \to \{\tau\}\times \Proj^1$ 
in (iii). 
We can identify $z\cH_-/\cH_-$ with 
the fiber $\hF_{(\tau,\infty)}$, 
$z\cH_-\cap \F_\tau$ with the global section 
$\Gamma(\Proj^1,\hF_\tau)$ 
and $\F_\tau/z\F_\tau$ with 
the fiber $\hF_{(\tau,0)}$.  
Since the maps in (ii) are induced from 
the restrictions, that one of them 
is an isomorphism implies that $\hF_\tau$ 
is a trivial holomorphic vector bundle. 
For a local co-ordinate $w=z^{-1}$ around $z=\infty$, 
we have $\nabla_{w\partial_w} = - \nabla_{z\partial_z}$. 
Hence the homogeneity implies  
$\nabla_{w\partial_w} (z\cH_-) \subset (z\cH_-)$, 
so $\nabla$ has a logarithmic singularity at $w=0$. 

(iii) $\Rightarrow$ (i). 
Note that $\cH$ is identified with 
the space of sections of $\hF_\tau$ over $\{\tau\}\times \C^*$. 
Because $\hF_\tau$ is trivial, 
that (\ref{eq:opposite}) is an isomorphism 
follows from the decomposition 
\[
\cO(\C^*) =z^{-1} \cO(\Proj^1\setminus\{0\}) \oplus \cO(\C). 
\]
The logarithmic singularity of $\nabla$ implies 
the homogeneity of $\cH_-$. 
\end{proof} 

By the isomorphism in (ii) of Lemma \ref{lem:opposite},  
a homogeneous opposite subspace $\cH_-$ gives 
a local trivialization of $F$. 
In fact, since $F|_{\{\tau\}\times \C}$ extends 
to a trivial vector bundle $\hF_\tau$ over 
$\{\tau\}\times \Proj^1$, 
we have 
\begin{equation}
\label{eq:triv_opposite}
F_{(\tau,z)} \cong 
\Gamma(\{\tau\}\times \Proj^1,\hF_\tau) 
\cong z \cH_- \cap \F_\tau 
\cong z \cH_-/\cH_-. 
\end{equation} 
The finite dimensional vector space 
$z\cH_-/\cH_-$ does not depend on $\tau$,  
so this defines a trivialization of $F$ 
over an open subset of $\tcM$. 
Under this trivialization, 
the flat connection $\nabla$ 
can be written as follows: 
\begin{align}
\label{eq:conn_under_opposite} 
\begin{split}
&\nabla_X = X + \frac{1}{z} \cA_X(\tau), \quad 
X\in T\cM, \\
&\nabla_{z\partial_z} = z\partial_z - 
\frac{1}{z} \cU(\tau) + \cV,  
\end{split}
\end{align} 
where $\cA(\tau)$ is an 
$\End(z\cH_-/\cH_-)$-valued 1-form, 
$\cU(\tau)$ is an $\End(z\cH_-/\cH_-)$-valued function,  
and $\cV$ is a constant operator in $\End(z\cH_-/\cH_-)$. 
Here $\cA(\tau),\cU(\tau)$ are independent of $z$ 
and defined on an open subset of $\tcM$.  
Note that $\cU(\tau) = \cA_E(\tau)$ by 
the definition of the Euler vector field $E$. 

\begin{figure}[htbp] 
\begin{center}
\begin{picture}(200,60) 
\put(85,0){\line(0,1){60}}
\put(75,55){\makebox(0,0){$\cH_-$}}
\put(185,50){\makebox(0,0){$\F_\tau$}} 
\put(125,10){\makebox(0,0){$\bullet$}} 
\put(135,10){\makebox(0,0){$\hat{v}_0$}}
\put(125,0){\line(0,1){60}} 
\put(127,65){\makebox(0,0){$\hat{v}_0 + \cH_-$}} 
\put(125,33){\makebox(0,0){$\bullet$}} 
\put(152,30){\makebox(0,0){$\hat{v}_0+\psi(\tau)$}} 
\put(85,20){\makebox(0,0){$\bullet$}} 
\put(77,25){\makebox(0,0){$0$}} 
\thicklines
\put(25,0){\line(3,1){150}}
\end{picture} 
\end{center} 
\caption{$J$-function $\hat{v}_0+\psi(\tau)$ and 
flat co-ordinates $[\psi(\tau)]\in \cH_-/z^{-1}\cH_-$.} 
\label{fig:flat} 
\end{figure}
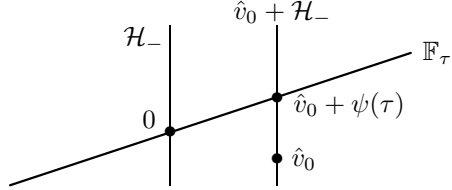

In order to have a Frobenius structure on $\tcM$,  
in addition to $\cH_-$, we need to choose 
an eigenvector $v_0\in z\cH_-/\cH_-$ of $\cV$
satisfying the miniversality condition\footnote{
The action of $\cV$ on $z\cH_-/\cH_-$ is induced from 
that of $\nabla_{z\partial_z}$ on $z\cH_-$.}: 
\begin{equation}
\label{eq:primitivesection} 
T_\tau \tcM \to z\cH_-/\cH_-, \quad 
X \mapsto \cA_X(\tau) v_0 \quad 
\text{is an isomorphism.} 
\end{equation} 
We call $v_0$ the \emph{dilaton shift}. 
The isomorphism $T_\tau \tcM \cong z\cH_-/\cH_-$ above 
defines a \emph{flat structure} on $\tcM$. 
A vector field $X$ is defined to be flat 
if $\cA_X(\tau)v_0$ is a constant element in $z\cH_-/\cH_-$. 
This flat structure is integrable. 
Let $\hat{v}_0 + \psi(\tau)$ be the unique intersection 
point of $\F_\tau$ and the affine subspace $\hat{v}_0 + \cH_-$, 
where $\hat{v}_0\in z\cH_-$ is an 
(arbitrarily fixed) lift of $v_0$ and $\psi(\tau)\in \cH_-$. 
See Figure \ref{fig:flat}. 
Then the map 
\[
\tcM\ni \tau\mapsto [\psi(\tau)]\in \cH_-/z^{-1}\cH_-
\] 
is a local isomorphism and gives a flat co-ordinate system. 
In fact, the differential of this map is identified 
with (\ref{eq:primitivesection}). 
Varying $\tau$, the intersection point 
$\hat{v}_0 + \psi(\tau) \in \F_\tau$ 
gives a section $s_0$ of $F$ 
which corresponds to $v_0\in z\cH_-/\cH_-$ 
in the trivialization (\ref{eq:triv_opposite}). 
(Note that $\hat{v}_0+\psi(\tau) \in z\cH_-\cap\F_\tau$.) 
This section $s_0$ is called a \emph{primitive section}. 
In Gromov-Witten theory, the corresponding vector 
$\hat{v}_0 + \psi(\tau)\in \cH$ is called 
the \emph{$J$-function}.

For a flat vector field $X$, 
we have $\cV(\cA_X v_0)= \cA_{(\alpha+1)X - [X,E]} v_0$ 
where $\alpha$ is the eigenvalue of $v_0$ 
with respect to $\cV$. 

When $\cH_-$ is isotropic, the pairing $(\cdot,\cdot)_\cH$ 
on $\cH$ induces a symmetric bilinear 
$\C$-valued pairing on $z\cH_-\cap \F_\tau \cong z\cH_-/\cH_-$. 
By pulling back this pairing on $z\cH_-/\cH_-$ 
to $T_\tau\tcM$ by the map (\ref{eq:primitivesection}), 
we obtain a $\C$-bilinear metric 
$g\colon T_\tau\tcM\times T_\tau \tcM \to \C$. 
The metric tensor of $g$ is constant in the 
flat co-ordinates above, so the metric $g$ is flat.

\begin{proposition}[{\cite[Proposition 2.12]{CIT}}]
\label{prop:Frobenius} 
Take an isotropic homogeneous opposite subspace 
$\cH_-$ and a dilaton shift $v_0\in z\cH_-/\cH_-$ 
satisfying (\ref{eq:primitivesection}) at some point $\tau$. 
Then the $F$-manifold structure $(\circ_\tau,e,E)$ 
in Proposition \ref{prop:Fmanifold} 
is lifted to the Frobenius manifold structure 
$(\circ_\tau,e,E,g)$ 
on the complement of a complex analytic 
subvariety in $\tcM$. 
These data satisfy:  

{\rm (i)} the Levi-Civita connection $\nabla^{\rm LC}$ 
of $g$ is flat;  

{\rm (ii)} $(T_\tau\cM,\circ_\tau,g)$ is a commutative 
Frobenius algebra; 

{\rm (iii)} the pencil of flat connections 
$\nabla^\lambda_X= \nabla^{\rm LC}_X + \lambda X\circ_\tau$ 
is flat;  

{\rm (iv)} the unit vector field $e$ is flat;  

{\rm (v)} the Euler vector field $E$ 
satisfies (\ref{eq:Euleraxiom}), 
$(\nabla^{\rm LC})^2 E =0$ and 
\begin{align*} 
E g(X,Y) = g([E,X],Y) + g(X,[E,Y]) + 2(\alpha+1) g(X,Y),  
\end{align*} 
where $\alpha\in \C$ is the eigenvalue of $v_0$: 
$\cV v_0 = \alpha v_0$. 
\end{proposition} 

\subsection{Opposite subspaces at cusps} 
\label{subsec:opp_cusp}
We regard the large radius limit point of $\cX_i$ 
as a \emph{cusp} of the global K\"{a}hler moduli space $\cM$
and $V_i$ as its neighborhood. 
Since the base space of $QDM(\cX_i)$ 
is a quotient of a vector space, 
$V_i$ is equipped with the standard 
Frobenius/flat structure as described in 
\cite{Manin, Dub:2D}.  
We will show that, under certain conditions, 
the Frobenius structure (or the corresponding opposite subspace)  
of $V_i$ can be uniquely characterized 
by the monodromy invariance and 
the compatibility with the Deligne extension.  
This means that \emph{there is a canonical choice of 
the Frobenius manifold structure at each cusp 
from a purely $D$-module theoretic viewpoint.} 
The characterization here 
was shown in the case $\cX=\Proj(1,1,1,3)$ in \cite{CIT}. 

Henceforth we study  
the global quantum $D$-module restricted to 
$V_i$ \emph{i.e.} $QDM(\cX_i)$. 
We omit the subscript $i$ 
and write $V,\cX$ for $V_i, \cX_i$ etc.  
The open set $U\subset H_{\rm CR}^*(\cX)$ in 
Assumption \ref{as:convergence} 
is identified with the universal cover of $V\cong U/H^2(\cX,\Z)$. 

\begin{definition}[Givental space \cite{CG,Giv:quant}] 
The \emph{Givental symplectic space} $\cH^{\cX}$ is defined 
to be a free $\cO(\C^*)$-module 
\[
\cH^{\cX} := H^*_{\rm CR}(\cX) \otimes \cO(\C^*), 
\]
endowed with an $\cO(\C^*)$-valued pairing 
$(\cdot,\cdot)_{\cH}$: 
\[
(f(z),g(z))_{\cH} = (f(-z),g(z))_{\rm orb}. 
\]
As an infinite dimensional vector space over $\C$, 
$\cH^\cX$ has the following symplectic form: 
\begin{equation}
\label{eq:symp_form_Giv}
\Omega(f,g) = \Res_{z=0} (f(-z), g(z))_{\rm orb} dz. 
\end{equation} 
We identify the Givental space $\cH^\cX$ 
with the space $\cH$ of flat sections of $QDM(\cX)$ 
over $U$ through the fundamental solution in 
Proposition \ref{prop:fundamentalsol}: 
\[
\cH^{\cX} \cong \cH, \quad 
\phi(z) \mapsto L(\tau,z) \phi(z).  
\]
This identification preserves the pairing.  
\end{definition} 

In terms of the Givental space, 
the semi-infinite Hodge structure 
$\F_\tau$ is identified with 
the Lagrangian subspace: 
\begin{equation}
\label{eq:Ftau}
\F_\tau = L(\tau,z)^{-1}(H^*_{\rm CR}(\cX)\otimes \cO(\C) ) 
\subset \cH^{\cX}, \quad \tau \in U. 
\end{equation} 
The Givental space has a standard opposite 
subspace $\cH_-^{\cX}$:  
\[
\cH_-^{\cX} := z^{-1} 
H^*_{\rm CR}(\cX) \otimes \cO(\Proj^1\setminus\{0\}) 
\subset \cH^{\cX}.
\]
In fact, this is opposite to $\F_\tau$ 
(\emph{i.e.} $\cH_-^{\cX} \oplus \F_\tau = \cH^{\cX}$)  
for every $\tau\in U$ 
because $L(\tau,z)$ is regular at $z=\infty$ 
and $L(\tau,z) = \id +O(z^{-1})$.  
\begin{proposition}
The standard opposite subspace $\cH_-^{\cX}$ 
is homogeneous and isotropic. 
This $\cH_-^{\cX}$ and the standard dilaton shift 
$v_0 = \unit \in z\cH_-^{\cX}/\cH_-^{\cX}$ 
endow the base space 
$V \cong U/H^2(\cX,\Z)$ 
of the quantum $D$-module 
with the standard Frobenius manifold structure 
coming from the linear structure on 
$U\subset H^*_{\rm CR}(\cX)$ and 
the orbifold Poincar\'{e} pairing on 
$T_\tau U \cong H^*_{\rm CR}(\cX)$. 
(See Proposition \ref{prop:Frobenius} for 
the construction of Frobenius manifolds.)
\end{proposition} 
\begin{proof}
It follows from Proposition \ref{prop:fundamentalsol} 
that $L(\tau,z)$ satisfies the differential equation 
$
\nabla_{z\partial_z} L(\tau,z) \phi = 
L(\tau,z) (\mu - \rho/z) \phi$ 
for $\phi \in H_{\rm CR}^*(\cX)$. 
This shows that the action of 
$\nabla_{z\partial_z}$ on the Givental space is 
given by 
\begin{equation}
\label{eq:nabla_z_H}
\nabla_{z\partial_z} = z\partial_z + \mu - \frac{\rho}{z} 
\quad \text{on $\cH^{\cX}$}. 
\end{equation} 
Therefore the standard opposite subspace is homogeneous 
$\nabla_{z\partial_z} \cH^{\cX}_- \subset \cH^{\cX}_-$. 
It is obvious that $\cH_-^{\cX}$ is isotropic. 
Because $L(\tau,z)^{-1}\phi = \phi + O(z^{-1})$ for 
$\phi\in H^*_{\rm CR}(\cX)$, we have 
$L(\tau,z)^{-1}\phi \in z\cH_-^{\cX} \cap \F_\tau$. 
Therefore, the constant section $\phi$ of $QDM(\cX)$ 
corresponds to (again) the constant element 
$\phi\in z\cH_-^{\cX}/\cH_-^{\cX}$ under 
the trivialization (\ref{eq:triv_opposite}). 
This means that $\cH_-^{\cX}$ yields    
exactly the \emph{given} trivialization of $QDM(\cX)$. 
In particular, the connection operators $\cA_X$, $\cU$, 
$\cV$ in (\ref{eq:conn_under_opposite}) 
are identified with $X\circ_\tau$, 
$E\circ_\tau$, $\mu$ and $\unit \in H^*_{\rm CR}(\cX)$
is the eigenvector of $\cV=\mu$ of eigenvalue 
$-\dim_\C\cX/2$. 
Now we only need to check that the corresponding 
flat metric $g$ is the orbifold Poincar\'{e} pairing. 
But this is obvious from 
$(L(\tau,-z)^{-1} \phi_1, L(\tau,z)^{-1}\phi_2)_{\rm orb}  
= (\phi_1,\phi_2)_{\rm orb}$. 
\end{proof} 

\noindent
\emph{The monodromy invariance of $\cH_-^{\cX}$}:  
We see that $\cH_-^{\cX}$ 
is invariant under the local monodromy 
(or Galois actions) around the large radius limit. 
The Galois action in Lemma \ref{lem:Galois} 
acts on the Givental space $\cH^{\cX}$ by $G^{\cH}(\xi)$: 
\[
G^{\cH}(\xi) = \bigoplus_{v\in \sfT}
e^{-2\pi\iu \xi_0/z} e^{2\pi\iu f_v(\xi)}, 
\quad 
\xi\in H^2(\cX,\Z),   
\]
where we used the decomposition 
$\cH^{\cX} = \bigoplus_{v\in \sfT} H^*(\cX_v)\otimes \cO(\C^*)$. 
Since $G^{\cH}(\xi)$ contains only negative powers in $z$, 
we have 
\begin{equation}
\label{eq:monodromy-invariance_opposite}
G^{\cH}(\xi)\cH_-^\cX \subset \cH_-^\cX. 
\end{equation} 
The Hodge structures are monodromy-equivariant:   
$G^{\cH}(\xi)\F_\tau = \F_{G(\xi)\tau}$. 
The monodromy-invariance of $\cH_-^\cX$ 
corresponds to the fact that 
the corresponding Frobenius manifold structure 
is well-defined\footnote{More precisely, 
we also need the fact that the vector 
$\unit\in z\cH_-^\cX/\cH_-^\cX$ is invariant 
under the Galois action.} 
on the quotient $V\cong U/H^2(\cX,\Z)$. 
The induced action of $G^{\cH}(\xi)$ 
on $z\cH_-^\cX/\cH_-^\cX$ 
is given by $\bigoplus_{v\in \sfT} e^{2\pi\iu f_v(\xi)}$. 
Because $f_v(\xi)$ is a rational number, 
there exists a positive integer $k_0>0$ such that 
\begin{equation}
\label{eq:power-identity}
(G^{\cH}(\xi))^{k_0} = \id \quad 
\text{on $z\cH_-^\cX/\cH_-^\cX$}. 
\end{equation} 
This corresponds to the fact that 
the monodromy of the Levi-Civita connection $\nabla^{\rm LC}$ 
of the flat metric $g$ (or the monodromy of the 
trivialization (\ref{eq:triv_opposite}))  
becomes trivial on a $k_0$-fold cover of $V$. 
In fact, one can see that 
the monodromy of $\nabla^{\rm LC}$ 
is trivial on the cover 
$U/H^2(X,\Z)\to U/H^2(\cX,\Z)\cong V$, 
where $X$ is the coarse moduli space of $\cX$.  

\noindent
\emph{Compatibility with Deligne's extension}: 
As we did at the end of Section \ref{subsec:QDM}, 
we can extend the quantum $D$-module on the cover $U/H^2(X,\Z)$ 
to a connection on $\ov{U/H^2(X,\Z)}$ 
with a logarithmic pole along $q^1\cdots q^r =0$ 
by choosing a nef basis $p_1,\dots,p_r$ of $H^2(X,\Z)/{\rm tors}$. 
This is a Deligne extension 
of $\nabla$ for a fixed $z\in \C^*$. 
A Deligne extension is given by 
the choice of a logarithm of the monodromy 
$M_a := G^{\cH}(p_a)=e^{-2\pi\iu p_a/z}$ 
around the axis $q^a=0$.  
In our case, we have the ``standard" logarithm 
$\Log(M_a) = -2\pi\iu p_a/z$ 
since $M_a$ is unipotent.  
Our Deligne extension can be described as follows. 
A section $s(\tau,z)$ of $F$ over $(U/H^2(X,\Z))\times \C^*$ 
is extendible to $\ov{U/H^2(X,\Z)}\times \C^*$ 
if the image $\iota_\tau(s)\in \cH^\cX$ of 
$s(\tau,\cdot)\in \Gamma(\{\tau\}\times \C^*,F)$ 
satisfies the following: 
\emph{the family of elements in $\cH^\cX$} 
\[ 
U/H^2(X,\Z)\ni 
[\tau] \mapsto 
\tilde{s}_\tau:=
\exp\left(\sum_{a=1}^r \frac{\log q^a}{2\pi\iu} \Log(M_a)
\right) \iota_\tau(s) 
\in \cH^\cX 
\]
\emph{extends holomorphically to $\ov{U/H^2(X,\Z)}$}, 
where we put $\tau=\tau_{0,2} +\tau'$ as in 
(\ref{eq:decomp_tau}) and $\tau_{0,2} = \sum_{a=1}^r p_a \log q^a$. 
Note that $\tilde{s}_\tau$ is single-valued on $U/H^2(X,\Z)$ 
since the exponential factor offsets the monodromy. 
Moreover, the limit of $s(\tau,z)$ at $q=\tau'=0$ is 
regular at $z=0$ if $\tilde{s}_\tau|_{q=\tau'=0}$ 
lies in the \emph{limiting Hodge structure} 
$\Flim$: 
\[
\Flim := \lim_{\substack{q\to 0 \\ \tau'\to 0}} 
\exp\left(\sum_{a=1}^r \frac{\log q^a}{2\pi\iu}
\Log(M_a)\right) \F_\tau,    
\]
where we put $\tau = \sum_{a=1}^r p_a \log q^a+\tau'$ 
as in (\ref{eq:decomp_tau}).    
By using (\ref{eq:Ftau})  
and the definition (\ref{eq:fundamentalsol}) of $L(\tau,z)$, 
one can check that $\Flim$ exists and 
\begin{equation}
\label{eq:limitingHodge}
\Flim = H_{\rm CR}^*(\cX) \otimes \cO(\C) \subset \cH^\cX.  
\end{equation} 
The existence of $\Flim$ 
is an analogue of the \emph{nilpotent orbit theorem}
\cite{Schm} in quantum cohomology. 
This means that the Hodge structure 
$\F_\tau$ is approximated by the nilpotent orbit 
$e^{-\sum_{a=1}^r \log q^a \Log(M_a)/(2\pi\iu)}\Flim$ 
as $q,\tau'\to 0$. 
The standard opposite subspace is opposite to $\Flim$: 
\begin{equation}
\label{eq:limiting_opposite}
\cH_-^{\cX} \oplus \Flim = \cH^\cX. 
\end{equation} 
This corresponds to the fact that the 
trivialization induced from $\cH_-^\cX$ 
is compatible with the Deligne extension 
at $q=0$, \emph{i.e.} 
a section which is constant 
in the trivialization (\ref{eq:triv_opposite}) 
is extendible across $q=0$ in the Deligne extension. 
Note that this is a stronger condition   
than that $\cH_-$ is opposite to $\F_\tau$ 
for every $\tau\in U$.  

For a multiplicative character 
$\alpha\colon H^2(\cX,\Z)\to \C^*$, we put 
\[
\sfT_\alpha := 
\{v\in \sfT\;;\; \exp(2\pi\iu f_v(\xi)) = \alpha(\xi), 
\forall \xi \in H^2(\cX,\Z)\}.  
\]
Because $e^{2\pi\iu \iota_v} = \alpha([-K_\cX])$ for 
$v\in \sfT_\alpha$, the age $\iota_v$ 
for $v\in \sfT_\alpha$ have the common fractional part  
for each $\alpha$. 
Consider the following two conditions.  
\begin{align} 
\label{eq:uniqueness_opposite} 
&\forall \alpha, \ \exists n_\alpha\in \Q  
\text{ such that } 
\forall v\in \sfT_\alpha \ 
(n_v + 2 \iota_v = n_\alpha \text{ or } 
\ n_\alpha+1).   
\\
\label{eq:uniqueness_dilaton}
&\iota_v = 0, \ v\neq 0 \ 
\Longrightarrow \ 
\exists \xi \in H^2(\cX,\Z)  
\text{ such that } f_v(\xi)>0.   
\end{align} 
Here $n_v := \dim_\C\cX_v$.  
The first condition is 
a rather weaker version of 
the Hard Lefschetz condition 
we will see later\footnote{
The condition (\ref{eq:uniqueness_opposite}) says that  
$V_\alpha = \bigoplus_{v\in \sfT_\alpha}H^{*-2\iota_v}(\cX_v)$ 
is bicentric HL in the sense of Definition 
\ref{def:gen_HL}. See also Remark \ref{rem:genHL_coarse}.}.  
(There we have $n_v+2\iota_v = \dim_\C \cX$ for all $v$.) 
When (\ref{eq:uniqueness_opposite}) is satisfied, we put 
\begin{equation}
\label{eq:decomp_Talpha}
\sfT_{\alpha,j} = \{v\in \sfT_\alpha\;;\; 
n_v + 2\iota_v = n_\alpha + j\}, \quad 
\sfT_\alpha =  \sfT_{\alpha,0} \sqcup \sfT_{\alpha,1}. 
\end{equation}

% \vspace{5pt}
% \noindent
% \emph{Example}:
\begin{example}
If $\cX$ is isomorphic to a quotient $[M/G]$ of a manifold $M$ 
by an \emph{abelian} Lie group $G$ as a topological orbifold,
the conditions (\ref{eq:uniqueness_opposite}), 
(\ref{eq:uniqueness_dilaton}) are satisfied since 
every $\sfT_\alpha$ consists of one element. 
In fact, there are sufficiently many line bundles on $[M/G]$ 
arising from characters of $G$ which ``separate"  
different inertia components. 
In particular, these hold for toric orbifolds. 
\end{example} 

\begin{theorem}
\label{thm:uniqueness_cusp}
Assume that the coarse moduli space $X$ of $\cX$ is projective.  
The standard opposite subspace $\cH_-=\cH^\cX_-$ 
and the standard dilaton shift $v_0=\unit$ 
are characterized as follows. 

(i) Under the condition 
(\ref{eq:uniqueness_opposite}), 
there exists a unique homogeneous opposite subspace satisfying 
the monodromy invariance (\ref{eq:monodromy-invariance_opposite}), 
(\ref{eq:power-identity}) and 
the compatibility with the Deligne extension 
(\ref{eq:limiting_opposite}).  

(ii) Under the condition 
(\ref{eq:uniqueness_dilaton}), there exists 
a unique vector $v_0\in z\cH_-^\cX/\cH_-^\cX$ 
(up to a scalar multiple) such that 
$v_0$ is an eigenvector of $\mu=\cV=[\nabla_{z\partial_z}]$ 
of the smallest eigenvalue $-\dim_\C\cX/2$ 
and invariant under every Galois action on 
$z\cH^\cX_-/\cH^\cX_-$. 

Thus under (\ref{eq:uniqueness_opposite}) and 
(\ref{eq:uniqueness_dilaton}), 
the above conditions determine  
a canonical Frobenius structure at the cusp 
up to a constant multiple of the flat metric. 

% This means, under the conditions 
% (\ref{eq:uniqueness_opposite}), (\ref{eq:uniqueness_dilaton}), 
% there exists only one ``reasonable" 
% Frobenius manifold structure at each cusp  
% determined up to a constant multiple of the flat metric $g$. 
\end{theorem}  
\begin{proof}
Let $\cH_-\subset \cH^\cX$ 
be any homogeneous opposite subspace 
satisfying (\ref{eq:monodromy-invariance_opposite}), 
(\ref{eq:power-identity}) and (\ref{eq:limiting_opposite}). 
We decompose the Galois action as 
\[
G^{\cH}(\xi) = e^{-2\pi\iu \xi_0/z} \circ G^{\cH}_0(\xi), 
\quad G^{\cH}_0(\xi) = \bigoplus_{v\in \sfT}
e^{2\pi\iu f_v(\xi)}. 
\] 

\noindent
\emph{Claim}: 
$\cH_-$ satisfies the following: 
\[
\xi_0 \cdot \cH_- \subset \cH_-, \quad 
G^{\cH}_0(\xi) \cH_- \subset \cH_-, \quad 
(z\partial_z +\mu) \cH_- \subset \cH_-. 
\]
Take a sufficiently big $k_0>0$ such that 
$(G^{\cH}_0(\xi))^{k_0}=\id$ 
and (\ref{eq:power-identity}) hold.  
Then $(G^{\cH}(\xi))^{k_0} = e^{-k_0 2\pi\iu \xi_0/z}$ 
preserves $\cH_-$ and acts trivially on $z\cH_-/\cH_-$.
Then $\Log((G^{\cH}(\xi))^{k_0}) = -k_0 2\pi\iu \xi_0/z$ 
sends $z\cH_-$ to $\cH_-$. This implies the first equation.   
The second equation follows from 
$G^\cH_0(\xi) = e^{2\pi\iu \xi_0/z} \circ G^\cH(\xi)$ 
and (\ref{eq:monodromy-invariance_opposite}).  
The third equation follows from 
$\nabla_{z\partial_z}\cH_-\subset \cH_-$, 
the formula (\ref{eq:nabla_z_H}) 
for $\nabla_{z\partial_z}$ 
and $(\rho/z) \cH_- \subset \rho \cH_- \subset \cH_-$. 

The third equation in the claim means that 
$\cH_-$ is homogeneous with respect to 
the usual grading on $H_{\rm CR}^*(\cX)$ together with 
$\deg z = 2$. 
The opposite property (\ref{eq:limiting_opposite}) 
and the formula (\ref{eq:limitingHodge}) for $\Flim$ 
imply that 
\begin{equation}
\label{eq:limiting_lift}
z\cH_-\cap \Flim \cong \Flim/z\Flim 
= H_{\rm CR}^*(\cX). 
\end{equation} 
Since $z\partial_z +\mu$ preserves 
$z\cH_- \cap \Flim$, this is an isomorphism of 
graded vector spaces.  
Also $G^{\cH}_0(\xi)$ preserves 
$z\cH_-\cap \Flim$ and (\ref{eq:limiting_lift}) 
is equivariant with respect to 
the action of $G^{\cH}_0(\xi)$. 
Therefore (\ref{eq:limiting_lift}) 
is decomposed into the sum 
of simultaneous eigenspaces of the commuting 
operators $G^{\cH}_0(\xi)$. 
Recall that the condition (\ref{eq:uniqueness_opposite}) 
gives the decomposition (\ref{eq:decomp_Talpha}). 
Take a multiplicative character 
$\alpha \colon H^2(\cX,\Z)\to \C^*$ 
and set 
\[
V_{\alpha,j} = \bigoplus_{v\in\sfT_{\alpha,j}} 
H^{*-2\iota_v}(\cX_v), \quad j=0,1, \quad 
V_\alpha = V_{\alpha,0}\oplus V_{\alpha,1}. 
\] 
Then $V_\alpha$ 
is the simultaneous eigenspace 
of $G^{\cH}_0(\xi)$ of eigenvalue $\alpha$.  
By (\ref{eq:limiting_lift}), for a homogeneous element 
$\phi\in V_{\alpha,j}$, there exists a unique lift 
$\hat{\phi}\in z\cH_-\cap \Flim$ 
such that   
\[
\hat{\phi} = \phi + O(z), \quad \deg \hat{\phi} = \deg \phi, 
\quad \hat{\phi} \in V_\alpha \otimes \cO(\C^*).  
\]
By the Claim above,  
the $H^2(\cX)$-action also preserves $z\cH_-\cap \Flim$. 
Therefore we have 
$\widehat{\omega\cdot \phi} = \omega \cdot \hat{\phi}$ 
for a K\"{a}hler class $\omega$. 
Because $X$ is K\"{a}hler, 
the cohomology ring $H^*(\cX_v)$ of every inertia component 
has the Hard Lefschetz property.  
Hence under the condition (\ref{eq:uniqueness_opposite}), 
the following holds with respect to the 
grading of the Chen-Ruan cohomology $H_{\rm CR}^*(\cX)$. 
\begin{equation}
\label{eq:HL_Valphaj}
\omega^i 
\colon V_{\alpha,j}^{n_\alpha+j-i} \to 
V_{\alpha,j}^{n_\alpha+j+i} \quad 
\text{ is an isomorphism} \quad j=0,1.   
\end{equation} 
We also have the Lefschetz decomposition of $V_{\alpha,j}$: 
\[
V_{\alpha,j} = \bigoplus_{k\ge 0} \bigoplus_{i=0}^k 
\omega^i PV_{\alpha,j}^{n_\alpha+ j -k} 
\]
where $PV_{\alpha,j}^{n_\alpha+j-k} = \Ker(\omega^{k+1}\colon 
V_{\alpha,j}^{n_\alpha+j-k} \to V_{\alpha,j}^{n_\alpha +j +k +2})$
is the primitive part. 
By the property $\widehat{\omega\cdot\phi} = 
\omega\cdot \hat{\phi}$, we only need to know 
$\hat{\phi}$ for $\phi\in PV_{\alpha,j}^{n_\alpha+j-k}$. 
For $\phi\in PV_{\alpha,j}^{n_\alpha+j-k}$, 
we can put 
\[
\hat{\phi} = \phi + z \phi_1 + z^2 \phi_2 +\cdots. 
\]
where $\phi_i \in V_\alpha^{n_\alpha+j-k-2i}$.  
Then $0=\widehat{\omega^{k+1} \phi} = 
\sum_{i\ge 1} z^i \omega^{k+1} \phi_i$. 
This implies $\omega^{k+1} \phi_i =0$. 
Note that $\phi_i \in V_{\alpha,0}^{n_\alpha-(k+2i-j)} 
\oplus V_{\alpha,1}^{n_\alpha+1-(k+2i+1-j)}$.  
Then the Hard Lefschetz (\ref{eq:HL_Valphaj}) 
for $V_{\alpha,*}$ implies $\phi_i=0$ and so 
$\hat\phi = \phi$. 
By the Lefschetz decomposition,  
we have $\hat\phi = \phi$ for every $\phi\in V_{\alpha,j}$. 
Therefore $z\cH_-\cap \Flim = H_{\rm CR}^*(\cX)$ 
and $z \cH_- = H_{\rm CR}^*(\cX)\otimes 
\cO(\Proj^1\setminus\{0\})$. 

It is easy to show the characterization of $v_0$.  
When $v_0$ is replaced with $\lambda v_0$ for some 
$\lambda\in\C$, 
the flat metric $g$ is multiplied by $\lambda^2$.  
\end{proof} 

\begin{remark}
The limiting Hodge structure $\Flim$ depends on the 
choice of co-ordinates $q^1,\dots,q^r$ 
on $\ov{U/H^2(X,\Z)}$. 
Another co-ordinate system 
$\hat{q}^a := c^a q^a \exp(F_a(q))$ with $F_a(0)=0$ 
changes $\Flim$ by the multiplication 
by $\exp(\sum_{a}\log c^a \Log(M_a)/(2\pi\iu))$. 
Under the monodromy invariance 
(\ref{eq:monodromy-invariance_opposite}) for $\cH_-$, 
$\cH_-$ being opposite to $\Flim$ (\ref{eq:limiting_opposite}) 
is independent of the choice of a co-ordinate system  
since $\Log(M_a)$ preserves $\cH_-$. 
\end{remark} 

\begin{remark} 
We can normalize the dilaton shift $v_0\in z\cH_-/\cH_-$ 
using the integral structure $F_\Z$. 
The dilaton shift $v_0$ defines a primitive section $s_0$ 
of the quantum $D$-module via the  
trivialization (\ref{eq:triv_opposite}). 
Under the condition (\ref{eq:uniqueness_dilaton}), 
there exists a one-dimensional subspace $\C A_0$ 
of the space $\Sol(\cX)$ of flat sections 
which is invariant under every Galois action 
and contained in the image of 
$(\id - G^{\Sol}(\xi))^n$ for some 
unipotent operator $G^{\Sol}(\xi)$ 
with the maximum unipotency $n=\dim_\C\cX$. 
(This can be seen from the cohomology framing. 
See (\ref{eq:cohframing_Galois}).) 
An integral generator $A_0$ of this 
subspace is determined up to sign: 
In fact, this is given by the structure 
sheaf of a non-stacky point 
$A_0 = \pm \cZ_K(\cO_{\rm pt})$. 
The choice $v_0=\pm\unit$ corresponds to the normalization 
$(s_0, A_0)_F \sim (2\pi \iu)^n/(2\pi z)^{\frac{n}{2}}$ 
in the large radius limit. 
\end{remark} 

\subsection{Symplectic transformation  
between Givental spaces}
\label{subsec:symp}
Here we see that 
Picture \ref{Pic:GQDM} gives rise to 
a symplectic transformation $\U$ 
between the Givental spaces $\cH^{\cX_1}$ and $\cH^{\cX_2}$. 
The transformation $\U$ was introduced 
in \cite{CIT} to describe relationships 
between the genus zero Gromov-Witten theories  
of $\cX_1$ and $\cX_2$. 
As we have seen, the genus zero theory defines 
a semi-infinite variation of Hodge structures 
$\F_\tau^{\cX_i} \subset \cH^{\cX_i}$ in the Givental spaces. 
We shall see in (\ref{eq:U_matches}) that 
they match under $\U$: 
$\U \F_\tau^{\cX_1} = \F_\tau^{\cX_2}$. 
This implies that 
Givental's Lagrangian cones $\cL_i\subset \cH^{\cX_i}$ 
\cite{CG} 
swept by the semi-infinite subspaces $z \F_\tau^{\cX_i}$ 
are mapped to each other under $\U$: 
\[
\U \cL_1 = \cL_2, \quad \text{where} \quad 
\cL_i := \bigcup_{\tau} z \F_\tau^{\cX_i} \subset \cH^{\cX_i}.  
\]
The Lagrangian cone $\cL_i \subset \cH^{\cX_i}$ 
can be also described as the graph of the 
genus zero descendant potential of $\cX_i$ \cite{CG} 
and encodes all the information on genus zero 
Gromov-Witten theory. 
% The equality $\U \cL_1 = \cL_2$ expresses  
% the relationships of the genus zero theories 
% in a more intuitive way. 
In the literature \cite{CIT, Coates-R, Coates}, 
the crepant resolution conjecture was formulated 
in this way and verified in several examples. 
See these references for more 
details and examples of $\U$.  

Take a path $\gamma\colon [0,1]\to \cM$ 
connecting two cusp neighborhoods $V_1$, $V_2$. 
Then we have the analytic 
continuation map (\ref{eq:analyticcont})
$P_\gamma\colon \Sol(\cX_1) \to \Sol(\cX_2)$ 
along the path $\hat{\gamma}=(\gamma,1)\colon 
[0,1] \to \cM\times \C^*$. 
Through the cohomology framing $\cZ_{\rm coh}$ 
(\ref{eq:cohframing}), the $P_\gamma$ 
induces the following isomorphism: 
\begin{equation}
\label{eq:U_coh}
\U_{\rm coh}\colon H_{\rm CR}^*(\cX_1) 
\to H_{\rm CR}^*(\cX_2), \quad 
\U_{\rm coh}=\cZ_{\rm coh}^{-1} P_\gamma \cZ_{\rm coh}. 
\end{equation}
Recall that the Givental space $\cH^{\cX_i}$ 
is identified with the space of (multi-valued) 
sections of $F$ over $V_i\times \C^*$ 
which are flat in the $V_i$ direction. 
Therefore, the analytic continuation along $\hat\gamma$ 
also induces the map between the Givental spaces: 
\begin{equation}
\label{eq:U_Givental} 
\U \colon \cH^{\cX_1} \to \cH^{\cX_2}.  
\end{equation} 
The map $\U$ is an $\cO(\C^*)$-linear isomorphism 
preserving the pairing $(\cdot,\cdot)_{\cH}$ on 
the Givental spaces. In particular, $\U$ 
is a symplectic transformation with respect to 
the symplectic form (\ref{eq:symp_form_Giv}).
Recall that the cohomology framing 
identifies $\phi\in H_{\rm CR}^*(\cX_i)$ 
with a flat section $L(\tau,z)z^{-\mu_i}z^{\rho_i} \phi$ 
of $QDM(\cX_i)$. 
Also recall that $\phi(z)$ in the Givental space
$\cH^{\cX}$ corresponds to the flat section 
$L(\tau,z)\phi(z)$. 
Therefore, one has the commutative diagram 
involving ``multi-valued" Givental spaces: 
\begin{equation}
\label{eq:diag:U_Ucoh}
\begin{CD}
H_{\rm CR}^*(\cX_1) @>{\U_{\rm coh}}>> H_{\rm CR}^*(\cX)  \\
@V{z^{-{\mu_1}}z^{\rho_1}}VV @V{z^{-\mu_2}z^{\rho_2}}VV \\ 
\cH^{\cX_1}\otimes_{\cO(\C^*)} \cO(\widetilde{\C^*}) 
@>{\U}>>
\cH^{\cX_2}\otimes_{\cO(\C^*)} \cO(\widetilde{\C^*}) 
\end{CD}
\end{equation} 
where $\rho_i=c_1(\cX_i)$ and $\mu_i$ is the Hodge grading 
operator of $\cX_i$. 

For a rational number $f\in [0,1)$, we set  
\begin{equation}
\label{eq:f-part}
H^*_{\rm CR}(\cX)_f := 
\bigoplus_{\fract{\iota_v}=f} H^{*-2\iota_v}(\cX_v) 
=\footnote{This equality holds since 
we ignore cohomology classes of odd parity.} 
\bigoplus_{\fract{p/2}=f} H^p_{\rm CR}(\cX). 
\end{equation} 
Here $\fract{\iota_v}$ is the fractional part of $\iota_v$. 
Correspondingly, we set 
\[
\cH_f^{\cX} := H_{\rm CR}^*(\cX)_f \otimes 
\cO(\C^*) \subset \cH^\cX. 
\]
We list basic properties of $\U_{\rm coh}$ and $\U$, 
some of which already appeared in \cite{CIT, Coates-R}. 
We will use these later. 

\begin{lemma}
Under Picture \ref{Pic:GQDM}, 
the analytic continuation maps $\U_{\rm coh}$ and 
$\U$ given in (\ref{eq:U_coh}), (\ref{eq:U_Givental}) 
satisfy the following: 
\begin{gather} 
\label{eq:rho_comm} 
\U_{\rm coh} \rho_1 = \rho_2 \U_{\rm coh}, \quad   
\U \rho_1 = \rho_2 \U, \\
\label{eq:fpart_preserved}
\U_{\rm coh} H_{\rm CR}^*(\cX_1)_f = H_{\rm CR}^*(\cX_2)_f, 
\quad  
\U \cH^{\cX_1}_f = \cH^{\cX_2}_f,  \\
\label{eq:U_Ucoh}
\U = z^{-\mu_2} \U_{\rm coh} z^{\mu_1}, \\ 
\label{eq:U_matches} 
\U \F^{\cX_1}_\tau = \F^{\cX_2}_\tau, 
\quad \tau \in \cM.  
\end{gather} 
Here the $\F^{\cX_i}_\tau \subset \cH^{\cX_i}\cong \cH$ 
is the semi-infinite Hodge structure (\ref{eq:Ftau}) 
at $\tau\in \cM$ considered as a subspace of the 
Givental space. 
The equation (\ref{eq:U_Ucoh}) shows that 
$\U$ is degree-preserving, 
where the grading on $\cH^\cX$ is given by 
the usual grading on $H_{\rm CR}^*(\cX)$ and $\deg z=2$. 

Assume that $\cX_1$ and $\cX_2$ are $K$-equivalent 
and related by the diagrams (\ref{eq:Kequiv}), (\ref{eq:contract})
such that $\pi_1 \circ p_1 = \pi_2 \circ p_2$. 
Let $\gamma$ be the path in (ii) of Picture \ref{Pic:GQDM}. 
Then for a class 
$\alpha\in H^2(Z,\C)$, 
\begin{equation}
\label{eq:pullback_comm}
\U_{\rm coh} (\pi_1^*\alpha) = 
(\pi_2^*\alpha) \U_{\rm coh}, \quad 
\U (\pi_1^*\alpha) = (\pi_2^*\alpha) \U. 
\end{equation} 
\end{lemma} 
\begin{proof}
The analytic continuation along $\hat{\gamma} = (\gamma,1)$ 
must be equivariant under the monodromy in $z\in \C^*$. 
A simple calculation shows that the monodromy in $z$ 
acts on $\Sol(\cX_i)\cong H^*_{\rm CR}(\cX_i)$ by 
\begin{equation}
\label{eq:mon_z}
M_i =(-1)^n e^{-2\pi\iu\rho_i} \bigoplus_{v\in \sfT_i} 
e^{2\pi\iu \iota_v}, \quad n=\dim\cX_i,  
\end{equation} 
where $\sfT_i$ is the index set of 
the inertia component of $\cX_i$. 
Then $M_2 \U_{\rm coh} = \U_{\rm coh} M_1$. 
Taking a sufficiently high powers of $M_i$, 
we have $e^{-k_0 2\pi\iu \rho_2} \U_{\rm coh} 
= \U_{\rm coh} e^{-k_0 2\pi\iu \rho_1}$.  
This shows the first equation of (\ref{eq:rho_comm}). 
Therefore we also have
$\U_{\rm coh} \bigoplus_{v\in\sfT_1} e^{2\pi\iu \iota_v} 
= \bigoplus_{v\in \sfT_2} e^{2\pi\iu \iota_v} \U_{\rm coh}$. 
This shows the first equation of (\ref{eq:fpart_preserved}). 
Since $\U_{\rm coh}$ commutes with $\rho_i$, 
$z^{\rho_i}$'s in the commutative diagram (\ref{eq:diag:U_Ucoh}) 
cancel each other. This shows (\ref{eq:U_Ucoh}) 
and in turn shows the second equations of 
(\ref{eq:rho_comm}), (\ref{eq:fpart_preserved}). 
The equation (\ref{eq:U_matches}) is a tautological 
relation since $\F_\tau^{\cX_1}$ and $\F_\tau^{\cX_2}$ 
arise from the same subspace $\F_\tau$ of $\cH$. 

When $\cX_1$ and $\cX_2$ are related by the 
birational correspondences (\ref{eq:Kequiv}), 
(\ref{eq:contract}), the analytic continuation 
$P_\gamma$ is equivariant under the monodromy 
(Galois) action coming from a line bundle $L$ on $Z$. 
By the formula (\ref{eq:cohframing_Galois}) 
of the Galois action in terms of $\cZ_{\rm coh}$, 
we have $\U_{\rm coh} e^{-2\pi\iu \pi_1^*c_1(L)} = 
e^{-2\pi\iu \pi_2^* c_1(L)} \U_{\rm coh}$ 
and (\ref{eq:pullback_comm}) follows. 
\end{proof}

\subsection{Hard Lefschetz condition} 
\label{subsec:HL}
We have seen under Picture \ref{Pic:GQDM} 
that quantum cohomology of $\cX_1$ and $\cX_2$ 
underlies the same $F$-manifold $\cM$  
(Proposition \ref{prop:Fmanifold}) 
and that the $F$-manifold structure 
can be (canonically) lifted over $V_i$ 
to a Frobenius manifold structure 
by the opposite subspace $\cH_-^{\cX_i}$ 
(Propositions \ref{prop:Frobenius} and 
Theorem \ref{thm:uniqueness_cusp}). 
Since a Frobenius structure is well-defined 
over the complement of an analytic subvariety of $\tcM$, 
we can compare the two Frobenius structures 
arising from different cusps $V_1, V_2$. 
However, there are some examples where 
they do not necessarily coincide \cite{ABK,CIT}. 
The \emph{Hard Lefschetz condition} 
introduced in \cite{CIT, BG} is a criterion 
for the two Frobenius structures to match. 
The point is that the monodromy action 
coming from line bundles on $Z$ 
uniquely fixes opposite subspaces 
under this condition. 

In this section, we consider the case where 
$\cX_1$ and $\cX_2$ are $K$-equivalent (\ref{eq:Kequiv}) 
and related by the birational correspondence: 
\[
\begin{CD}
\cX_1 @>{\pi_1}>> Z @<{\pi_2}<< \cX_2 
\end{CD}
\]
such that $\pi_1\circ p_1 = \pi_2 \circ p_2$. 
\begin{definition}
Assume that $H_{\rm CR}^*(\cX_i)$ is 
graded by integers. 
We say that $\pi_i\colon \cX_i \to Z$ 
satisfies the \emph{Hard Lefschetz condition} if 
the map 
\[
(\pi_i^*\omega_Z)^k \colon 
H_{\rm CR}^{n-k}(\cX_i) \to H_{\rm CR}^{n+k}(\cX_i)  
\]
is an isomorphism for a class $\omega_Z$ of 
an ample line bundle on $Z$.  
\end{definition} 

\begin{remark}
In the context of crepant resolution conjecture, 
one can take $\cX_1=\cX$, $Z$ to be 
the coarse moduli space $X$ of $\cX$ and 
$\cX_2$ to be a crepant resolution $Y$ of $X$.  
The Hard Lefschetz condition was originally 
discussed in \cite{CIT,BG} for 
the natural map $\cX \to X$.  
As was observed in \cite{Fern}, 
the Hard Lefschetz condition 
for $\cX\to X$ is equivalent to 
\[
\iota_v = \iota_{\inv(v)} \quad \forall v\in \sfT. 
\]
This definition applies to the case where $\cX$ 
is non-compact. It is important to consider 
non-compact cases, but unfortunately,  
the discussion in this section 
does not apply to a non-compact $\cX$. 
\end{remark} 

\begin{remark}
\label{rem:semismall}
Cataldo-Migliorini \cite{Cat-Mig} showed 
that when $\cX_i=Y$ is a smooth projective variety, 
$\pi \colon Y\to Z$ satisfies the Hard Lefschetz 
condition if and only if $\pi$ is semismall. 
Here a proper morphism $\pi\colon Y\to Z$ is said to be 
\emph{semismall} if 
$\dim Z^k + 2k \le \dim Y$, where 
$Z^k = \{ z\in Z\;;\; \dim \pi^{-1}(z) = k\}$. 
\end{remark} 

We will consider a generalization of 
the Hard Lefschetz condition, where 
we do not assume the integer grading 
and also include the ``bicentric" case. 

\begin{definition}
\label{def:gen_HL} 
(i) We say that a pair $(V,\omega)$ of 
a $\Q$-graded complex vector space $V$ and  
a nilpotent endomorphism $\omega\in \End(V)$ of degree 2 
is \emph{bicentric HL} 
if there exists a rational number $n\in \Q$ 
and a graded decomposition $V = V_0 \oplus V_1$ such that 
$V^p= 0$ unless $p\in n+\Z$ and 
\[
\omega^k \colon V_j^{n+j-k} \to V_j^{n+j+k} 
\text{ is an isomorphism for $j=0,1$ and all $k\ge 0$.}
\]
We call the set $\{n,n+1\}$ the \emph{bicenter}. 
Note that this definition contains 
the ``mono-centric" case 
where $V_0$ or $V_1$ vanishes. 

(ii) 
We say that a proper morphism 
$\pi\colon \cX \to Z$ 
satisfies the \emph{generalized Hard Lefschetz condition} 
if for every rational number $f\in [0,1)$, 
the pair $(H_{\rm CR}^*(\cX)_f, \pi^*\omega_Z)$ 
is bicentric HL, where $H_{\rm CR}^*(\cX)_f$ 
is the graded subspace of $H_{\rm CR}^*(\cX)$  
defined in (\ref{eq:f-part}) 
and $\omega_Z$ is a class of an ample line bundle on $Z$. 
\end{definition} 

\begin{remark} 
\label{rem:genHL_coarse}
When $\pi$ is the natural map $\cX\to X$ 
to the coarse moduli,  
the generalized Hard Lefschetz condition for 
$\pi$ reads as follows: 
For every rational number $f\in [0,1)$, 
there exists $n_f\in \Q$ such that 
\[
\fract{\iota_v} = f \ 
\Longrightarrow \ 
\dim_\C \cX_v + 2\iota_v = n_f \text{ or } n_f+1.  
\]
Here $\{n_f,n_f+1\}$ is the bicenter of 
$(H_{\rm CR}^*(\cX)_f,\omega_X)$.    
\end{remark}

\begin{theorem}
\label{thm:HL} 
Let $\cX_1,\cX_2$ be $K$-equivalent  
smooth Deligne-Mumford stacks 
related by the diagrams (\ref{eq:Kequiv}), 
(\ref{eq:contract}) such that 
$p_1^*K_{\cX_1} = p_2^*K_{\cX_2}$ and 
$\pi_1\circ p_1 = \pi_2 \circ p_2$.  
Assume that $\pi_1 \colon \cX_1 \to Z$ satisfies the 
(generalized) Hard Lefschetz condition.  
Under Picture \ref{Pic:GQDM}, 
the standard opposite subspaces $\cH_-^{\cX_1}$, 
$\cH_-^{\cX_2}$ coincide under the analytic continuation  
along the path $\gamma$ in (ii) of Picture \ref{Pic:GQDM}, 
i.e. $\U(\cH_-^{\cX_1}) = \cH_-^{\cX_2}$. 
Moreover, 

(i) If $\cX_1$ or $\cX_2$ does not have generic stabilizers, 
the Frobenius manifold structures on $\cM$ 
coming from the quantum cohomology of 
$\cX_1$ and $\cX_2$ coincide up to 
a scalar multiple of the flat metric 
though the analytic continuation along $\gamma$. 

(ii) There is a graded isomorphism 
$(H_{\rm CR}^*(\cX_1),\pi_1^*\omega_Z) \cong 
(H_{\rm CR}^*(\cX_2),\pi_2^*\omega_Z)$ preserving 
the actions of $\omega_Z$. 
In particular, $\pi_2\colon \cX_2\to Z$ 
also satisfies the (generalized)
Hard Lefschetz condition. 
\end{theorem} 

This theorem is a generalization of a result in 
\cite{CIT}. 
We use the following lemma in the proof.  
\begin{lemma}
\label{lem:bicHL}
Let $V_{i}$, $i=1,2$ be $\Q$-graded 
vector spaces and $\omega_i\in \End(V_i)$ 
be nilpotent endomorphisms of degree two.  
Assume that $V_1$ and $V_2$ are isomorphic 
as graded vector spaces 
and that there exists a (not necessarily graded) 
linear isomorphism $\U\colon V_1\to V_2$ such that 
$\U \omega_1 = \omega_2 \U$.  
If $(V_1,\omega_1)$ is bicentric HL, 
then there exists a (not canonical) 
graded isomorphism $\varphi\colon V_1\to V_2$ 
such that $\varphi \omega_1 = \omega_2 \varphi$. 
In particular, $(V_2, \omega_2)$  
is also bicentric HL. 
\end{lemma} 
\begin{proof}
Let $V$ be a $\Q$-graded vector space 
and $\omega$ be a nilpotent operator on $V$ 
of degree 2. 
Let $a_1\ge a_2\ge \cdots\ge a_l$ 
be lengths of the Jordan cells appearing 
in the Jordan normal form of $\omega$. 
Then we can take a basis of $V$ 
of the form 
\begin{equation}
\label{eq:omegaseq}
\{\omega^k \phi_j\;;\; 1\le j\le l, \ 
0\le k\le a_j\}, \quad a_1\ge a_2 \ge \cdots \ge a_l  
\end{equation} 
such that $\omega^{a_j+1} \phi_j=0$. 
Here we can assume that $\phi_j$ is homogeneous. 
Set $\deg \phi_j=-a_j+\lambda_j$ 
for some $\lambda_j\in \Q$. 
By rearranging the basis, we can assume 
that $\lambda_j\ge \lambda_{j+1}$ if $a_j=a_{j+1}$. 
The sequence $\{(a_j,\lambda_j)\}_{j\ge 1}$ 
is uniquely determined by $(V,\omega)$ 
and we call it the \emph{type} of $(V,\omega)$. 
It suffices to show that $(V_i,\omega_i)$, $i=1,2$ 
have the same type. 
Let $\{(a_j^{(i)},\lambda_j^{(i)})\}_{j\ge 1}$ be 
the type of $(V_i,\omega_i)$.  
Since $\omega_1$ and $\omega_2$ are conjugate, 
we have $a_j := a_j^{(1)} = a_j^{(2)}$. 
%Let us assume $a_1=\cdots=a_k>a_{k+1}$. 
Because $(V_1,\omega_1)$ is bicentric HL, 
there exists $n\in \Q$ such that 
$\lambda_j^{(1)} = n$ or $n+1$ for all $j$. 
Then the degree spectrum of $V_1$ 
is contained in $[-a_1+n, a_1+n+1]$.  
Since $V_1$ and $V_2$ are isomorphic as 
graded vector spaces, we know that 
$[-a_j+\lambda^{(2)}_j, a_j+\lambda^{(2)}_j]\subset 
[-a_1+n,a_1+n+1]$. 
Therefore, $\lambda_j^{(2)} = n$ or $n+1$ if $a_j=a_1$. 
Take $k>0$ such that 
$a_1 = \cdots =a_k > a_{k+1}$. 
We calculate 
\begin{align*}
\dim V_1^{a_1+n+1} + \dim V_1^{-a_1+n}  
&= k \\ 
\dim V_2^{a_1+n+1} + \dim V_2^{-a_1+n} 
&= k + \sharp\{j>k\;;\; -a_j+\lambda_j^{(2)} = -a_1+n \} \\
&+ \sharp\{j>k\;;\; a_j+\lambda_j^{(2)} = a_1 + n+1\}.  
\end{align*}
Since these are equal, we have 
$[-a_j+\lambda^{(2)}_j, a_j+\lambda^{(2)}_j]\subset 
(-a_1+n,a_1+n+1)$ if $j>k$. 
Therefore, 
\begin{align*} 
\sharp\{j\le k\;;\; \lambda_j^{(1)}=n+1\} 
&= \dim V_1^{a_1+n+1} = \dim V_2^{a_1+n+1} \\ 
& = \sharp\{j\le k\;;\; \lambda_j^{(2)}=n+1\}.   
\end{align*} 
Hence $\lambda_j^{(1)} = \lambda_j^{(2)}$ 
for $j\le k$. This shows that $(V_1,\omega_1)$ 
and $(V_2,\omega_2)$ contains an isomorphic 
graded subspace $(V',\omega')$ of the type 
$\{(a_j,\lambda_j^{(1)})\}_{1\le j\le k}$. 
By taking the quotient by this subspace, 
one can proceed by the induction on dimensions. 
\end{proof} 

\begin{proof}[Proof of Theorem \ref{thm:HL}]
Take a path $\gamma\colon [0,1]\to \cM$ 
satisfying the condition (ii) of Picture \ref{Pic:GQDM}. 
The analytic continuation map $P_\gamma$ 
(\ref{eq:analyticcont}) along the path 
$\hat\gamma=(\gamma,1)$ 
induces maps $\U_{\rm coh}$ (\ref{eq:U_coh}) 
and $\U$ (\ref{eq:U_Givental}). 
Recall that $\U_{\rm coh}$ 
splits into isomorphisms 
$\U_{{\rm coh},f} \colon H^*_{\rm CR}(\cX_1)_f
\to H^*_{\rm CR}(\cX_2)_f$ 
for each $f\in [0,1)$ by (\ref{eq:fpart_preserved}).  
By (\ref{eq:pullback_comm}), we have 
\begin{equation}
\label{eq:Ucoh_omegaZ}
\U_{{\rm coh},f} (\pi_1^*\omega_Z) =  
(\pi_2^*\omega_Z) \U_{{\rm coh},f}. 
\end{equation} 
for an ample class $\omega_Z$ on $Z$. 
On the other hand, 
by the theorem of Lupercio-Poddar \cite{Lup-Pod} 
and Yasuda \cite{Yasuda1, Yasuda2}, 
$H^*_{\rm CR}(\cX_1)$ and  
$H^*_{\rm CR}(\cX_2)$ are isomorphic as 
graded vector spaces 
when $\cX_1$ and $\cX_2$ are $K$-equivalent.  
Thus $H^*_{\rm CR}(\cX_1)_f$ 
and $H^*_{\rm CR}(\cX_2)_f$ are also 
isomorphic as graded vector spaces. 
By Lemma \ref{lem:bicHL} and 
(\ref{eq:Ucoh_omegaZ}),  
we know that there is a graded 
isomorphism 
\[
\varphi\colon 
(H^*_{\rm CR}(\cX_1)_f, 
\pi_1^*\omega_Z) \to 
(H^*_{\rm CR}(\cX_2)_f, 
\pi_2^*\omega_Z) 
\]
and $(H^*_{\rm CR}(\cX_2)_f,\pi_2^*\omega_Z)$
is also bicentric HL. 

In general, a nilpotent operator $\omega$ 
on a vector space $V$ 
defines a unique (increasing) 
\emph{weight filtration} $W_i(V)$ of $V$ such that 
$\omega W_i(V) \subset W_{i-2}(V)$ and 
that $\omega^i \colon \Gr^W_i(V) \to \Gr^W_{-i}(V)$ 
is an isomorphism. 
Here $\Gr^W_i(V) = W_i(V)/W_{i-1}(V)$. 
% , one can check that 
% each filter $W_i(V)$ is homogeneous 
% with respect to the grading of $V$. 
When $V$ is a graded vector space, 
$\omega$ is of degree two and 
$(V,\omega)$ is bicentric HL with a 
graded decomposition $V=V_0\oplus V_1$ 
and a bicenter $\{n,n+1\}$ 
(as in Definition \ref{def:gen_HL}),   
the weight filtration of $V$ is given by 
\[
W_k(V) = V_0^{\ge n-k} \oplus V_1^{\ge n+1-k}. 
\]
Consider the case $(V,\omega) = 
(H^*_{\rm CR}(\cX_i)_f, \pi_i^*\omega_Z)$. 
Since the isomorphism $\U_{{\rm coh},f}$ 
preserves the weight filtration (by (\ref{eq:Ucoh_omegaZ}))
and $(H_{\rm CR}^*(\cX_i)_f,\pi_i^*\omega_Z)$ 
is bicentric HL, we have  
\begin{equation}
\label{eq:Ucoh_increase_deg}
\U_{{\rm coh},f}(H_{\rm CR}^p(\cX_1)_f) \subset 
H_{\rm CR}^{\ge p-1}(\cX_2)_f.  
\end{equation} 
When $\phi\in H^p_{\rm CR}(\cX_1)$, 
this together with the formula (\ref{eq:U_Ucoh})  
implies that $\U \phi$ cannot contain positive 
powers in $z$. 
Therefore a matrix representation $U(z)$ of $\U$ 
with respect to a basis of $H_{\rm CR}^*(\cX_i)$ 
does not contain positive powers in $z$. 
Since $\U$ preserves the pairing $(\cdot,\cdot)_\cH$, 
the same is true for the inverse $U(z)^{-1}$ 
which is the adjoint of $U(-z)$ with respect to 
the Poincar\'{e} pairing. Thus 
we have $\U \cH_-^{\cX_1} \subset \cH_-^{\cX_2}$ 
and $\U^{-1}\cH_-^{\cX_2} \subset \cH_-^{\cX_1}$. 
Hence $\U\cH_-^{\cX_1} = \cH_-^{\cX_2}$. 

Now we assume $\cX_i$ does not have  
generic stabilizers. Let $\cH_-\subset \cH$ 
be the common opposite subspace. 
Then the dilaton shift $v_0 \in z\cH_-/\cH_-$ 
is characterized up to a constant 
by the condition that $v_0$ is an eigenvector
of $\nabla_{z\partial_z}$ on $z\cH_-/\cH_-$
of the smallest eigenvalue. This shows (i).  
The rest of the statements 
follows from what we already showed. 
\end{proof} 

\begin{remark}
We used the theorem of Lupercio-Poddar and Yasuda 
\cite{Lup-Pod, Yasuda1, Yasuda2} in the proof. 
However, as \cite{CIT} did, we can deduce 
the graded isomorphism 
$H_{\rm CR}^*(\cX_1)\cong H_{\rm CR}^*(\cX_2)$ 
from Picture \ref{Pic:GQDM} and 
certain additional assumptions.  
For example, we can show this under the assumption 
that \emph{$\cH_-^{\cX_2}$ is opposite 
to the limiting Hodge structure $\Flim^{\cX_1}$ 
at the cusp of $V_1$}, \emph{i.e.} 
$\U(\Flim^{\cX_1}) \oplus \cH_-^{\cX_2} = \cH^{\cX_2}$. 
This assumption was conjectured to hold 
for a general crepant resolution  
$\cX_2 = Y \to X \leftarrow \cX_1$ in \cite{Coates-R}. 
Interestingly, under the generalized 
Hard Lefschetz condition, this assumption 
is a consequence of Picture \ref{Pic:GQDM}.   
\end{remark} 

By Theorem \ref{thm:HL} and Cataldo-Migliorini's theorem 
\cite{Cat-Mig} (see Remark \ref{rem:semismall}),  
Picture \ref{Pic:GQDM} has  
the following interesting consequences:  
\begin{itemize} 
\item Let $\cX$ be a Gorenstein orbifold and $Y\to X$ 
be a crepant resolution. 
Then $\cX$ satisfies the Hard Lefschetz 
condition if and only if $Y\to X$ is semismall. 

\item Let $X_1$ and $X_2$ be 
$K$-equivalent smooth projective varieties 
related by the diagrams (\ref{eq:Kequiv}), (\ref{eq:contract}) 
with $\pi_1\circ p_1 = \pi_2 \circ p_2$. 
Then $X_1 \to Z$ is semismall if and only if 
$X_2\to Z$ is semismall. 
\end{itemize} 
The author learned from Tom Coates that 
the first statement has been conjectured 
by Jim Bryan \cite{B:email}. 

% Every nilpotent operator $\omega$ on a vector space $V$ 
% defines a unique (increasing) \emph{weight filtration} 
% $W_i(V)$ of $V$ such that 
% $\omega W_i(V) \subset W_{i-2}(V)$ and 
% that $\omega^i \colon \Gr^W_i(V) \to \Gr^W_{-i}(V)$ 
% is an isomorphism. 
% Here $\Gr^W_i(V) = W_i(V)/W_{i-1}(V)$. 
% When $V$ is a graded vector space and $\omega$ 
% is of degree two, one can check that 
% each filter $W_i(V)$ is homogeneous 
% with respect to the grading of $V$. 
% By considering a homogeneous splitting 
% of $V=W_m(V)\supset W_{m-1}(V) 
% \supset \cdots \supset W_0(V)\supset W_{-1}(V)$, 
% one can obtain a homogeneous basis of $V$ of the form 

\subsection{Integral periods (Central charges)} 
\label{subsec:intperiod}
Up to now, we have not used the 
integral structure $F_\Z$ of the 
global quantum $D$-module. 
In this section, we will see that 
the integral structure defines 
an integral co-ordinate --- integral period --- 
on the global K\"{a}hler moduli space. 
This is called a \emph{central charge} (see
(\ref{eq:centralcharge})) in physics. 
For example, using this, 
we can give a ``reason" why 
the specialization value of quantum parameters 
should be a root of unity in the crepant resolution 
conjecture \cite{I:real}. 
% This was used in \cite{I:real} 
% to give a reason why the specialization value 
% of quantum parameters should be a root of unity 
% in the crepant resolution conjecture. 
% Hopefully, the integral structure 
% would also explain the  
% (quasi)-modularity of the Gromov-Witten
% potentials observed in physics. 
In this section, we restrict our attention to 
the case of crepant resolution $\cX_1=\cX\to X\leftarrow Y=\cX_2$. 
Also we assume that $Y$ and $\cX$ are Calabi-Yau. 
The case where $c_1(\cX)$ is semi-positive 
can be discussed in a similar way 
by using the \emph{conformal limit}\footnote{This is 
very close to Y. Ruan's quantum corrected cohomology 
ring of $Y$ which has the quantum correction only  
from the exceptional locus \cite{Ruan:crepant2}; 
In the abstract Hodge theory, this is also known 
as a graded quotient by the 
Sabbah filtration \cite{Sab:hypergeom,HS}.} 
introduced in \cite{I:real}. 
See \cite{I:real} for semi-positive case. 

Let $\cX$ be a Calabi-Yau Gorenstein orbifold 
of dimension $n$ and 
$\pi\colon Y\to X$ be a crepant resolution of 
the coarse moduli space $X$. 
Note that the Gorenstein assumption implies that 
$H_{\rm CR}^*(\cX)$ is graded by even integers. 
In Calabi-Yau case, the base space of the 
quantum $D$-module has a distinguished locus 
where the Euler vector field $E$ vanishes. 
By the formula (\ref{eq:Euler}), this 
is exactly the small (orbifold) quantum cohomology locus 
$H^2_{\rm CR}(\cX)$ or $H^2(Y)$. 
Recall that the Euler vector field 
is globally defined on $\cM$ by 
Section \ref{subsec:Fmanifold}. 

\begin{assumption}
The locus $\cM_0\subset \cM$ where the Euler vector 
field vanishes is connected. 
Also the path $\gamma\colon [0,1]\to \cM$ in 
(ii) of Picture \ref{Pic:GQDM} can be chosen 
so that it is contained in $\cM_0$. 
\end{assumption} 

In Calabi-Yau case ($\rho=0$), the situation 
is greatly simplified. 
The monodromy in $z\in \C^*$ is almost trivial 
and given by $(-1)^n$ by (\ref{eq:mon_z}). 
Over the locus $\cM_0$, the global quantum $D$-module 
gives rise to a \emph{finite dimensional 
variation of Hodge structures (VHS)}. 
The finite dimensional VHS arises from 
the filtration of flat sections 
by the pole/zero orders at $z=0$. 
The space $\Sol$ of multi-valued $\nabla$-flat sections 
of $F$ is single-valued in $w= z^{1/2}$ 
since the monodromy in $z$ is $\pm 1$. 
Moreover, over the locus $\cM_0$, 
the flat connection $\nabla$ has a logarithmic pole 
at $z=0$ since $\cU=\cA_{E}(\tau,0)$ in (\ref{eq:nabla_z_localframe}) 
is zero. Therefore, a $\nabla$-flat section 
$s(\tau,z)\in \Sol$ is at most meromorphic at $w=z^{1/2}=0$. 
This introduces the decreasing filtration 
$\Sol=F^0_\tau(\Sol)\supset F^1_\tau(\Sol) \supset \cdots \supset 
F^n_\tau(\Sol) \supset 0$ for $\tau\in \cM_0$:  
\[
F^p_\tau(\Sol)= \{s\in \Sol\;;\; 
z^{\frac{n}{2}-p}s(\tau,z) \text{ is regular at } z=0\}. 
\] 
Note that the factor 
$z^{\frac{n}{2}}$ cancels with 
the monodromy of $s(\tau,z)$ in $z$.   
On the neighborhoods $V_1, V_2$ of cusps, 
$\Sol$ is identified with $\Sol(\cX), \Sol(Y)$ 
and $F^p(\Sol)$ can be described as follows. 
Because $E=0$ on $\cM_0$, 
$\nabla_{z\partial_z}= z\partial_z + \mu$ 
for quantum $D$-modules 
and we have 
\begin{align}
\label{eq:Fp_nearcusp}
\begin{split}
F^p_\tau(\Sol) &\cong  
\{ s \in \Sol(\cX) \;;\; 
s(\tau,z) = z^{-\mu} \phi, \ \exists \phi\in 
H^{\le 2n-2p}_{\rm CR}(\cX)\} \\
& \cong  \{ s\in \Sol(Y) \;;\; 
s(\tau,z) = z^{-\mu} \phi, \ \exists \phi\in H^{\le 2n-2p}(Y)\} 
\end{split}
\end{align} 
on $V_1\cap \cM_0$ and $V_2\cap \cM_0$ respectively. 
The usual Griffiths transversality 
and Hodge-Riemann bilinear relation hold for 
$F^p_\tau(\Sol)$: 
\[
d F^p_\tau(\Sol) \subset F^{p-1}_\tau(\Sol)\otimes 
\Omega_{\cM_0}^1, \quad 
(F^p_\tau(\Sol), F^{n-p+1}_\tau(\Sol))_\Sol = 0.    
\] 
Here the pairing $(\cdot,\cdot)_\Sol$ is defined 
in the same way as in the case of quantum $D$-modules 
(see Definition \ref{def:spaceofflatsections}). 
The \seminf VHS $\F_\tau$ at $\tau\in \cM_0$ 
can be recovered from $F^p_\tau(\Sol)$ as follows: 
\[
\F_\tau = 
(z^{-\frac{n}{2}} F^n_\tau(\Sol) 
+ z^{-\frac{n}{2}+1} F^{n-1}_\tau(\Sol)  
+ \cdots 
+ z^{\frac{n}{2}} F^{0}_\tau(\Sol) )\otimes \cO(\C). 
\]

We introduce an integral period on $\cM_0$ 
corresponding to an element of $\Sol_\Z$, \emph{i.e.} 
a section of the integral local system $F_\Z$. 
This coincides with the central charge 
introduced in (\ref{eq:centralcharge}) 
for quantum $D$-modules. 
Recall that the analytic continuation map 
$\Sol(\cX)\cong \Sol \cong \Sol(Y)$ 
along the path $\hat\gamma$ in Picture \ref{Pic:GQDM} 
is equivariant under the Galois action of line bundles 
of the coarse moduli space $X$. 
% Under Picture \ref{Pic:GQDM} for $\cX$ and $Y$, 
% we have the identification 
% $\Sol(\cX)\cong \Sol\cong \Sol(Y)$ 
% via the analytic continuation along the path 
% $\hat\gamma$. 
% This identification is equivariant under the 
% Galois action of line bundles of the coarse moduli space $X$. 
Take an ample line bundle $L$ on $X$ and 
consider the corresponding Galois action $M=G^\Sol([L])$ 
on $\Sol$. 

\begin{lemma}
\label{lem:A0} 
{\rm (i)}  $F^n_\tau(\Sol)\subset \Sol$ 
is a one dimensional subspace 
for a generic $\tau\in \cM_0$. 

{\rm (ii)} There exists a unique (up to sign) 
integral vector $A_0\in \Sol_\Z$ 
contained in the image of $(\Log(M)-1)^n$.  
Under the $K$-group framing  
(\ref{eq:Kgroupframing}) 
$\cZ_K \colon K(\cX) \to \Sol(\cX)$ 
(or $K(Y) \to \Sol(Y)$), $A_0$ is identified 
with the structure sheaf of a non-stacky point 
$A_0 = \pm \cZ_K(\cO_{\rm pt})$. 
\end{lemma}
\begin{proof}
Since $\dim F^n_\tau$ is upper semi-continuous, 
(i) follows from the description (\ref{eq:Fp_nearcusp}) 
of $F^n_\tau(\Sol)$ near the cusps. 
The operator $M$ corresponds to the unipotent 
operator $e^{-2\pi\iu c_1(L)}$ on $H^*_{\rm CR}(\cX)$
through the cohomology framing (\ref{eq:cohframing}), 
thus $\Image (\Log(M)-1)^n \cong \Image c_1(L)^n = 
H^{2n}(\cX)$ is one-dimensional. 
This contains an integral vector $\cZ_K(\cO_{\rm pt})$.  
\end{proof} 

% Take a double cover $\C^*_w\to \C^*=\C^*_z$ 
% of the $z$-plane, where $w=z^{1/2}$. 
% Note that there is a natural map 
% \[
% \jmath_\tau \colon \Gamma(\{\tau\}\times \C^*_w, F) \cong \Sol 
% \]
% which extends a section over $\{\tau\}\times \C^*_w$ 
% to a flat section.  
By Lemma \ref{lem:A0}, the following definition makes sense. 

\begin{definition}
\label{def:int_period} 
Let $\C^*_w\to \C^*=\C^*_z$ be the double 
cover of the $z$-plane with a co-ordinate $w=z^{1/2}$. 
Take a flat section $A_0\in \Sol_\Z$ in Lemma \ref{lem:A0}. 
A \emph{normalized primitive section} 
is a section 
$\tilde{s}_0\in \Gamma(\cM_0\times \C^*_w,F)$ 
satisfying  
\begin{itemize}
\item For every $\tau\in \cM_0$,  
$\tilde{s}_0(\tau,z)$ is the restriction of 
an element of $F^n_\tau(\Sol)$ to $\{\tau\}\times \C^*_w$. 
\item $(\tilde{s}_0(\tau,e^{\pi\iu} z),A_0(\tau,z))_F =1$.  
\end{itemize}
This $\tilde{s}_0$ is unique up to sign. 
An \emph{integral period} $\Pi_A$ 
associated to $A\in \Sol_\Z$ is 
the function on $\cM_0$ defined by 
\begin{equation}
\label{eq:int_period}
\Pi_A(\tau):=(\tilde{s}_0(\tau, e^{\pi\iu}z), A(\tau,z))_F, 
\quad \tau\in \cM_0. 
\end{equation} 
\end{definition} 

We compute the normalized primitive section 
and integral periods for 
the quantum $D$-modules of $\cX$ and $Y$.  
Using the fundamental solution $L(\tau,z)$ 
in Proposition \ref{prop:fundamentalsol}, 
we define the \emph{$J$-function} by 
\[
J(\tau,-z) := L(\tau,z)^\dagger \unit,  
\]
where $L(\tau,z)^\dagger$ is the adjoint with respect to 
the Poincar\'{e} pairing. 
The $J$-function has the following expression: 
\begin{align*}
&J(\tau,-z) = e^{-\tau_{0,2}/z} 
\Bigg (1- \frac{\tau'}{z} + \\
& + \sum_{\substack{
d\in \Eff_\cX, 1\le k\le N\\ d=0 \Rightarrow m\ge 2}} 
\corr{\tau',\dots,\tau', 
\frac{\phi_k}{z(z+\psi_{m+1})}}_{0,m+1,d} 
\frac{e^{\pair{\tau_{0,2}}{d}}}{m!} \phi^k \Bigg).  
\end{align*}
Here $\tau = \tau_{0,2} + \tau'$ is 
the decomposition in (\ref{eq:decomp_tau}).  
(This can be derived from (\ref{eq:fundamentalsol}) 
and the String equation.) 
When $\cX$ is Calabi-Yau and $\tau\in H^2_{\rm CR}(\cX)$, 
the $J$-function is homogeneous of degree zero 
and is of the form 
\begin{equation}
\label{eq:J_CY}
J(\tau,-z) = 1 - \frac{\tau}{z} + \sum_{k\ge 2} 
\frac{\alpha_k(\tau)}{z^k}, \quad 
\alpha_k(\tau)\in H^{2k}_{\rm CR}(\cX).
\end{equation} 
\begin{proposition} 
(In this proposition, $\cX$ can be $Y$.)  
The normalized primitive section 
of the quantum $D$-module is given by 
\[
\tilde{s}_0(\tau,z) = 
\frac{(2\pi z)^{\frac{n}{2}}}{(-2\pi)^n} \unit. 
\]
Therefore, the integral period $\Pi_A$ (\ref{eq:int_period}) 
associated to an integral flat section 
$A=\cZ_K(V)$, $V\in K(\cX)$ 
equals the central charge $Z(V)$ (\ref{eq:centralcharge}).  
This is a component of the $J$-function: 
\[ 
\Pi_A = Z(V) = (2\pi)^{-\frac{n}{2}} \iu^{-n}  
(J(\tau,-1), \Psi(V))_{\rm orb},  \quad \tau \in H^2_{\rm CR}(\cX),   
\]
where $\Psi(V)$ was defined in (\ref{eq:Kgroupframing}) 
and $J(\tau,-z)$ is the $J$-function.  
\end{proposition} 
\begin{proof}
By (\ref{eq:Fp_nearcusp}), $\tilde{s}_0$ satisfies 
the first condition in Definition \ref{def:int_period}. 
From $A_0= \cZ_K(\cO_{\rm pt}) = L(\tau,z) 
((2\pi\iu)^n/(2\pi z)^{\frac{n}{2}}) [{\rm pt}]$ 
and the formula (\ref{eq:J_CY}) for the $J$-function, 
the second condition follows. 
The rest of the statements just follows from the 
definition (\ref{eq:int_period}) of $\Pi_A$ 
with the formulas 
(\ref{eq:Kgroupframing}), (\ref{eq:centralcharge}), 
(\ref{eq:J_CY})
and $\mu^\dagger = -\mu$. 
\end{proof} 

\begin{remark}
The above calculation shows that the ``normalized" 
primitive section is (up to a function in $z$) 
nothing but the primitive section $s_0=\unit$ 
associated to the standard 
opposite subspace and dilaton shift 
(see Section \ref{subsec:opp_Frob}). 
The existence of a canonical (normalized) primitive section 
along the locus $\cM_0$ does not mean 
that the Frobenius manifold structures of $\cX$ and $Y$ 
are the same. 
In fact, the primitive sections $s_0$ 
of $\cX$ and $Y$ may differ 
outside the locus $\cM_0\subset \cM$. 
\end{remark} 

\begin{corollary} 
\label{cor:match_period} 
Under the Picture \ref{Pic:GQDM} 
and Conjecture \ref{conj:Ktheoryisom}, 
the central charges of the corresponding 
$K$-group elements define  
the same function (up to sign) on $\cM_0$: 
\[
Z^{Y}(V) = \pm Z^{\cX}(\U_K^{-1}(V)), \quad V\in K(Y),  
\]
where $Z^{\cX}$ and $Z^Y$ are the central 
charges (\ref{eq:centralcharge}) of $\cX$ and $Y$ 
respectively and $\U_K=\U_{K,\gamma}\colon K(\cX) \cong K(Y)$ 
is the isomorphism 
in Conjecture \ref{conj:Ktheoryisom}. 
The sign $\pm$ depends on the sign of 
$\U_K(\cO_{\rm pt}) = \pm \cO_{\rm pt}$ 
(conjecturally plus). 
\end{corollary} 

It is interesting to study what 
integral periods are \emph{affine linear}  
functions on $H^2_{\rm CR}(\cX)$ or $H^2(Y)$.  
For example, there exists an affine 
co-ordinate system on 
$H^2(\cX)\oplus \bigoplus_{\codim \cX_v = 2} H^0(\cX_v)
\subset H^2_{\rm CR}(\cX)$ or on $H^2(Y)$  
consisting of integral periods \cite[Proposition 6.3]{I:real}.  
If we have a stratum $\cX_v$ of 
codimension $\ge 3$ with $\iota_v=1$, 
the corresponding linear projection  
$H^2_{\rm CR}(\cX) \to H^0(\cX_v)=\C$ 
may not be written as an affine 
linear combination of integral periods. 
Also, an affine linear integral period 
on $H^2(Y)$ 
may not correspond to an affine linear 
integral period on $H^2_{\rm CR}(\cX)$. 
In the next section, 
we will examine some local examples.  
% For an ample line bundle $L$ on $X$, 
% the Galois action $M=G^\Sol([L])$ on $\Sol(\cX)$ 
% is unipotent by (\ref{eq:cohframing_Galois}). 
% Thus the logarithm $\Log(M)$ is a nilpotent 
% operator which is represented by 
% $-2\pi\iu c_1(L)$ under the cohomology framing 
% (\ref{eq:cohframing}) 
% $\cZ_{\rm coh}\colon H^*_{\rm CR}(\cX)\cong \Sol(\cX)$.  
% Thus this defines a canonical weight filtration 
% \[
% W_{-n}(\Sol)\subset W_{-n+1}(\Sol) 
% \subset \cdots \subset W_{n}(\Sol) = \Sol. 
% \]
% See the proof of Theorem \ref{thm:HL} for 
% the definition of weight filtration. 
% Since the monodromy logarithm $\Log(M)$ is 
% defined over $\Q$, the weight filtrations 
% are actually defined on $\Sol_\Q:=\Sol_\Z \otimes \Q$. 
% In terms of the cohomology framing 
% $H_{\rm CR}^*(\cX) \cong \Sol(\cX) \cong \Sol$, 
% the weight filtration is given by 
% \[
% W_{-k}(\Sol) = \bigoplus_{v\in \sfT} H^{\ge n_v+k}(\cX_v), \quad 
% n_v = \dim \cX_v. 
% \]
% Under Gorenstein assumption, 
% there is no codimension one stratum $\cX_v$ 
% and the codimension two stratum $\cX_v$ satisfies 
% $\iota_v=1$. 
% Therefore the fist few terms of the weight filtrations 
% are identified as 
% \[
% W_{-n}= W_{-n+1}= H^{2n}(\cX)  
% \subset 
% W_{-n+2} = H^{2n-2}(\cX) \oplus 
% \bigoplus_{\iota_v=1} H^{2n_v}(\cX_v). 
% \]
% Because the one-dimensional space 
% $W_{-n}$ is defined over $\Q$, it contains 
% an integral generator $A_0\in W_{-n}(\Sol)$ 
% unique up to sign. 

\subsection{Local examples}
\label{subsec:locex}
We consider the crepant resolution conjecture for 
$\cX=[\C^n/G]$ where $G\subset SL(n,\C)$  
is a finite subgroup and $n=2$ or $3$. 
A standard crepant resolution of $X=\C^n/G$
is given by the $G$-Hilbert scheme \cite{BKR}: 
\[
\pi\colon 
Y := G\operatorname{-Hilb}(\C^n) \to X=\C^n/G. 
\] 
Moreover, an equivalence of 
derived categories $D(Y) \cong D(\cX) := D^G(\C^n)$ 
is given by the Fourier-Mukai transformation 
$\Phi\colon D(Y) \to D(\cX)$ \cite{BKR}: 
\[
\begin{CD}
\Phi = \bR q_* \circ p^*, \quad  
Y @<{p}<< \cZ @>{q}>> \C^n. 
\end{CD}
\]
where $\cZ\subset Y\times \C^n$ 
is the universal subscheme and 
$p$ and $q$ are natural projections. 
It would be natural to conjecture that 
our $K$-group isomorphism $\U_K$ 
comes from this derived equivalence: 
\[
\U_K^{-1} \colon K_E(Y) \cong K_0^{G}(\C^n), \quad 
[V] \longmapsto [\bR q_* (p^*V)],    
\]
where $E=\pi^{-1}(0) \subset Y$ is the exceptional set. 
Recall that we need to use compactly supported $K$-groups  
in order to get well-defined central charges.  
For a rational curve $\Proj^1\cong C \subset E$ 
in the exceptional set, 
the central charge of the class 
$[\cO_C(-1)] \in K_{E}(Y)$ is given by (\emph{c.f.} 
Example \ref{ex:integral}) 
\[
Z^Y(\cO_C(-1)) = - \frac{1}{2\pi\iu} \tau \cap [C] 
\]
for $\tau\in H^2(Y)$. 
Let $\tau_C := \tau \cap [C]$, $\tau\in H^2(Y)$ 
be the co-ordinate on $H^2(Y)$ and  
$\varrho_C$ be the virtual representation of $G$ 
given by the Fourier-Mukai transform 
$[\varrho_C \otimes \cO_0] = [\bR q_* (p^*\cO_C(-1))]$. 
Corollary \ref{cor:match_period} gives 
the following conjecture: 

\begin{conjecture}  
\label{conj:localex}
The small quantum cohomology (or $D$-modules) of 
$\cX$ and $Y$ are isomorphic under the co-ordinate change 
\begin{equation}
\label{eq:conj_coordinatechange}
\tau_C = - 2\pi\iu Z^\cX(\cO_0\otimes \varrho_C) 
\end{equation} 
where the right-hand side is 
the central charge function on $H^2_{\rm CR}(\cX)$.  
See (\ref{eq:charge_C2}) and (\ref{eq:charge_C3}) 
for formulas of $Z^\cX(\cO_0\otimes \varrho_C)$. 
In particular, the quantum variable $q_C= \exp(\tau_C)$ 
specializes to $\exp(-2\pi\iu (\dim \varrho_C)/|G|)$ 
at the large radius limit point of $\cX$. 
\end{conjecture} 

\begin{remark}
(i) Because $\cX$ is not compact, 
the characterization of the vector $A_0$ in 
Lemma \ref{lem:A0} does not hold.  
However, we can expect that the conclusion of 
Corollary \ref{cor:match_period} still holds 
because the $K$-group class\footnote{This corresponds to 
$[\cO_0 \otimes \varrho_{\rm reg}]$ in $K_0^G(\C^n)$.} 
$[\cO_{\rm pt}]$ of a non-stacky point
should correspond to each other 
under a birational transformation. 

(ii) Since $H^2$-variables do not carry the degree, 
we expect that the co-ordinate change above 
is also correct for $\C^*$-equivariant quantum cohomology. 
Here $\C^*$ acts on $\C^n$ diagonally. 
In dimension two, the non-equivariant 
quantum product is constant in $\tau$, 
so only the equivariant version 
is interesting. 

(iii) The specialization of $q_C$ to a root of unity 
comes from the fact that the central charges  
(\ref{eq:charge_C2}), (\ref{eq:charge_C3}) 
of $[\cO_0 \otimes \varrho_C]=\U_K^{-1}[\cO_C(-1)]$ 
take rational values 
at the orbifold large radius limit point $\tau=0$.  
In \cite{I:real}, the rationality of the 
central charge of $\U_K^{-1}[\cO_C(-1)]$ 
at the large radius limit 
was also discussed without assuming the precise 
form of the $K$-group framing. 
When the coarse moduli space $X$ is projective, 
under the assumption that $H^*(\cX)$ 
is generated by $H^2(\cX)$ and 
the condition (\ref{eq:uniqueness_dilaton}),   
the rationality here is forced only 
by the monodromy consideration 
\cite{I:real}.  
\end{remark} 

We have two cases. 

\noindent
{\bf (Case 1)} \emph{When the Hard Lefschetz condition 
holds for $\cX\to X$}. 
Then we have \cite[Lemma 3.4.1]{B-Gh:poly} 
\begin{itemize}
\item $n=2$ \emph{or}
\item $n=3$ and $G$ is conjugate 
to a subgroup of $SL(2,\C)$ \emph{or}  
\item 
$n=3$ and $G$ is conjugate to a subgroup 
of $SO(3,\R)$. 
\end{itemize} 
In these cases, every inertia component 
has age $\iota_v=1$ and 
the small quantum cohomology is already ``big"   
(ignoring the unit direction), so 
the above conjecture determines 
the full relationships of quantum cohomology. 
Because all the central charges 
$Z^\cX(\cO_0\otimes \varrho)$ 
are affine linear on $H^2_{\rm CR}(\cX)$ 
(the third term in (\ref{eq:charge_C3}) does not exist),  
the co-ordinate change (\ref{eq:conj_coordinatechange}) 
preserves the flat structure on the base 
and the Frobenius structures match.  
Each irreducible component $C$ of the exceptional 
set $E$ is a rational curve and corresponds 
to a non-trivial irreducible representation 
$\varrho_C$ under the Fourier-Mukai transformation\footnote
{The author thanks Samuel Boissiere 
for explaining this for $G\subset SO(3,\R)$.} 
(see \cite{Ito-Nak, GNS,Boi-Sar}).  
The formula (\ref{eq:conj_coordinatechange}) 
agrees with the conjecture of Bryan-Gholampour 
\cite{B-Gh:ADE, B-Gh:poly,BG}. 
The conjecture has been proved 
for $A_n$ surface singularities 
$\cX=[\C^2/\Z_n]$ \cite{CCIT:comp} 
and for $\cX=[\C^3/\Z_2\times \Z_2]$ and  
$[\C^3/A_4]$ \cite{B-Gh:A4} 
(where $G=A_4$ is the alternating group; 
this is the only case where the non-abelian crepant 
resolution conjecture has been proved).

\vspace{5pt} 
\noindent
{\bf (Case 2)} \emph{When the Hard Lefschetz condition fails 
for $\cX\to X$}. 
This happens only when $n=3$. 
In this case, since we have the component with 
age $\ge 2$, the above conjecture does not 
give a full co-ordinate change between 
Frobenius manifolds (see 
Remark \ref{rem:fulldetermine} below).   
As we can see from (\ref{eq:charge_C3}),  
integral periods can be non-linear functions 
on $H^2_{\rm CR}(\cX)$, so the co-ordinate change 
(\ref{eq:conj_coordinatechange}) 
can be also non-linear. 
Consider the case $\cX=\C^3/\Z_3$, 
where $\Z_3$ acts on $\C^3$ by 
the weight $\frac{1}{3}(1,1,1)$. 
Then $Y$ is the total space of 
the canonical bundle of $\Proj^2$ 
with the exceptional set $E=\Proj^2$. 
The Fourier-Mukai transformation is given 
by the diagram 
\[
Y=\cO_{\Proj^2}(-3) \overset{p}{\longleftarrow} \cZ = 
\cO_{\Proj^2}(-1) \overset{q}{\longrightarrow} \C^3. 
\]
Let $\varrho_1, \varrho_2$ 
be the representations of $\Z_3$ 
such that $\varrho_k(1\mod 3) = e^{2\pi\iu k/3}$.  
For a degree one rational curve 
$\Proj^1\cong C \subset E$, 
the Fourier-Mukai transform of $\cO_C(-1)$ 
gives the representation 
$\varrho_C=2 \varrho_1 \oplus \varrho_2$. 
Thus the predicted co-ordinate change is 
\begin{equation}
\label{eq:coord_localP2}
\tau_C = -2\pi\iu - 
\frac{2\pi \sqrt{3}}{3\Gamma(\frac{2}{3})^3}\alpha^2 t 
+ \frac{2\pi \sqrt{3}}{\Gamma(\frac{1}{3})^3} \alpha 
\parfrac{F_0^{\cX}}{t}, 
\end{equation}
where $t$ is a co-ordinate 
on the twisted sector $H^2_{\rm CR}(\cX)$ 
dual to $\unit_{\frac{1}{3}}$, $\alpha=e^{2\pi\iu/3}$  
and $F_0^{\cX}$ is the genus zero potential of $\cX$ 
(see (\ref{eq:C3_pot})). 
Since we have \cite{CIT,CCIT:comp}: 
\[
F_0^{\cX}(t) = 
\frac{1}{3\cdot 3!} t^3
-\frac{1}{3^3\cdot 6!} t^6  
+ \frac{1}{3^2\cdot 9!} t^9 
- \frac{1093}{3^5\cdot 12!} t^{12}
+\cdots,   
\]
the co-ordinate change (\ref{eq:coord_localP2}) 
is quite non-linear. 
This (\ref{eq:coord_localP2}) 
agrees with the computation in \cite{CIT,Coates} 
up to the Galois actions $\tau_C\mapsto \tau_C +2\pi\iu,
t \mapsto \alpha^2 t$. 

\begin{remark}
\label{rem:fulldetermine} 
In the second case, we can predict   
the full relationships between the small quantum cohomology 
by considering the central charges of $[\cO_S]\in K_E(Y)$ 
associated to surfaces $S\subset E$ 
in Corollary \ref{cor:match_period}. 
Note that $Z^Y(\cO_S)$ contains 
the information of the derivative of the potential 
$F_0^Y$ (see Example \ref{ex:integral}, (ii)). 
The co-ordinate change of big quantum cohomology 
can be also determined by $\U_K$ in principle, 
but the formula could be very complicated.    
\end{remark} 

%%%%%%%%%%%%%%%%%%%%%%%%%%%%%%%%%
% References
%%%%%%%%%%%%%%%%%%%%%%%%%%%%%%%%%

\end{document}